%% file: main.tex
\documentclass[preprint,12pt]{elsarticle}

\usepackage{graphicx}
\usepackage{amsmath}
\usepackage{amssymb}
\usepackage{amsthm}

\usepackage{algorithm}
\usepackage{algpseudocode}
\algrenewcommand\textproc{}

\usepackage{xcolor}

\usepackage{subcaption}

\usepackage{ulem}

\usepackage{multicol}
\usepackage{multirow}
\usepackage{adjustbox}

\journal{Elsevier}

\input{preamble.pap.tex}

\input{preamble.macros.tex}

\input{preamble.tikz.tex}

\begin{document}

\begin{frontmatter}



\title{Model reduction of convection-dominated viscous conservation laws using implicit feature tracking and landmark image registration}


\author[inst1]{Victor Zucatti}
\author[inst1]{Matthew J. Zahr}

\affiliation[inst1]{organization={University of Notre Dame},
            city={Notre Dame},
            postcode={IN 46556},
            country={United States of America}}

\begin{abstract}
	
Reduced-order models (ROMs) remain generally unreliable for convection-dominated problems, such as those encountered in hypersonic flows, due to the slowly decaying Kolmogorov $n$-width of linear subspace approximations, known as the Kolmogorov barrier. This limitation hinders the accuracy of traditional ROMs and necessitates impractical amounts of training data during the offline phase. To address this challenge, we introduce a novel landmark-based registration procedure tailored for ROMs of convection-dominated problems. Our approach leverages limited training data and incorporates a nonlinear transformation of the data using a landmark-based registration technique combined with radial basis function (RBF) interpolation. During the offline phase, we align dominant convective features in a reference domain, resulting in a rapid decay of error relative to the reduced space dimension. Landmarks are generated through a three-step process: (1) detecting shocks via edge detection techniques, (2) sampling using Monte Carlo methods, and (3) domain partitioning with $k$-means clustering, where cluster centroids serve as landmarks. Accurate landmark correspondence is achieved by minimizing pairing distances for similar features. The online phase integrates standard minimum-residual ROM methodologies, extending the optimization space to include admissible domain mappings. We validate our approach on two test cases: a space-time Burgers’ equation parameterized by the initial condition, and a hypersonic viscous flow over a cylinder parameterized by the Mach number. Results demonstrate the efficacy of the proposed method in overcoming the Kolmogorov barrier and enhancing the reliability of ROMs for convection-dominated problems. 

\end{abstract}


\end{frontmatter}

\section{Introduction}
\label{sec:intro}

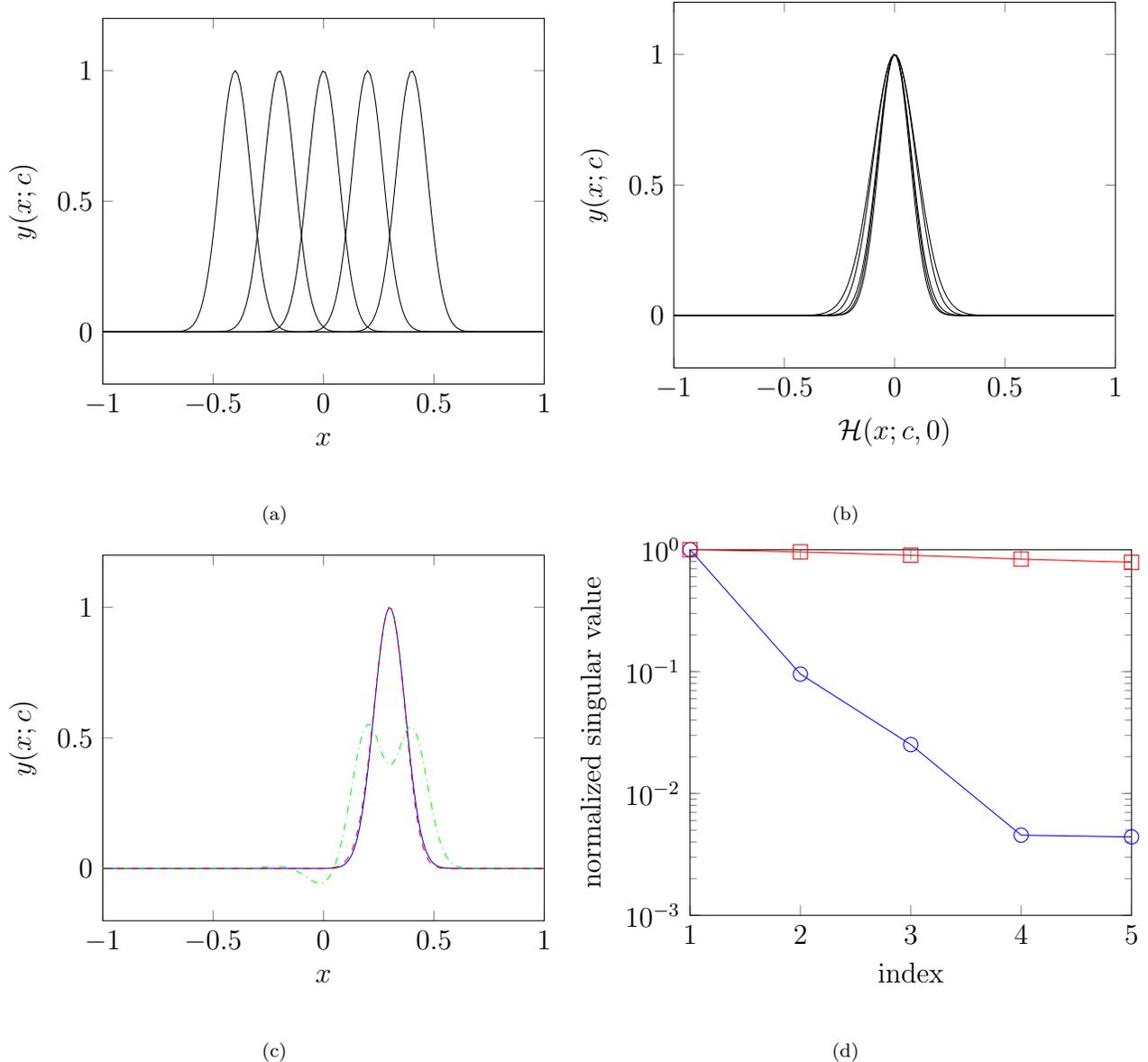
\begin{figure}[hbt!]
	\centering
	\begin{subfigure}{.49\textwidth}
		\centering
		\input{figs/methodology/bump1d_motivation/bump_1d_snapshots_xphys.tikz}
		\caption{\label{fig:bump_1d_snap}}
	\end{subfigure}
	\hfill
	\begin{subfigure}{.49\textwidth}
		\centering
		\input{figs/methodology/bump1d_motivation/bump_1d_snapshots_xref.tikz}
		\caption{\label{fig:bump_1d_snap_align}}
	\end{subfigure}
	\\
	\begin{subfigure}{0.49\textwidth}
		\centering
		\input{figs/methodology/bump1d_motivation/bump_1d_rom.tikz}
		\caption{\label{fig:bump_1d_snap_rom}}
	\end{subfigure}
	\hfill
	\begin{subfigure}{0.49\textwidth}
		\centering
		\input{figs/methodology/bump1d_motivation/bump_1d_sig_val.tikz}
		\caption{\label{fig:bump_1d_sig_val}}
	\end{subfigure}
	\caption{
		Top left: five different instances of the Gaussian bump function Eq. \eqref{eq:bump_1d} with equally spaced centers \( c \in [-0.4, 0.4] \).
		Top right: the different instances are aligned at \( x = 0 \) using the parametric mapping given by Eq. \eqref{eq:par_map_bump_1d}.
		Bottom left: we approximate a Gaussian bump centered at \( c = 0.3 \) (\ref{line:bump_1d_rom_ref}) by optimally combining the unaligned samples (\ref{line:bump_1d_rom_r}), as typically done in reduced-order modeling, and by composing the bump centered at \( c = 0 \) with the mapping \( x_t = 0.3 \) (\ref{line:bump_1d_rom_rft}).
		Bottom right: decay of the normalized singular values for aligned (\ref{line:bump_1d_sig_val_rft}) and unaligned (\ref{line:bump_1d_sig_val_r}) bumps, with the aligned case exhibiting significantly faster decay.
	}
	\label{fig:bump_1d}
\end{figure}

The extensive data from high-dimensional models (HDMs) can be used to construct a reduced-order model (ROM) through a common two-phase process, involving offline and online steps \cite{Rowley2004,Kevin01,lui_wolf_2019}. In the offline phase, a reduced-dimensional representation is extracted from HDM training data to create a simplified model, a pre-computation step that, while performed only once, can be costly due to the reliance on high-dimensional data. Conversely, the online stage involves solving the resulting system of equations at a reduced computational cost.
Despite significant research efforts, reduced-order models (ROMs) still face numerous challenges that make them generally unreliable. This is especially true when dealing with convection-dominated problems, such as those encountered in hypersonic flow. The primary issue with generalizing ROMs for these problems is the slowly decaying Kolmogorov $n$-width of linear subspace approximations, often referred to as the Kolmogorov barrier \cite{2016_Ohlberger}. This slow decay in error relative to the dimension of the reduced space limits the practical accuracy of ROMs and demands a large amount of training data, which is often impractical to gather in the offline phase \cite{Charbel_2024}. 
We reproduce this issue in Fig. \ref{fig:bump_1d}. Let $y(x;c)$ be the following parametrized Gaussian function defined as
\begin{equation}
	y(x; c) = \exp{(-100(x - c)^2)}
	\mbox{ ,}
	\label{eq:bump_1d}
\end{equation}
where $c$ defines the center.
We align the different instances at $x=0$, as illustrated in Fig. \ref{fig:bump_1d_snap_align}, using the following parametric mapping:
\begin{equation}
	\mathcal{H}(x;c,x_t) = 
	x_t \frac{(x-1)(x+1)}{(c -1)(c + 1)} +
	\frac{(x+1)(x-c)}{2(1 - c)} -
	\frac{(x-1)(x-c)}{2(c+1)}
	\mbox{ ,}
	\label{eq:par_map_bump_1d}
\end{equation}
where $x_t$ denotes the position at which the centers are aligned, with $x_t=0$ in this case.

The Kolmogorov barrier can be circumvented, for instance, by exploiting the local low-rank structure of convection-dominated problems.
In \cite{Amsallem2012_local,Amsallem2015_local}, this is achieved by dividing the training data in subregions with a local reduced-order basis constructed and assigned to each subregion.
Also taking advantage of localized low-rank structure, adaptive reduced order models (AROMs) \cite{2020_Peherstorfer_AADEIM,HUANG_2023,zucatti_arom} combine HDM and ROM operations on the fly to generate predictive ROMs. While AROMs may not achieve the high speedup factors typically associated with the two-step ROM approach \cite{lui_wolf_2019}, numerical experiments demonstrate that AROMs can effectively accelerate numerical solutions for problems where traditional order reduction techniques completely fail.

Nonlinear order reduction methods can also overcome the Kolmogorov barrier.
Quadratic manifolds have been employed in both physics-based approaches \cite{2022_Barnett} and data-driven techniques \cite{2022_Geelen_quad} for model order reduction.
A nonlinear manifold can also be generated through artificial neural networks (ANNs) and subsequently integrated into a reduced-order framework, resulting in models that exceed the performance of linear alternatives \cite{kevin_04_autoencoder,Romor2023,DIAZ_2024,KIM_2022,HOANG_2022}.
Alternatively, a nonlinear transformation can be applied directly to the data to eliminate the convection-dominated nature of the solutions \cite{2020_Tommaso,2023_Mirhoseini,nair_2019,mojgani_2017, Welper_2020,Nonino_2023,Nonino2024}.
This approach has been investigated through various methods, including leveraging concepts from optimal transport theory. Notable examples include framing the alignment as a Monge–Kantorovich mass transfer problem \cite{Iollo_2014_opt} or solving the Monge–Ampère equation \cite{Heyningen_2024}.
In \cite{2023_Mirhoseini,2024_Mirhoseini}, the mappings are derived through $r$-adaptivity with residual minimization, utilizing principles from the high-order implicit shock tracking (HOIST) method \cite{ZAHR_2018_hoist,HUANG_2022_hoist}.
Image registration is a methodology widely utilized for generating mappings that align snapshots \cite{Tommaso_2022, Tommaso_2023, Tommaso_2024, Ohlberger_2013,RAZAVI2025}.
The process starts by collecting multiple snapshots, selecting one as the template or reference, and then identifying mappings that optimally align the snapshots according to a chosen metric, such as mean squared error (MSE).

We introduce a landmark-based registration procedure for ROMs of convection-dominated problems that utilizes limited training data. Our nonlinear approach avoids the issues commonly encountered by traditional ROMs by first employing landmark-based registration, coupled with radial basis function (RBF) interpolation \cite{DEBOER_2007}, to align the dominant convective features during the offline phase in a reference domain. The transformation from convective to stationary solutions in the reference domain results in a rapid decay of error relative to the dimension of the reduced space. The landmarks are generated through a three-step process. First, we identify shocks using edge detection techniques (e.g., shock sensors, artificial viscosity). Next, samples are generated using Monte Carlo sampling methods, such as rejection or importance sampling \cite{casella_generalized_accept_reject}. The final step involves partitioning the domain using $k$-means clustering \cite{bishop_2006}, where the landmarks consist of the centroids of the clusters. Another important element of this registration method is determining landmark correspondence. For similar features and small displacements, pairing cluster centers by minimizing the total pairing distance results in accurate outcomes. Our online phase resembles standard minimum-residual reduced-order models; however, rather than solely minimizing the residual over the affine subspace of partial differential equation (PDE) states, our method enhances the optimization space to also incorporate admissible domain mappings similarly to \cite{2023_Mirhoseini,2024_Mirhoseini}.
We demonstrate the viability of our procedure on two viscous flow problems with shocks.
This is one of the first registration-based ROMs to successfully predict viscous hypersonic flow with limited training.
A novel residual-weighting approach proved critical to accurately predict flow fields containing both shock waves and boundary layers.

The remainder of this paper is structured as follows.
In Section \ref{sec:hdm}, we first introduce a general system of conservation laws, followed by its transformed counterpart on a fixed reference domain. We then present the high-dimensional model used to solve these equations.
Section \ref{sec:sparse_reg} introduces our shock alignment approach, beginning with landmark learning and the challenges of placing points on a continuous feature, followed by the formulation of the correspondence problem and strategies for matching points across different figures, and concluding with a radial basis function mesh warping method.
Next, Section \ref{sec:rom} introduces the reduced-order model with implicit feature alignment.
In Section \ref{sec:num_exp}, we compare the performance of our method with the standard fixed-mesh reduced-order model approach through two viscous flow problems. The first is a space-time Burgers' problem parameterized by the initial condition, while the second involves hypersonic viscous flow over a cylinder, parameterized by the Mach number. The presence of a boundary layer in the latter makes it particularly challenging.
Finally, \ref{sec:conclusions} presents the main conclusions and explores potential directions for future research.

\section{High-dimensional model}
\label{sec:hdm}

In this section, we first introduce the parameterized system of conservation laws (Section \ref{sec:claw}), then recast this system on a fixed reference domain (Section \ref{sec:claw_trans}), and finally, we discuss the discretization of the transformed equations (Section \ref{sec:claw_trans_disc}).

\subsection{Parametrized system of conservation laws}
\label{sec:claw}

A $\mu$-parameterized system of $m$ conservation laws, defined in a spatial domain $\Omega \subset \mathbb{R}^d$, can be expressed as
\begin{equation}
    \nabla \cdot f (q, \nabla q; \mu) = s(q, \nabla q; \mu)
    \mbox{ ,}
    \label{eq:claw}
\end{equation}
where $f : \mathbb{R}^m \times \mathbb{R}^{m \times d}  \times \mathcal{D} \rightarrow \mathbb{R}^{m \times d}$ is the flux function, $s : \mathbb{R}^m \times \mathbb{R}^{m \times d} \times \mathcal{D} \rightarrow \mathbb{R}^m$ is the source term, $q(x; \mu) \in \mathbb{R}^m$ is the vector of conservative variables implicitly defined as the solution of Eq. \eqref{eq:claw} at $x \in \Omega$, $\nabla \coloneqq [\partial_{x_1}, \ldots, \partial_{x_d}]$ is the gradient operator on the domain $\Omega \subset \mathbb{R}^d$, and $\mu \in \mathcal{D}$ is the parameter and $\mathcal{D}$ is the parameter domain. The domain boundary $\partial \Omega$ has an outward unit normal $n : \partial \Omega \rightarrow \mathbb{R}^d$.
This formulation is broad enough to encapsulate both steady second-order partial differential equations (PDEs) defined in a $d$-dimensional spatial domain and time-dependent PDEs defined in a $(d-1)$-dimensional spatial domain, which corresponds to a $d$-dimensional space-time domain. In general, we assume the flux function can be split into a viscous and inviscid term as 
\begin{equation}
	f (q, \nabla q; \mu) = 
	f^i (q; \mu) - f^v (q, \nabla q; \mu)
	\mbox{,}
	\label{eq:flux_decomp}
\end{equation}
where $f^i : \mathbb{R}^m \times \mathbb{R}^{m \times d}  \times \mathcal{D} \rightarrow \mathbb{R}^{m \times d}$ is the inviscid flux function and $f^v : \mathbb{R}^m \times \mathbb{R}^{m \times d}  \times \mathcal{D} \rightarrow \mathbb{R}^{m \times d}$ is the viscous flux function.

\subsection{Transformed system of conservation laws on a fixed reference domain}
\label{sec:claw_trans}

The proposed reduced-order model relies on mesh warping and employs the approach from \cite{PERSSON_2009}, where the governing equation \eqref{eq:claw} is reformulated on a fixed reference domain $\Omega_0 \coloneqq \tilde{\mathcal{G}}^{-1} (\Omega) \subset \mathbb{R}^d$ using a smooth, invertible mapping $\tilde{\mathcal{G}} : \mathbb{R}^d \rightarrow \mathbb{R}^d$.
Let $\mathbb{G}$ be any collection of bijections from the reference domain $\Omega_0$ to the physical domain $\Omega$. Then, for any $\mathcal{G} \in \mathbb{G}$, \eqref{eq:claw} can be written as a PDE on the reference domain as
\begin{equation}
    \nabla_0 \cdot F (Q, \nabla_0 Q; \mathcal{G}, \mu) = S(Q, \nabla_0 Q; \mathcal{G}, \mu)
    \mbox{ ,}
    \label{eq:claw_trans}
\end{equation}
where $F : \mathbb{R}^m \times \mathbb{R}^{m \times d} \times \mathbb{G} \times \mathcal{D} \rightarrow \mathbb{R}^{m \times d}$ is the transformed flux function, $S : \mathbb{R}^m \times \mathbb{R}^{m \times d} \times \mathbb{G} \times \mathcal{D} \rightarrow \mathbb{R}^m$ is the transformed source term, $\nabla_0 \coloneqq [\partial_{X_1}, \ldots, \partial_{X_d}]$ is the gradient operator on the domain $\Omega_0 \subset \mathbb{R}^d$, and $Q(X; \mathcal{G}, \mu) \in \mathbb{R}^m$ is the solution in the reference domain at point $X \in \Omega_0$.
The boundary of the domain $\Omega_0$ has an outward unit normal $N \colon \Omega_0 \rightarrow \mathbb{R}^d$. 
The physical solutions and their gradients are connected to the reference solutions and their gradients as follows
\begin{subequations}
    \begin{equation}
        Q(X; \mathcal{G}, \mu) = q(\mathcal{G} (X); \mu)
        \mbox{ ,}
    \end{equation}
    \begin{equation}
        \nabla_0 Q (X; \mathcal{G}, \mu) = 
        \nabla q (\mathcal{G} (X); \mu) G_{\mathcal{G}} (X)
        \mbox{ ,}
    \end{equation}
\end{subequations}
where $G_{\mathcal{G}} \coloneqq \nabla_0 \mathcal{G}$ is the mapping Jacobian and $g_\mathcal{G} \coloneqq det (\nabla_0 \mathcal{G})$ is its determinant. 
The reference and physical flux functions and source terms are related as
\begin{subequations}
    \begin{equation}
        F (w, \nabla_0 w; \mathcal{G}, \mu) =
        g_{\mathcal{G}} f (w, \nabla_0 w \cdot G_{\mathcal{G}}^{-1}; \mu) G_{\mathcal{G}}^{-T}
        \mbox{ ,}
    \end{equation}
    \begin{equation}
        S (w, \nabla_0 w; \mathcal{G}, \mu) = 
        g_{\mathcal{G}} s (w, \nabla_0 w \cdot G_{\mathcal{G}}^{-1}; \mu)
        \mbox{ ,}
    \end{equation}
    \label{eq:claw_trans_rel}
\end{subequations}
where $w: \Omega_0 \rightarrow \mathbb{R}^m$ represents a $m$-valued function defined over the reference domain. Finally, the systems of conservations laws in \eqref{eq:claw} and \eqref{eq:claw_trans} are equivalent. 
In other words, if $Q (X; \mathcal{G}, \mu)$ is the solution of \eqref{eq:claw_trans} and satisfies \eqref{eq:claw_trans_rel}, then $q(x; \mu)$ is the solution to \eqref{eq:claw}, and the reserve is also true.

\subsection{Discretization of the transformed conservation law}
\label{sec:claw_trans_disc}

The governing equations in \eqref{eq:claw_trans} are discretized using an appropriate numerical method
to yield the following nonlinear system of algebraic equations
\begin{equation}
    R (Q_h ; x_h , \mu) = 0
    \mbox{ ,}
    \label{eq:hdm}
\end{equation}
where $R : \mathbb{R}^N \times \mathbb{R}^{N_x} \times \mathcal{D} \rightarrow \mathbb{R}^N$ is the nonlinear residual function arising from the spatial discretization of \eqref{eq:claw_trans}, $x_h \in \Rbb^{N_x}$ are the degrees of freedom defining the domain mapping $\Gcal$, and $Q_h(X_h, \mu)\in\Rbb^{N}$  is the solution to the discretized partial differential equation (PDE)
implicitly defined as the solution of \eqref{eq:hdm}.
In this work, we use Hartmann's interior penalty discontinuous Galerkin (HIPDG) \cite{HARTMANN2008} to define the discrete state $Q_h$ and residual function $R$.
Furthermore, let $X_h\coloneqq((X_h)_1,\dots,(X_h)_{\hat{N}_x})\in\Rbb^{N_x}$ denote the nodes of the reference mesh on which the HIPDG discretization is defined, and let $x_h\coloneqq((x_h)_1,\dots,(x_h)_n)\in\Rbb^{N_x}$ denote the position of the corresponding nodes in the physical domain, i.e., $(x_h)_I = \Gcal((X_h)_I)$ for $I=1,\dots,\hat{N}_x$.
Then, the mapping defined from the discrete degrees of freedom is $\Gcal(X) = \psi_I(X) (x_h)_I$ for $I=1,\dots,\hat{N}_x$, where $\{\psi_I\}$ is a nodal, piecewise polynomial basis over the mesh.

\section{Shock alignment by sparse registration}
\label{sec:sparse_reg}

In this section, we begin by discussing the extraction of landmarks from edge-like features, such as shocks (Section \ref{sec:landmark}), and the process of matching landmarks between similar, yet distinct, shapes (Section \ref{sec:correspondence}). We then proceed to describe how the mesh is deformed through the displacement of these landmarks (Section \ref{sec:mesh_warp}).

\subsection{Landmark learning}
\label{sec:landmark}

Sparse image registration is used in image processing and computer vision to align two images by matching and comparing specific key points or features rather than using the entire image data. This method focuses on sparse sets of interest points instead of dense pixel-by-pixel analysis \cite{THIRION1998,Beg2005}, making it computationally efficient and suitable for scenarios where only key features matter.
In this approach \cite{Lowe_1999,Lowe2004,SuperPoint2018,Mikolajczyk2005,Schmid2000}, distinctive points in the image, such as corners, edges, or blobs, are identified. Afterwards, correspondences between key points in the two images are established using similarity metrics. Finally, the matching key points are used to estimate a transformation that aligns one image with the other.

In this work, our primary focus is on aligning edges, such as shocks. Unfortunately, placing point-based landmarks on an edge is a non-trivial task given its continuous nature. This is the reason many computer vision methods typically rely on detecting and matching features such as corners which are more distinctive and informative.
Moreover, the problems we are currently addressing are visually too simple to yield sufficient landmark points, if any, using traditional methods.
Taking all of this into account, we propose a method for edge landmark extraction, which involves the following steps:
\begin{enumerate}
    \item Apply an edge detection method to identify regions with steep gradients.
    \item Use rejection sampling \cite{casella_generalized_accept_reject} to generate sample points along the detected edges.
    \item Cluster the samples into $k$ groups using the $k$-means algorithm \cite{bishop_2006}, with the cluster centroids serving as landmarks.
\end{enumerate}

The first step identifies edges, regions of significant intensity change that indicate transitions between objects or regions, typically characterized by steep gradients and detected using first-order \cite{canny_1986} or second-order \cite{Marr_1980} derivative methods.
For instance, we can use the following gradient-based edge sensor $\hat{s} : \mathbb{R}^d\times\Dcal \rightarrow \mathbb{R}$ to identify shocks
\begin{equation}
	\hat{s} (x; \mu) = \| \nabla q (x; \mu) \|_{F}^2
	\mbox{ ,}
	\label{eq:shk_sensor_grad}
\end{equation}
where $\| \cdot \|_F$ is the Frobenius norm.
More recent advancements leverage machine learning techniques, including random forests \cite{Piotr_2013} and deep convolutional neural networks \cite{liu_2019}, to enhance edge detection \cite{Bertasius_2015, xie2017holistically, Zhuo_2021}.
Alternatively, incorporating a state-based artificial viscosity model \cite{barter_2010} into the system of conservation laws \eqref{eq:claw} not only enhances stability but also functions effectively as an edge detector, as demonstrated through the numerical experiments presented in Section \ref{sec:num_exp_cy}.

In the next step, we employ rejection sampling \cite{casella_generalized_accept_reject} to generate random samples in regions of $\Omega$ where $\hat{s}(x;\mu)$ is concentrated.
This method relies on a proposal distribution $v(x; \mu)$, which is straightforward to sample from, and an acceptance criterion based on the ratio of $\hat{s}(x; \mu)$ to $v(x; \mu)$.
In this work, we take
\begin{equation}
 v(x; \mu) = \begin{cases} 1, & x \in \Omega_\mu' \\ 0, & x \notin\Omega_\mu', \end{cases}
\end{equation}
where the subdomain $\Omega_\mu^{\prime}\subset\Omega_0$ is defined as $\Omega_\mu' = \{ x \in \Omega \mid \hat{s}(x; \mu) < 10^{-8}\}$, and require that
\begin{equation}
	\hat{s}(x; \mu) \leq C_\mu v(x; \mu), \quad \forall x \in \Omega_\mu',
\end{equation}
where $C_\mu > 1$ is a constant ensuring that $C_\mu v(x; \mu)$ bounds $\hat{s}(x; \mu)$ from above across the domain $\Omega_\mu'$.
This procedure is summarized in Algorithm \ref{alg:rejec_samp}.

\begin{algorithm}[hbt!]
\caption{Rejection sampling \cite{casella_generalized_accept_reject}}
    \begin{algorithmic}[1]
    	\State \textbf{Input:} Target density $\hat{s}(x; \mu)$, subdomain $\Omega_\mu'$, and constant $C_\mu > 1$.
		\State \textbf{Output:} Samples $\Xcal_\mu = \{x_1,\dots,x_{N_s'}\}$
	\State Initialize $\Xcal_\mu = \emptyset$
        \For{$i = 1, \ldots, N_s$}
            \State Uniformly sample $x_i \in \Omega_\mu'$
            \State Uniformly sample $u \in [0,1]$
            \If{$u \leq \frac{\hat{s}(x_i; \mu)}{C_\mu}$}
                \State $\Xcal_\mu \leftarrow \Xcal_\mu \cup \{ x_i \}$
            \EndIf
        \EndFor
    \end{algorithmic}
    \label{alg:rejec_samp}
\end{algorithm}

\begin{figure}[hbt!]
	\centering
	\input{figs/methodology/density_sampling/density_sampl.tikz}
	\caption{
	The probability density function of a two-dimensional multivariate normal distribution, defined by the mean vector $\mu_N = [0,0]$ and a covariance matrix $\Sigma = \text{diag} (10^{-3}, 7)$, along with samples generated via rejection sampling (\ref{line:density_sampl_sampl}) and their corresponding centroids (\ref{line:density_sampl_C}).
	}
	\label{fig:dens_samp}
\end{figure}
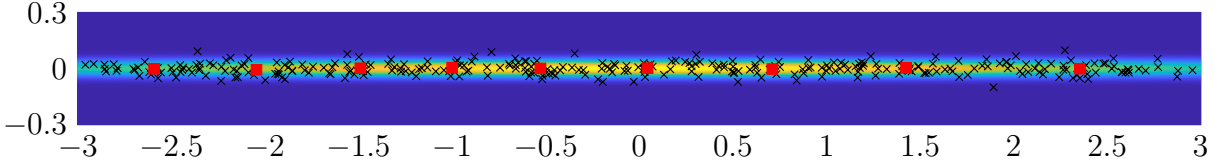

The final step employs $k$-means \cite{bishop_2006}, a widely-used clustering algorithm in unsupervised machine learning, to partition the samples $\Xcal_\mu$ obtained using rejection sampling in the previous step into $k$ distinct clusters $\Xcal_\mu = S_1^\mu \cup \cdots \cup S_k^\mu$ with $S_i^\mu\cap S_j^\mu=\emptyset$ if $i \neq j$ by solving the following optimization problem
\begin{equation}
	\operatorname*{\textrm{arg min}}_{S_1^\mu,\dots,S_k^\mu}
	\sum_{i=1}^{k} \sum_{x \in S_i^\mu} \| x - \bar{X}_i^\mu \|^2
	\mbox{,}
\end{equation}
where $\bar{X}_i^\mu$ is the centroid of points in cluster $S_i^\mu$.
This straightforward and effective method uncovers structure within the data, making it ideal for generating edge-based landmarks. Fig. \ref{fig:dens_samp} illustrates the second and third steps.
The centroids $\bar{X}_i^\mu$ will be used as radial basis function control points, and therefore will be referred to as \textit{control points} in the remainder.

\begin{remark}
    Other points of interest can provide additional landmarks. In this work, we also include the edge endpoints, that is, we use the start and end of each edge segment as landmarks.
\end{remark}

\begin{remark}
	The sampling process can be made more efficient in several ways. First, by selecting the optimal value of $C_\mu = \max_x \hat{s}(x; \mu)$, which is easily computable for low-dimensional problems. Additionally, algorithms such as importance sampling \cite{bishop_2006} can be used to reduce the rejection rate and improve convergence speed compared to traditional methods.
\end{remark}

\subsection{Correspondence problem}
\label{sec:correspondence}

\begin{figure}[hbt!]
    \centering
    \input{figs/methodology/parabola_correspondence/parabola_correspondence.tikz}
    \caption{Correspondence (\ref{line:par_corres_corres}) of landmarks (\ref{line:par_corres_C1}) between two distinct curves.}
    \label{fig:correspondence}
\end{figure}
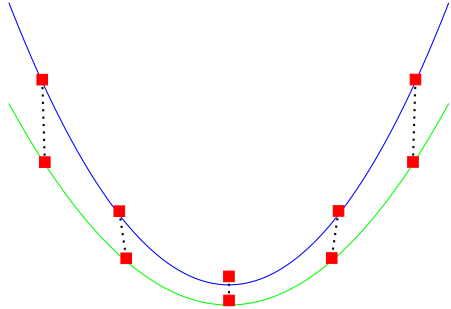

After placing an equal number of landmarks on different figures, the next step is to establish their correspondence. This involves identifying how specific parts of one image relate to corresponding parts of the other. Essentially, we need to map the landmarks from one figure to the appropriate locations in the other, ensuring that each point in the first image aligns with the correct point in the second. 
For similar features and small displacements, paring cluster centers by minimizing the total paring distance leads to accurate results. Let $\bar\Xcal_\mu = \{\bar{X}_1^\mu, \dots, \bar{X}_k^\mu\}$ denote the collection of control points for a given $\mu\in\Dcal$. Furthermore, let $\hat\mu\in\Dcal$ denote a \textit{reference} parameter to which features in solutions for any other $\mu\in\Dcal$ will be aligned. To establish the correspondence between the two control point collections, we seek the bijection $\tau_\mu: \bar{\mathcal{X}}_\mu \rightarrow \bar{\mathcal{X}}_{\hat\mu}$ that solves the following minimization problem
\begin{equation}
	\operatorname*{\textrm{minimize}}_{\tau_\mu}
	\sum_{\bar{X}_i^{\mu} \in \bar{\mathcal{X}}_{\mu}} \| \bar{X}_i^{\mu} - \tau_\mu(\bar{X}_i^{\mu}) \|_2
	\mbox{ .}
	\label{eq:correspondence}
\end{equation}
%
%
In Fig \ref{fig:correspondence}, correspondence can be obtained by either ordering the points horizontally or by minimizing total paring distance.

\begin{remark}
    The Hungarian algorithm \cite{kuhn_1955} provides a polynomial-time solution to Eq. \eqref{eq:correspondence}.
\end{remark}

\begin{remark}
    Landmarks can sometimes be easily matched if their distribution follows a certain pattern, e.g., if landmarks are expected to be ordered in a certain direction. For example, in Fig. \ref{fig:correspondence}, correspondence is trivially established after ordering the landmarks horizontally.
\end{remark}

\subsection{Mesh warping by radial basis functions}
\label{sec:mesh_warp}

A mapping from the  reference domain $X\in\Omega_0$ to the physical domain $x : \Dcal \rightarrow \Omega$ can be defined as
\begin{equation}
    x_i(\mu) = 
    X_i + \Delta X_i(X; \mu)
    \mbox{ ,}
\end{equation}
where $\Delta X_i : \Omega_0\times\Dcal\rightarrow\Rbb$ is the weighted sum of $n_c$ radial basis functions \cite{DEBOER_2007}, each associated to a unique control point $\bar{X}_j^\mu$ (cluster centroid) scaled by a coefficient $w_j^i(\mu)$, expressed by 
\begin{equation}
    \Delta X_i(X; \mu) = 
    \sum_{j=1}^{n_c} w_j^i(\mu) \phi \left( \frac{\| X - \bar{X}_j^\mu \|}{r} \right)
    \mbox{ ,}
    \label{eq:rbf_DX}
\end{equation}
where $r$ is a normalizing parameter controlling the support radius and $\phi : \Rbb_{\geq 0} \rightarrow \Rbb$ is a radial basis function kernel. This mapping function can be evaluated at any point $X$ once the weights have been determined.
In this work, we use the commonly adopted Wendland C2 basis function \cite{DEBOER_2007,wendland_1995},
\begin{equation}
    \phi (\eta) = 
    \begin{cases}
        (1 - \eta)^4 (4 \eta + 1) \mbox{,} & 0 \leq \eta \leq 1 \\
        0 \mbox{,} & \text{otherwise} \mbox{ .}
    \end{cases}
\end{equation}

The weights $w^i(\mu) = (w_1^i(\mu), \dots, w_{n_c}^i(\mu))$ are computed by solving the following linear system
\begin{equation}
    M_\mu w^i(\mu) = d^i(\mu)
    \mbox{ ,}
\end{equation}
where $d^i(\mu) = (d_1^i(\mu), \ldots, d_{n_c}^i(\mu))$ is the vector of control point displacement in the $i$th coordinate direction, with
\begin{equation}
	d_j^i(\mu) = [\tau_\mu(\bar{X}_j^\mu)]_i - [\bar{X}_j^\mu]_i
	\mbox{ ,}
\end{equation}
and $M_\mu \in \mathbb{R}^{n_c \times n_c}$ is defined as
\begin{equation}
    (M_\mu)_{ab} = \phi \left( \frac{\| \bar{X}_a^\mu - \bar{X}_b^\mu \|}{r} \right)
    \mbox{ .}
\end{equation}
For multi-dimensional mesh deformation, this process must be repeated for each dimension, $i=1,\dots,d$.

\begin{figure}[hbt!]
	\centering
	\begin{subfigure}{0.49\textwidth}
		\centering
		\input{figs/methodology/arc/arc0.tikz}
		\caption{}
	\end{subfigure}
	\hfill
	\begin{subfigure}{0.49\textwidth}
		\centering
		\input{figs/methodology/arc/arc1.tikz}
		\caption{}
	\end{subfigure}
	\\
	\begin{subfigure}{0.49\textwidth}
		\centering
		\input{figs/methodology/arc/arc0_cntl_pnts.tikz}
		\caption{}
	\end{subfigure}
	\hfill
	\begin{subfigure}{0.49\textwidth}
		\centering
		\input{figs/methodology/arc/arc0_match_cntl_pnts.tikz}
		\caption{}
	\end{subfigure}
	\caption{
		The figure in the top-left can be aligned with the figure in the top-right by placing control points, or landmarks, on the shock (\ref{line:arc0_match_arc}), the shock endpoints (\ref{line:arc0_match_end}), and the boundary (\ref{line:arc0_match_bnd}), then displacing them horizontally by the appropriate amount.
		This example uses the Wendland C2 basis function with a support radius of \( r = 10 \).
	}
	\label{fig:arc-main}
\end{figure}
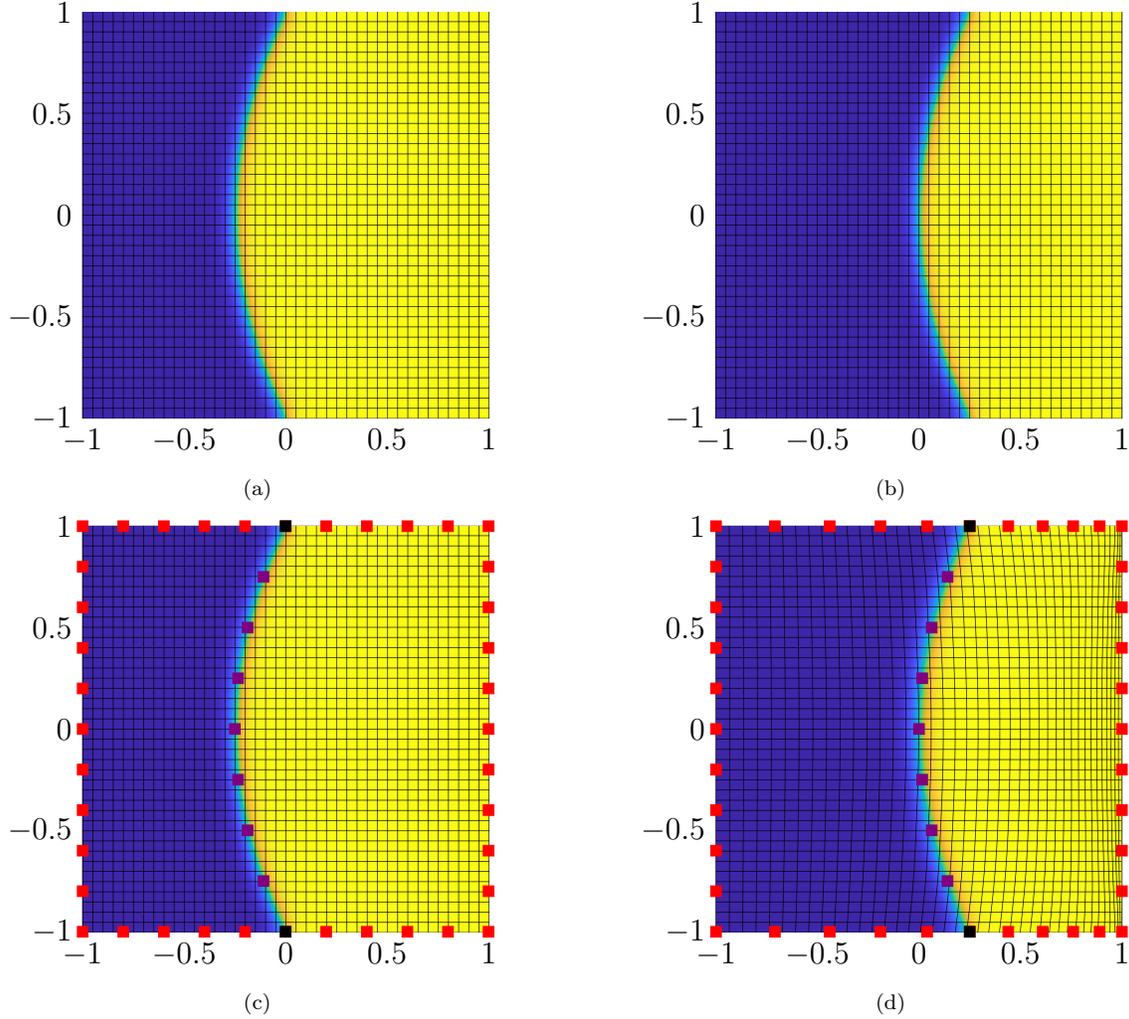

From these definitions, we define the mapping degrees of freedom $x_h : \Dcal \rightarrow \Rbb^{N_x}$ that uses landmark-based registration to align an image $Q_h(x_h(\mu); \mu)$ with the reference image $Q_h(X_h; \hat\mu)$ as
\begin{equation}
[x_h(\mu)]_I = [X_h]_I + \Delta X([X_h]_I; \mu)
\end{equation}
for $I = 1,\dots,\hat{N}_x$, where $\Delta X$ is defined in Eq. \eqref{eq:rbf_DX}.
Figure \ref{fig:arc-main} illustrates how features can be aligned through the positioning and displacement of control points in this framework.

\begin{remark}
	Including control points on the boundary is essential for preventing boundary nodes from leaving the domain, particularly when large deformations are applied or a larger support radius $r$ is used. While a smaller radius can sometimes eliminate the need for boundary control points, a larger radius is often preferred for smoother approximations, which reinforces the importance of maintaining boundary constraints.
\end{remark}

\section{Reduced-order model with implicit feature alignment}
\label{sec:rom}

Following a standard projection-based reduced-order modeling approach, $Q_h$ is approximated in a low-dimensional subspace of dimension $n_q$, with $n_q \ll N$. Specifically, we write
\begin{equation}
    Q_h (x_h(\mu), \mu) \approx 
    \hat{Q}(\mu) = 
    Q_r(\mu) + \Phi a (\mu)
    \mbox{ ,}
\end{equation}
where $\hat{Q} : \Dcal \rightarrow \mathbb{R}^N$ is a reduced basis approximation to $Q_h$, $\Phi \in \mathbb{R}^{N \times n_q}$ is the basis for the reduced subspace, $Q_r : \Dcal \rightarrow \mathbb{R}^{N}$ is a reference state, $a : \Dcal \rightarrow \mathbb{R}^{n_q}$ are the corresponding reduced coefficients, and $x_h(\mu)$ is the mapping degrees of freedom such that the features of $Q_h (x_h(\mu), \mu)$ align with the features of $Q_h(X_h; \hat\mu)$ (Section~\ref{sec:mesh_warp}).
At the high-dimensional model (HDM) stage, the snapshots $Q_h(X_h,\mu_i)$ are obtained by solving the discretized system of equations on a fixed grid, where the reference mesh equals the physical mesh. Consequently, these solutions are not aligned due to parameter-dependent variation in the solution features in the physical domain, making them ill-suited as reduced basis snapshots. We use the landmark-based registration approach of Section~\ref{sec:mesh_warp} to compute the mapping degrees of freedom $x_h(\mu_i)$ that aligns the original snapshots $Q_h(X_h,\mu_i)$ with a reference snapshot $Q_h(X_h,\hat\mu)$ to obtained the aligned snapshots $Q_h(x_h(\mu_i), \mu_i)$.  The reduced basis $\Phi$ is then constructed using proper orthogonal decomposition (POD) \cite{Sirovich_snap1} applied to the aligned perturbed states, defined as
\begin{equation}
	\Phi = 
	{\tt POD}_{N,n}^{n_q} ([Q_h(x_h(\mu_1), \mu_1) - Q_r(\mu_1),  \ldots, Q_h(x_h(\mu_{n}), \mu_{n}) - Q_r(\mu_{n})])
	\mbox{ ,}
	\label{eq:pod}
\end{equation}
where ${\tt POD}_{N,n}^{n_q} : \mathbb{R}^{N \times n} \rightarrow \mathbb{R}^{N \times n_q}$ performs the singular value decomposition (SVD) on the aligned perturbed snapshot matrix and selects the $n_q$ left singular vectors. 

%

Similarly, we approximate $x_h(\mu)$ as
\begin{equation}
    x_h(\mu) \approx 
    \hat{x}(\mu) = 
    X_h + \Psi y (\mu)
    \mbox{ ,}
\end{equation}
where $\hat{x} : \Dcal \rightarrow \Rbb^{N_x}$ is the low-dimensional parametrization of the domain mapping coefficients, $y : \Dcal \rightarrow \mathbb{R}^{n_x}$ are the reduced coefficients, and $\Psi \in \mathbb{R}^{N_x \times n_x}$ is the mesh reduced basis, constructed from the mesh displacements as
\begin{equation}
    \Psi = [x_h(\mu_1) - X_h, \ldots, x_h(\mu_n) - X_h]
    \mbox{ .}
\end{equation}

For a given parameter $\mu$, the reduced coefficients for both the state $a(\mu)$ and the mesh $y(\mu)$ are determined as the solution to the following optimization problem
\begin{equation}
    (a(\mu), y(\mu)) =
    \operatorname*{\textrm{arg min}}_{(a,y) \in \Rbb^{n_q} \times \Rbb^{n_x}}
    \frac{1}{2} 
    \| R (Q_r(\mu) + \Phi a; X_h + \Psi y, \mu)  \|_{\Theta}^2
    \mbox{ ,}
    \label{eq:lspg}
\end{equation}
where $\Theta \in \mathbb{R}^{N \times N}$ is a symmetric positive definite (SPD) matrix defining the norm $\| R \|_{\Theta} = \sqrt{R^T \Theta R}$, assumed to be the identity matrix $I \in \mathbb{R}^{N \times N}$ unless stated otherwise. The optimization problem given by Eq. \eqref{eq:lspg} is provably equivalent to a Petrov-Galerkin reduced-order model \cite{grimberg2020stability}.

\begin{remark}
	Despite dimensionality reduction ($n_q \ll N$), the optimization problem in Eq. \eqref{eq:lspg} still scales with $N$, a well-known issue in nonlinear reduced-order modeling. Recent work \cite{2024_Mirhoseini} on implicit feature tracking addresses this bottleneck through hyperreduction.
\end{remark}

\begin{remark}
Eq. \eqref{eq:pod} suggests the HDM solution is recomputed from scratch according to Eq. \eqref{eq:hdm} after the alignment mapping $x_h(\mu)$ is available. Another option is to use the mapping $x_h(\mu)$ to deform the mesh of the reference domain then interpolate the solution $Q_h(X_h; \mu)$ onto this deformed grid. While both approaches are valid, we use the latter in this work because it is simple and does not require additional HDM solutions.
\end{remark}

\section{Numerical experiments}
\label{sec:num_exp}

In this section, we solve two types of ROMs, one based on a fixed mesh and the other leveraging implicit feature tracking, to address two convection-dominated viscous flow problems. We begin by solving the space-time Burgers' equation (Section \ref{sec:num_exp_burgers}) parameterized by the initial condition. The second numerical experiment (Section \ref{sec:num_exp_cy}) is hypersonic viscous flow over a cylinder parameterized by the Mach number. This problem features both a shock and a boundary layer, making it more challenging.

\subsection{Space-time Burgers shock advection}
\label{sec:num_exp_burgers}

Consider the one-dimensional time-dependent Burgers' equation,
\begin{equation}
	\pder{q}{t} + \pder{}{x}\left(\frac{1}{2}q^2 - \nu\pder{q}{x}\right) = 0, 
	\quad q(x,0) = q_0 (x)
	\mbox{,}
\end{equation}
where $\nu \in \mathbb{R}^+$ is the diffusion coefficient and the parametrized initial condition $q_0(x;\mu)$ is defined as
\begin{equation}
	q_0 (x; \mu) = 
	\begin{cases}
		\mu \mbox{,} & x < - 0.5 \\
		0 \mbox{,} & \text{otherwise} \mbox{ }
	\end{cases}
\end{equation}
where $\mu \in \Dcal\coloneqq\mathbb{R}$. Boundary conditions are consistent with the initial condition at inflow boundaries.
We adopt the spatial-temporal domain $\Omega \times \mathcal{T} \coloneqq (-1,1) \times (0,1)$ and diffusion coefficient $\nu = 10^{-3}$. Varying $\mu$ affects both the step height and the wave speed.
The spatial domain is discretized with 100 elements, which are extruded in time and split in 50 time instances to form a linear ($p_m = 1$) mesh with $4900$ quadratic elements with a quartic ($p = 4$) solution.

\begin{figure}[hbt!]
	\centering
	\begin{subfigure}{0.49\textwidth}
		\centering
		\includegraphics[width=\textwidth,trim={0mm 0mm 0mm 0mm},clip]{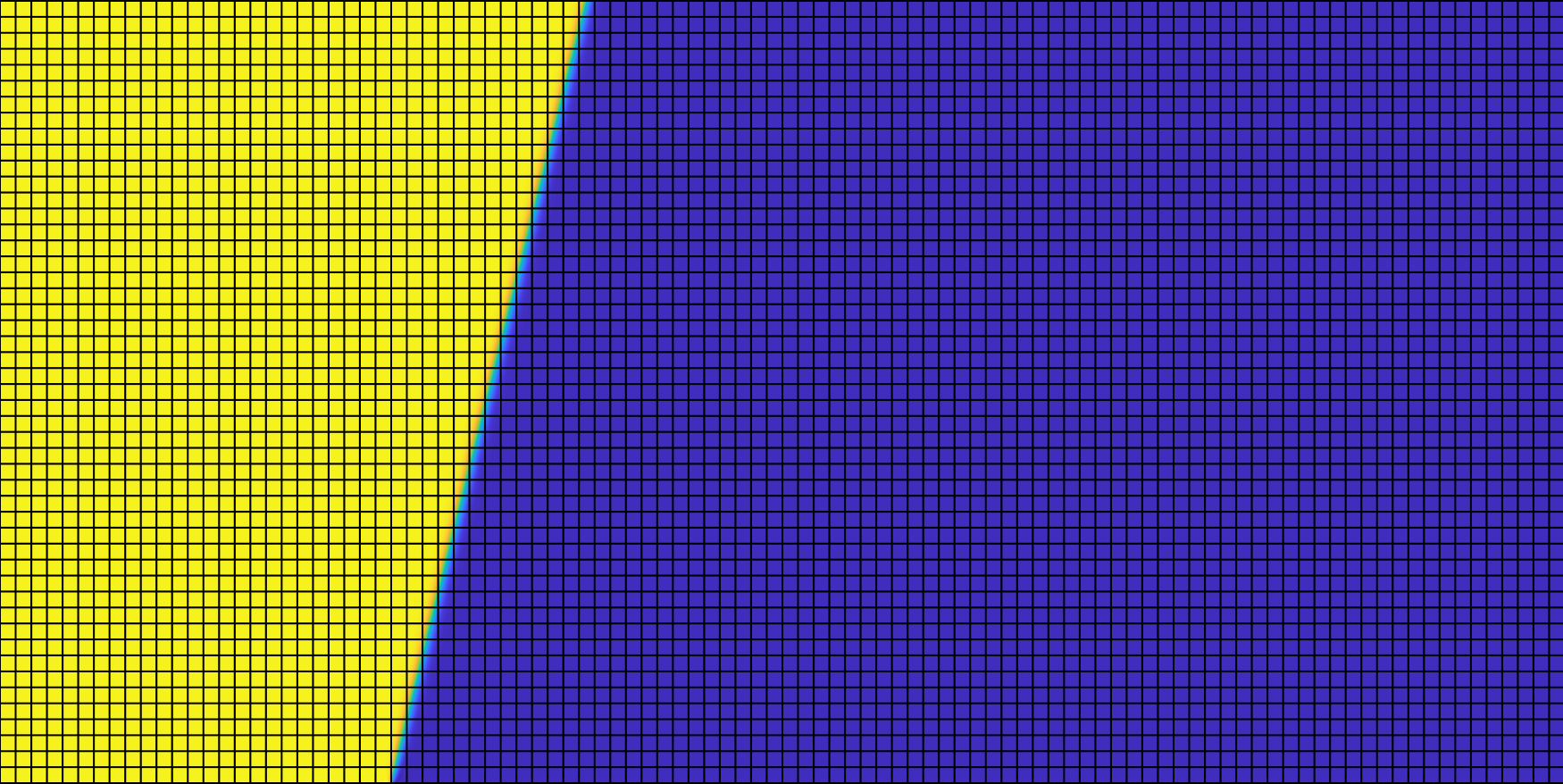}
		\caption{\( \mu = 0.5 \)}
		\label{fig:burgers_Q1}
	\end{subfigure}
	\hfill
	\begin{subfigure}{0.49\textwidth}
		\centering
		\includegraphics[width=\textwidth,trim={0mm 0mm 0mm 0mm},clip]{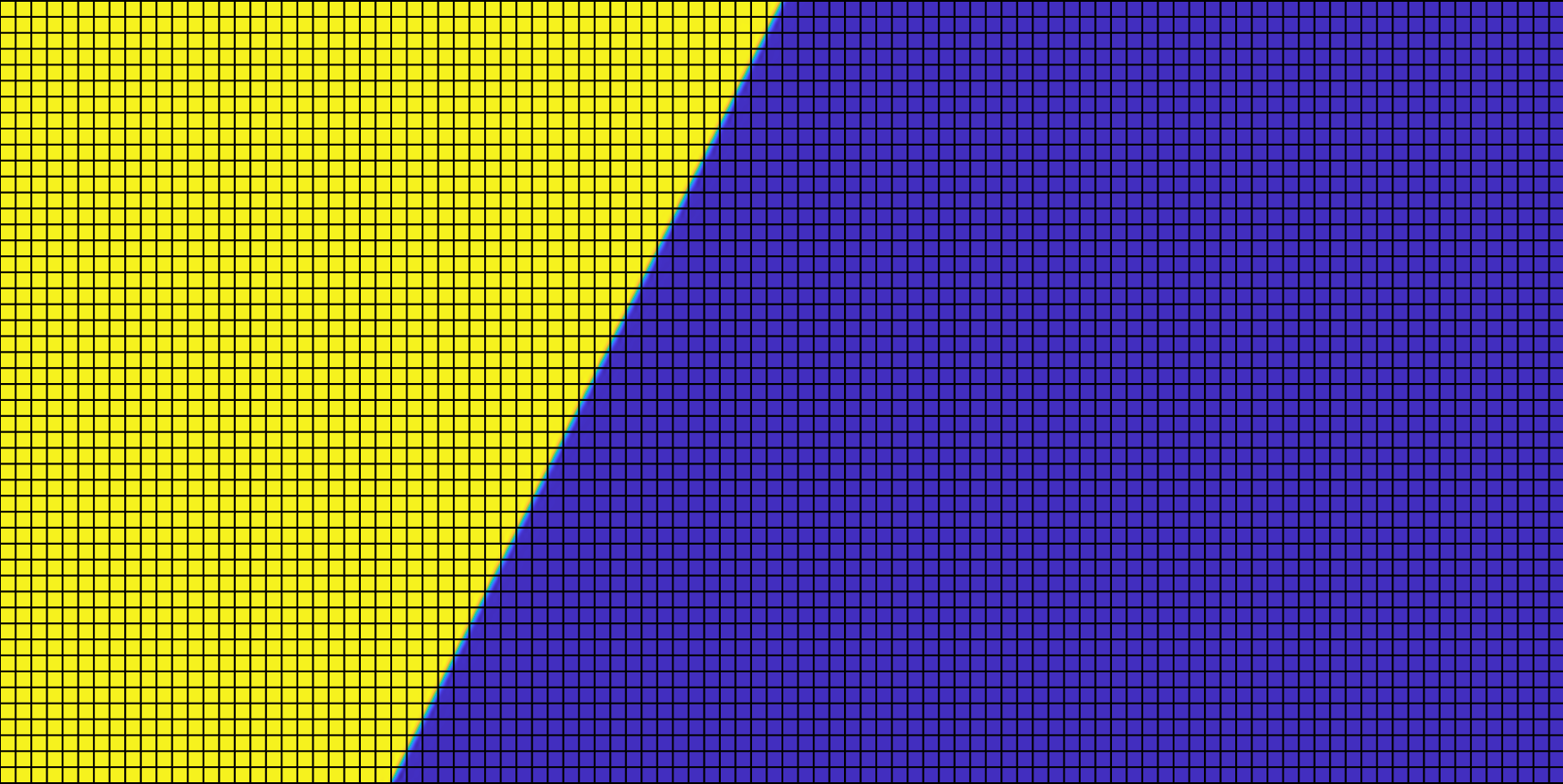}
		\caption{\( \mu = 1 \)}
		\label{fig:burgers_Q2}
	\end{subfigure}
	\caption{Burgers' equation: HDM solution training data.}
	\label{fig:burgers_data}
\end{figure}

\begin{figure}[hbt!]
	\centering
	\begin{subfigure}{0.49\textwidth}
		\centering
		\includegraphics[width=\textwidth,trim={0mm 0mm 0mm 0mm},clip]{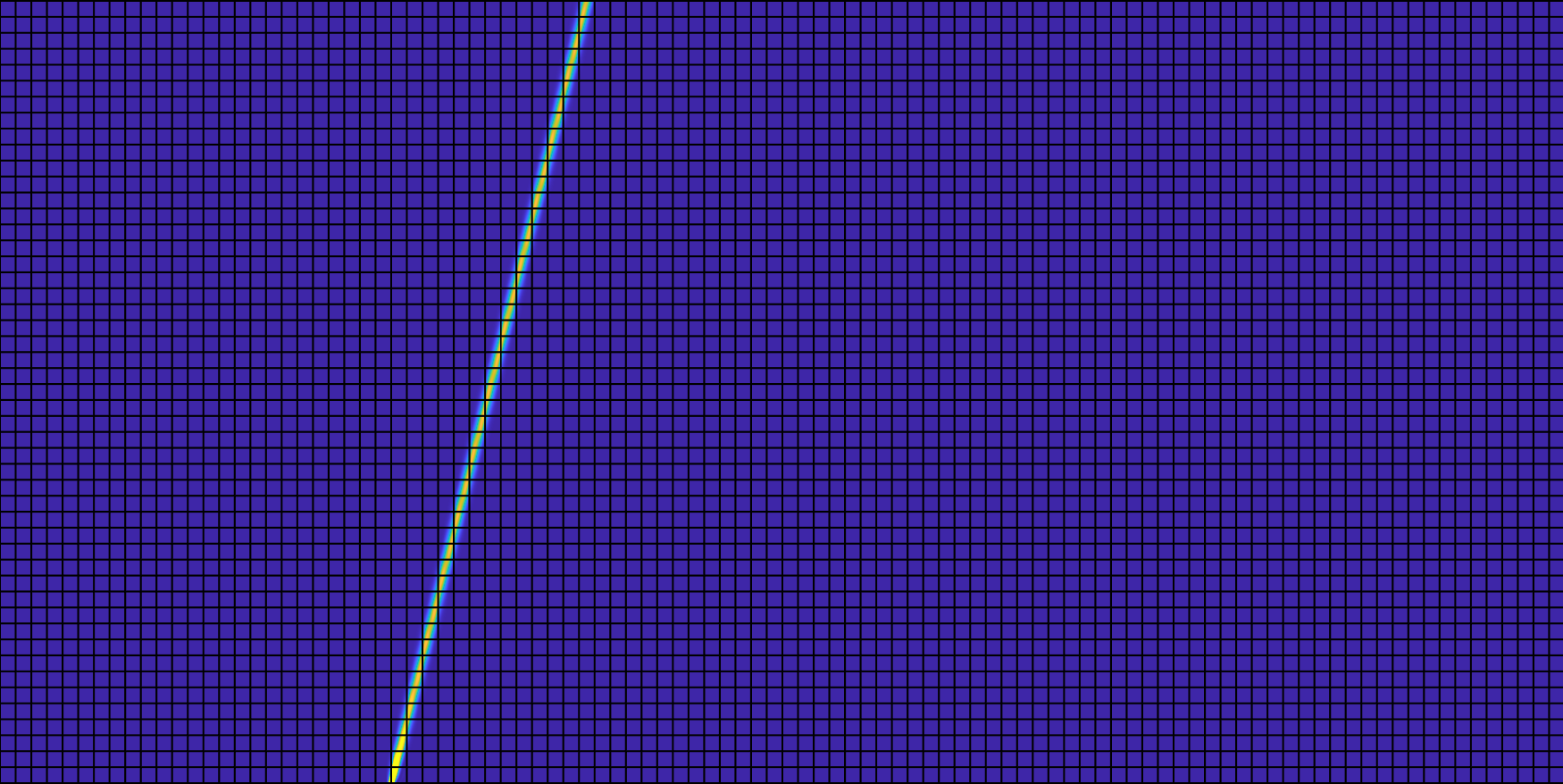}
		\caption{\( \mu = 0.5 \)}
		\label{fig:burgers_sensor_Q1}
	\end{subfigure}
	\hfill
	\begin{subfigure}{0.49\textwidth}
		\centering
		\includegraphics[width=\textwidth,trim={0mm 0mm 0mm 0mm},clip]{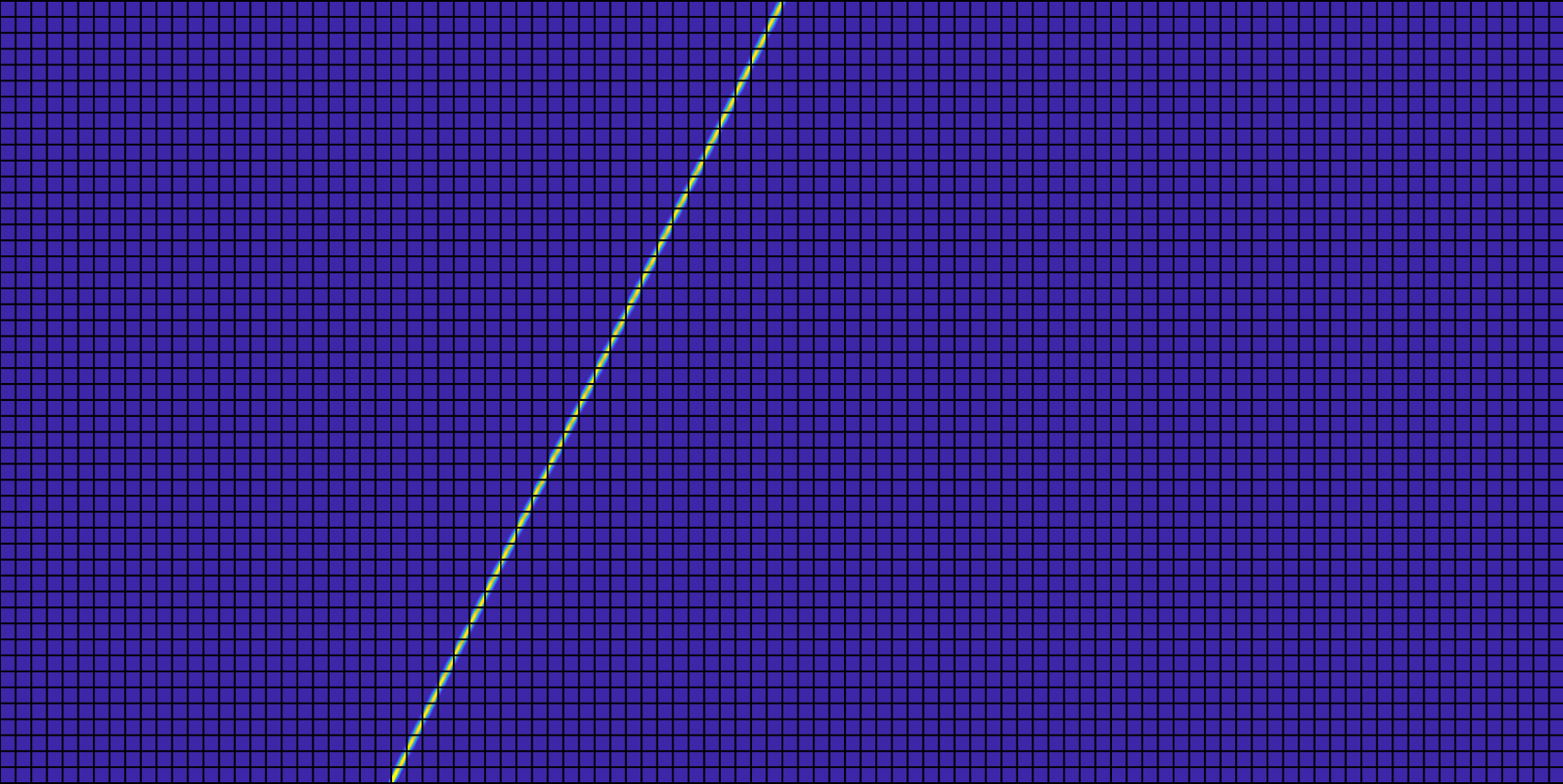}
		\caption{\( \mu = 1 \)}
		\label{fig:burgers_sensor_Q2}
	\end{subfigure}
	\caption{Burgers' equation: Edge sensor solution as defined by Eq.\eqref{eq:shk_sensor_grad}.}
	\label{fig:burgers_sensor}
\end{figure}

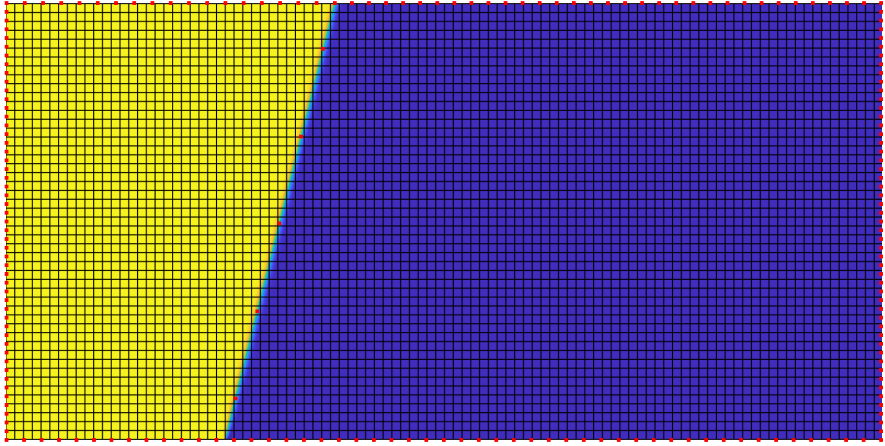
\begin{figure}[hbt!]
	\centering
	\input{figs/burgers/control_points/burgers_xc.tikz}
	\caption{Burgers' equation: Placement of control points (\ref{line:burgers_xc}) on the reference snapshot.
	}
	\label{fig:burgers_xc}
\end{figure}

We start by solving the HDM on a fixed mesh (i.e., $x_h = X_h$) for two parameter values, $\mu \in \{ 0.5, 1\}$, with the results shown in Fig. \ref{fig:burgers_data}.
Next, the shocks are identified using the gradient-based edge sensor described in Eq. \ref{eq:shk_sensor_grad}, with the corresponding edge sensors shown in Fig. \ref{fig:burgers_sensor}. Control points are then positioned along the shocks ($k = 5$) and boundaries, as depicted in Fig. \ref{fig:burgers_xc}.
We use the solution for $\mu = 0.5$ as the reference, aligning the other solution ($\mu = 1$) to match the shock position of this reference. The radius of support for the RBF functions set to $r = 100$.


For our numerical experiments, we construct two types of ROMs:  one based on a fixed mesh and the other leveraging implicit feature tracking. Furthermore, these ROMs are evaluated at two distinct parameter values, $\mu \in \{ 0.75, 2\}$, with one value lying within the training window and the other outside it.
For the reduced basis $\Phi$, the fixed mesh ROM is constructed using the two snapshots collected during the training step. In contrast, the implicit feature tracking ROM uses only the snapshot corresponding to $\mu = 0.5$. This choice is made because, for this problem, the snapshots become nearly linearly dependent after alignment. Additionally, for both ROMs, the reference state is defined as the zero vector, $Q_r(\mu) = [0, \ldots, 0]^T$.

\begin{figure}[hbt!]
	\centering
	\begin{subfigure}{0.49\textwidth}
		\centering
		\includegraphics[width=\textwidth,trim={0mm 0mm 0mm 0mm},clip]{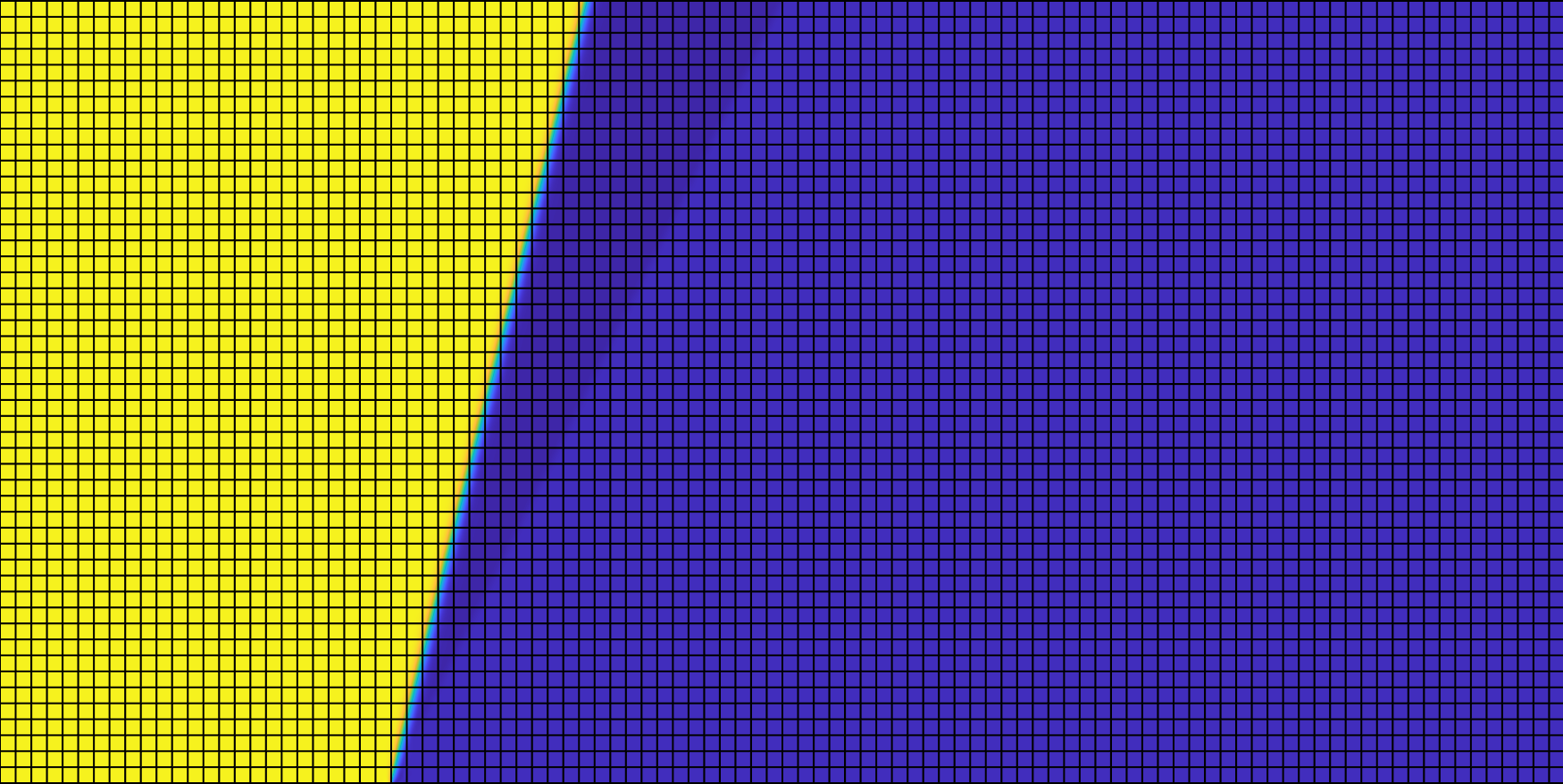}
		\caption{Fixed mesh ROM.}
		\label{fig:burgers_rom_Q75}
	\end{subfigure}
	\hfill
	\begin{subfigure}{0.49\textwidth}
		\centering
		\includegraphics[width=\textwidth,trim={0mm 0mm 0mm 0mm},clip]{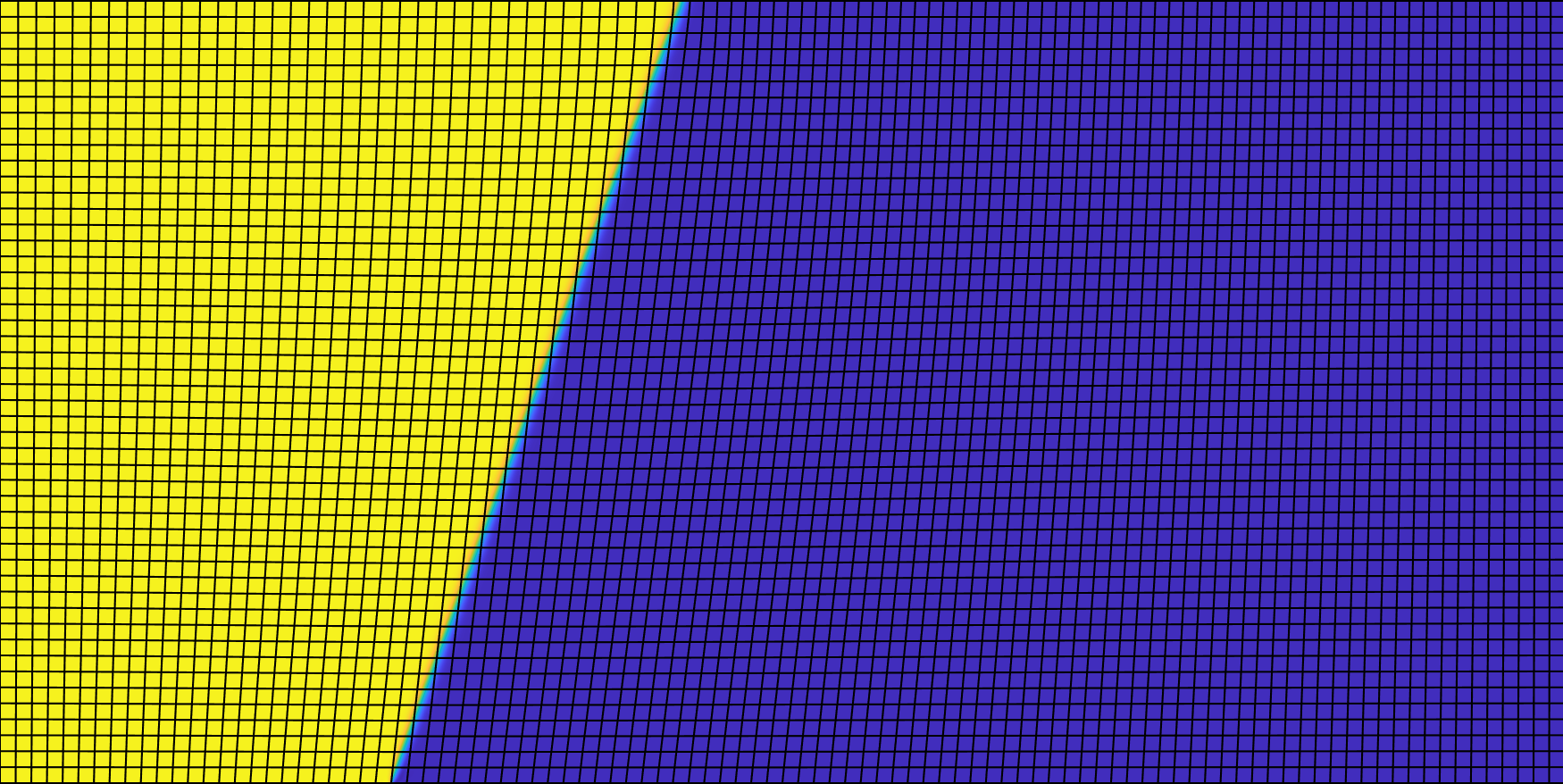}
		\caption{Implicit feature tracking ROM.}
		\label{fig:burgers_romft_Q75}
	\end{subfigure}
	\caption{Burgers' equation: Space-time reduced-order model solutions for the parameter \( \mu = 0.75 \).}
	\label{fig:burgers_rom_surf_75}
\end{figure}
\begin{figure}[hbt!]
	\centering
	\begin{subfigure}{0.49\textwidth}
		\centering
		\includegraphics[width=\textwidth,trim={0mm 0mm 0mm 0mm},clip]{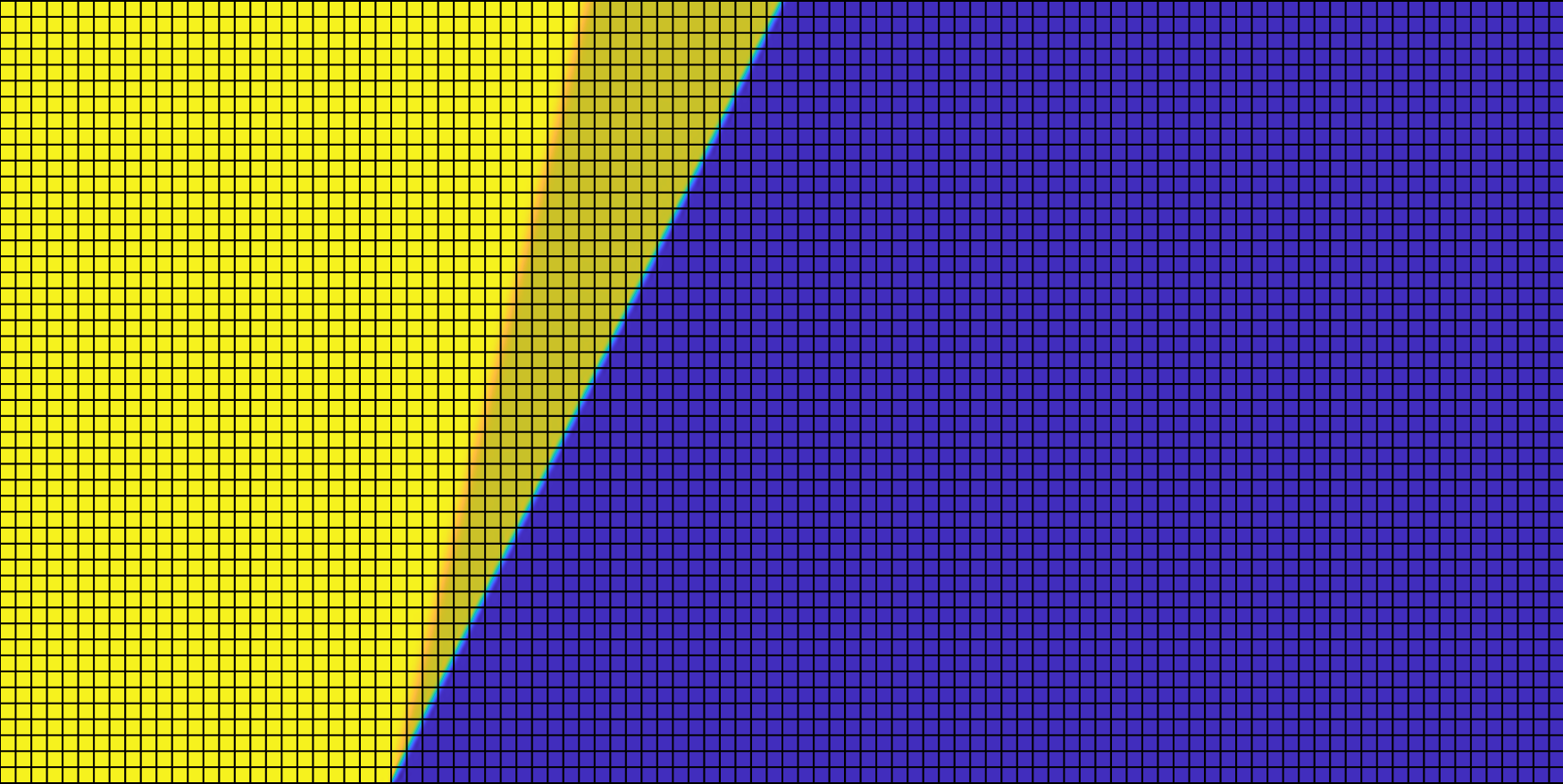}
		\caption{Fixed mesh ROM.}
		\label{fig:burgers_rom_Q200}
	\end{subfigure}
	\hfill
	\begin{subfigure}{0.49\textwidth}
		\centering
		\includegraphics[width=\textwidth,trim={0mm 0mm 0mm 0mm},clip]{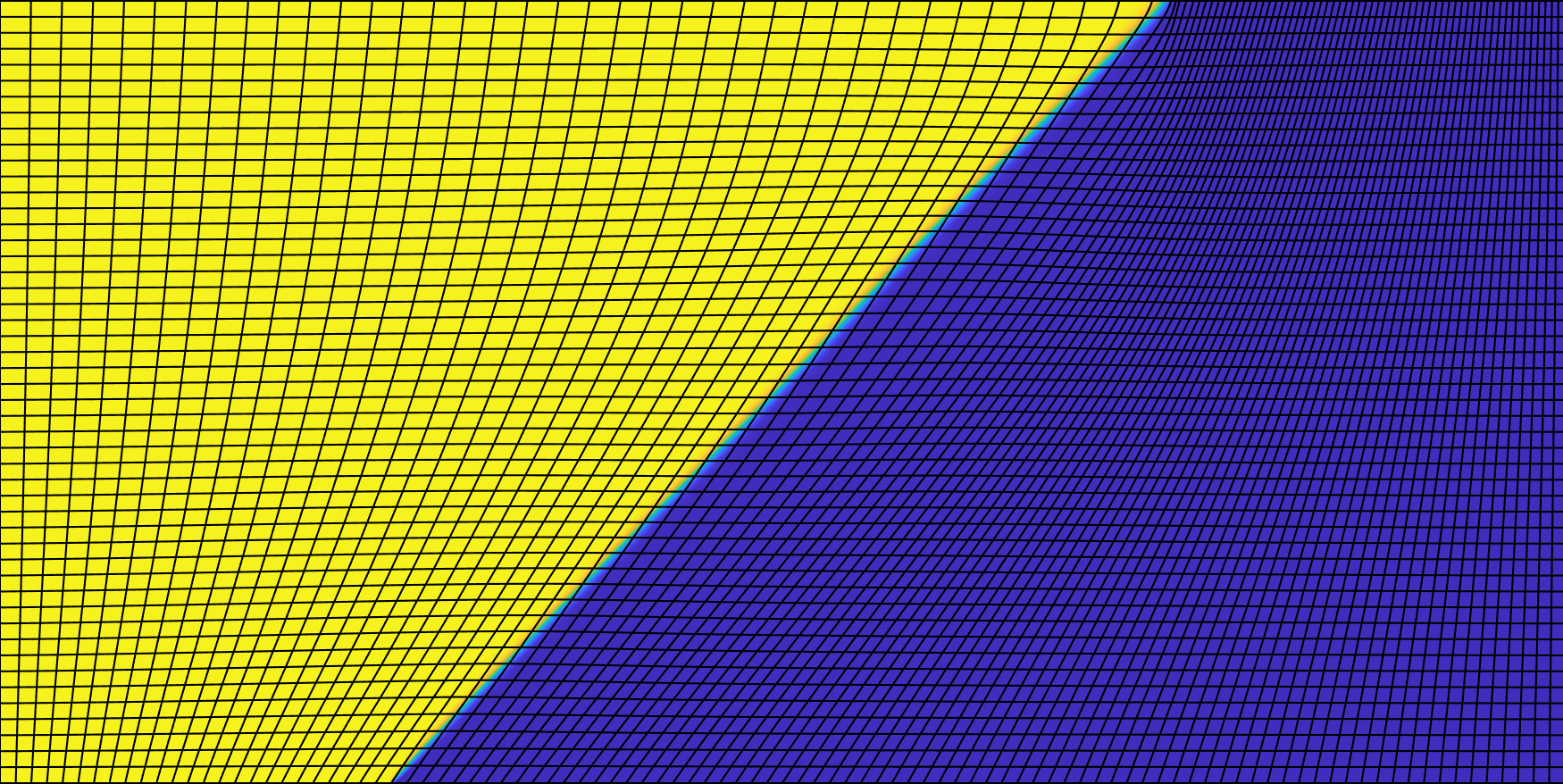}
		\caption{Implicit feature tracking ROM.}
		\label{fig:burgers_romft_Q200}
	\end{subfigure}
	\caption{Burgers' equation: Space-time reduced-order model solutions for the parameter \( \mu = 2 \).}
	\label{fig:burgers_rom_surf_200}
\end{figure}

Figures \ref{fig:burgers_rom_surf_75} and \ref{fig:burgers_rom_surf_200} present the space-time solutions for $\mu = 0.75$ and $\mu = 2$, respectively, with the corresponding solutions at $t=1$ shown in Figures \ref{fig:burgers_rom_slice_75} and \ref{fig:burgers_rom_slice_200}.
%
%
As expected and previously discussed, the fixed mesh ROMs produce staircase-like solutions, being particularly inaccurate when the parameter lies outside the training window. While using more data could improve accuracy within the training window, as the two snapshots used here are insufficient, additional data would not make the ROM predictive outside the training window.
In contrast, implicit feature ROMs successfully capture the convective nature of the solution, resulting in accurate predictions both inside and outside the training window.


\begin{figure}[hbt!]
	\centering
	\input{figs/burgers/slice/75/burgers_rom_75_slice.tikz}
	\caption{Burgers' equation: Reduced-order model solutions for the parameter $\mu = 0.75$ at $t=1$. Results are presented for the fixed mesh ROM solution (\ref{line:burgers_rom_75_slice_ref}), the implicit feature tracking ROM solution (\ref{line:burgers_romft_75_slice}), and the position the initial discontinuity would have reached if traveling at a speed of $\mu / 2$ (\ref{line:burgers_75_slice_ref}).
	}
	\label{fig:burgers_rom_slice_75}
\end{figure}


\begin{figure}[hbt!]
	\centering
	\input{figs/burgers/slice/200/burgers_rom_200_slice.tikz}
	\caption{Burgers' equation: Reduced-order model solutions for the parameter $\mu = 2$ at $t=1$. Results are presented for the fixed mesh ROM solution (\ref{line:burgers_rom_75_slice_ref}), the implicit feature tracking ROM solution (\ref{line:burgers_romft_75_slice}), and the position the initial discontinuity would have reached if traveling at a speed of $\mu / 2$ (\ref{line:burgers_75_slice_ref}).
	}
	\label{fig:burgers_rom_slice_200}
\end{figure}

\subsection{Hypersonic flow over a cylinder}
\label{sec:num_exp_cy}

In this section, we first present the compressible Navier-Stokes equations and the PDE-based shock-capturing method utilized in the analysis. We then introduce a sensor designed to compute the weights $\omega$ for a boundary layer-focused inner product matrix. We close with a direct comparison between a fixed-mesh ROM and implicit feature tracking for Mach-parametrized hypersonic flow over a cylinder.

\subsubsection{Governing equations}

Compressible, viscous flow of an ideal gas is governed by the Navier-Stokes equations, which can be written as a conservation law of the form Eq. \eqref{eq:claw} with
\begin{equation}
	q = \begin{bmatrix} \rho \\ \rho v \\ \rho E\end{bmatrix}\mbox{, }\quad
	f^{i} (q) = 
	\begin{bmatrix}
		\rho v^T \\
		\rho v v^T + P I_{d\times d} \\
		(\rho E + P) v^T
	\end{bmatrix}
	\mbox{, }\quad
	f^{v} (q, \nabla q) = 
	\begin{bmatrix}
		0 \\
		\tau \\
		v^T \tau + \kappa \nabla T^T
	\end{bmatrix}.
\end{equation}
The density $\rho(x)\in\Rbb_{\ge 0}$, velocity $v(x)\in\Rbb^d$, and total energy $E(x)\in\Rbb$ of the fluid are implicitly defined as the solution of Eq. \eqref{eq:claw}. Furthermore, $P(x) \in \Rbb_{\ge 0}$ is the pressure, $\kappa$ is the thermal conductivity, $T = P/(\rho R)$ is the temperature (ideal gas assumption), and $R$ is the gas constant. The pressure is related to the state vector by the equation of state, $P = (\gamma - 1) \rho (E - v \cdot v / 2)$, where $\gamma$ is the ratio of specific heats. The shear stress is, $\tau = \mu \left(\nabla v + \nabla v^T  - \frac{2}{3} \mathrm{tr }(\nabla v) \right)$, where $\mu$ is the dynamic viscosity. The Reynolds number is given by $\text{Re} = \rho v L / \mu$, where $L$ is the characteristic length. The Mach number is defined as $M = \norm{v} / c$, with the speed of sound given by $c = \sqrt{\gamma P / \rho}$.

Shocks are stabilized using PDE-based artificial viscosity \cite{barter_2010}. In this approach, the following artificial viscosity equation is appended to the compressible Navier-Stokes equations
\begin{equation}
	\nabla \cdot \left( \eta \nabla \epsilon \right) +
	\left[ \frac{\bar{h}(x)}{p \lambda_{\text{max}}(q)} S_k(q) - \epsilon \right] = 0
	\mbox{,}
	\label{eq:av}
\end{equation}
where $\epsilon : \mathbb{R}^d \times \mathbb{R}^+ \to \mathbb{R}$ is the artificial viscosity, $\eta \in \mathbb{R}^{d \times d}$ is the diffusivity, $S_k(q)$ is a shock indicator, $\bar{h} (x)$ is a measure of the mesh element size, $p$ is the polynomial degree of the DG approximation to be used, and $\lambda_\mathrm{max}(q)$ is the maximum wave speed associated with the state $q$.
The dynamic viscosity is modified as $\mu = \mu_0 + \rho \epsilon$, where $\mu_0$ is the physical dynamic viscosity of the fluid and $\rho \epsilon$ is the artificial dynamic viscosity.
In this work, we use a dilatation-based shock sensor $S_k(q) = -\min\{\mathrm{tr}(\nabla v), 0\}$, similar to \cite{moro2016dilation}.
Further details and discussion can be found in \cite{barter_2010}.
In the end, the Navier-Stokes equations with PDE-based artificial viscosity is a conservation law of the form in Eq. \eqref{eq:claw} with
\begin{equation}
\begin{aligned}
	q = \begin{bmatrix} \rho \\ \rho v \\ \rho E \\ \epsilon\end{bmatrix}&\mbox{, }\quad
	&&f^{i} (q) = 
	\begin{bmatrix}
		\rho v^T \\
		\rho v v^T + P I_{d\times d} \\
		(\rho E + P) v^T \\
		0
	\end{bmatrix}
	\mbox{, } \\
	f^{v} (q, \nabla q) = 
	\begin{bmatrix}
		0 \\
		\tau \\
		v^T \tau + \kappa \nabla T^T \\
		\eta \nabla \epsilon
	\end{bmatrix}
	&\mbox{, }\quad
	&&s(q, \nabla q) =
	\begin{bmatrix}
	0 \\ 0_d \\ 0 \\ \dfrac{\bar{h}(x)}{p \lambda_\mathrm{max}(q)} S_k(q) - \epsilon
	\end{bmatrix}.
\end{aligned}
\end{equation}

\subsubsection{Inner product matrix guided by boundary layer sensor}

For problems with boundary layers, such as hypersonic flow over a cylinder considered next, we have found it advantageous to use an inner product matrix $\Theta$ that is biased toward the boundary layer.
To this end, we introduce the sensor $s_b (x) : \mathbb{R}^d \rightarrow \mathbb{R}$ to detect boundary layers, defined as
\begin{equation}
	s_b (x) = f_w (x) \frac{\hat{s} (x)}{\hat{s} (x_w)}
	\mbox{ ,}
	\label{eq:bl_sensor}
\end{equation}
where $x_w$ represents the position of the nearest wall and $f_w$ is a predefined function that equals $1$ at the wall and quickly approaches zero as the distance from the wall increases.
%
%
Additionally, we define the inner product matrix $\Theta$ as a diagonal matrix, given by
\begin{equation}
    \Theta_{ii} (x; \mu) = 
    \begin{cases}
        \sqrt{\omega} \mbox{,} & \mbox{if } i \in \mathcal{S}_b \mbox{ ,}\\
        1 \mbox{,} & \text{otherwise} \mbox{ ,}
    \end{cases}
\end{equation}
where $\mathcal{S}_b = \{ i \in \{1,\dots,\hat{N}_x\} \mid s_b((X_h)_i) \geq 0.01\}$ is the collection of boundary layer nodes, and $\omega$ is a weight parameter.

\subsubsection{Numerical results}

\begin{figure}
	\centering
	\includegraphics[clip,width=1.0\textwidth]{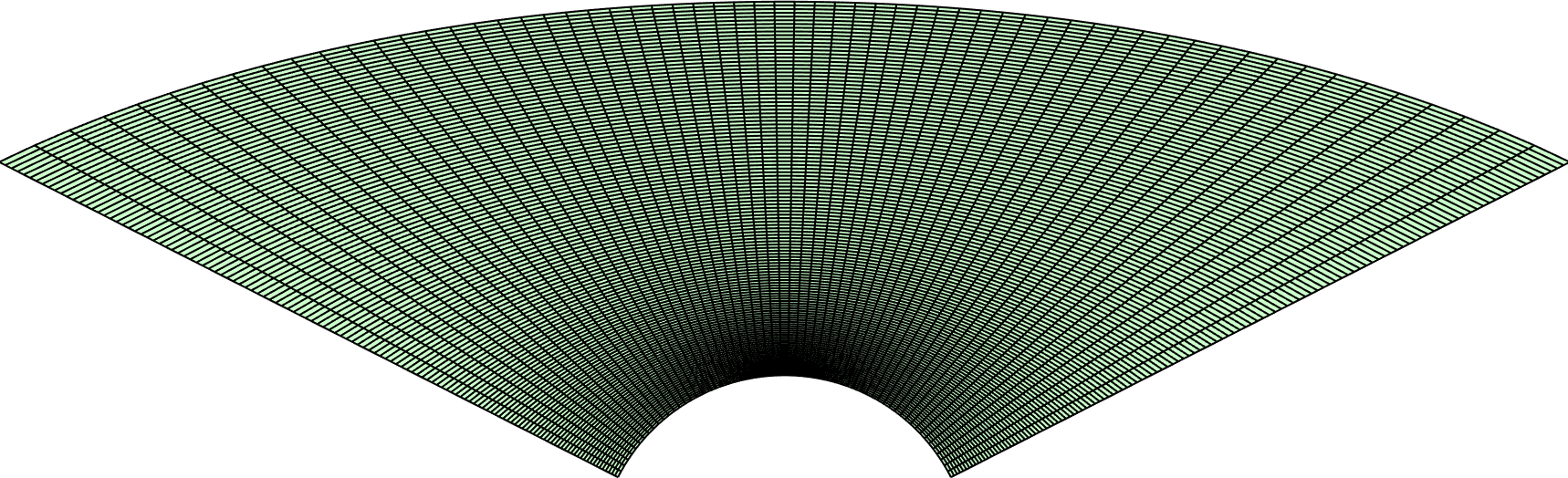}
	\caption{Hypersonic cylinder: Reference computational mesh.}
	\label{fig:cy_msh}
\end{figure}

Now, we compare the proposed implicit feature tracking method to fixed-mesh ROMs on hypersonic flow over a cylinder, parametrized by the inflow Mach number ($\mu = M_\infty$). The flow domain $\Omega\subset\Rbb^2$ an angular slice (polar angle $\theta \in [3\pi/20, 17\pi/20]$) between a circle and an ellipse, representing a segment of the cylinder and the surrounding flow domain. The circle, centered at the origin $(0, 0)$ with a radius of $1$, forms the cylinder wall. The ellipse, also centered at $(0, 0)$, has a semi-major axis of $6$ along the streamwise $x$-direction and a semi-minor axis of $3$ along the vertical $y$-direction, defines the outer farfield boundary. The computational domain is discretized with $7{,}503$ quadratic ($p_m = 2$) quadrilateral elements, supporting a state cubic solution ($p = 3$), as depicted in Fig. \ref{fig:cy_msh}. An isothermal no-slip condition is enforced along the cylinder wall at $r = 1$, with velocity $v=0$ and a constant wall temperature ($T_w = 2.5$). Supersonic flow enters the domain on the outer elliptical boundary and exits along the lateral sides of the flow domain.

\begin{figure}[hbt!]
	\centering
	\begin{subfigure}{0.49\textwidth}
		\centering
		\includegraphics[width=\textwidth,trim={0mm 0mm 0mm 0mm},clip]{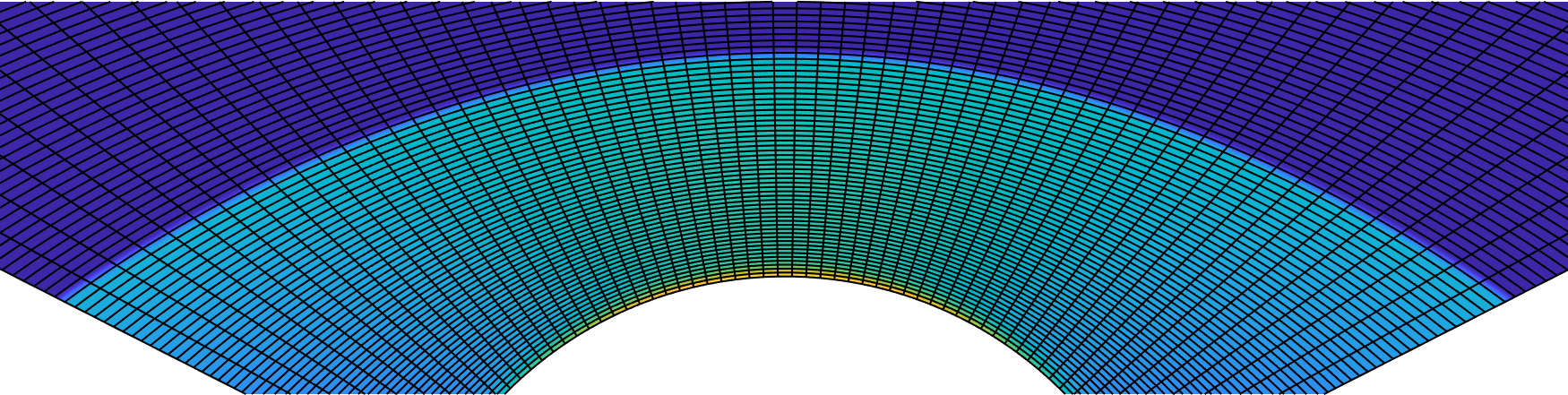}
		\caption{}
		\label{fig:cy_Uphy_1}
	\end{subfigure}
	\hfill
	\begin{subfigure}{0.49\textwidth}
		\centering
		\includegraphics[width=\textwidth,trim={0mm 0mm 0mm 0mm},clip]{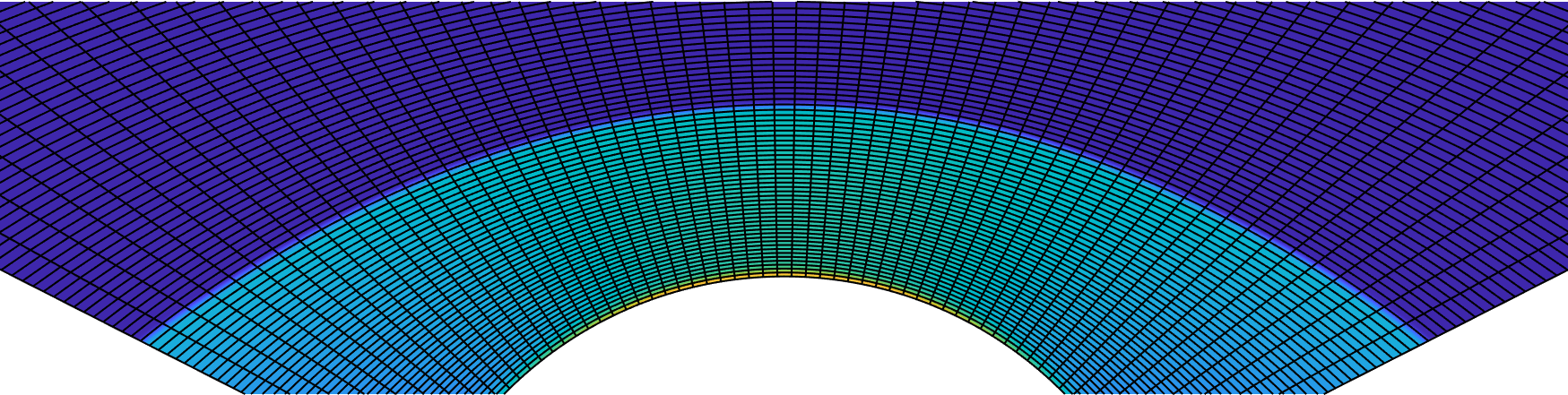}
		\caption{}
		\label{fig:cy_Uref_1}
	\end{subfigure}
	\\
	\begin{subfigure}{0.49\textwidth}
		\centering
		\includegraphics[width=\textwidth,trim={0mm 0mm 0mm 0mm},clip]{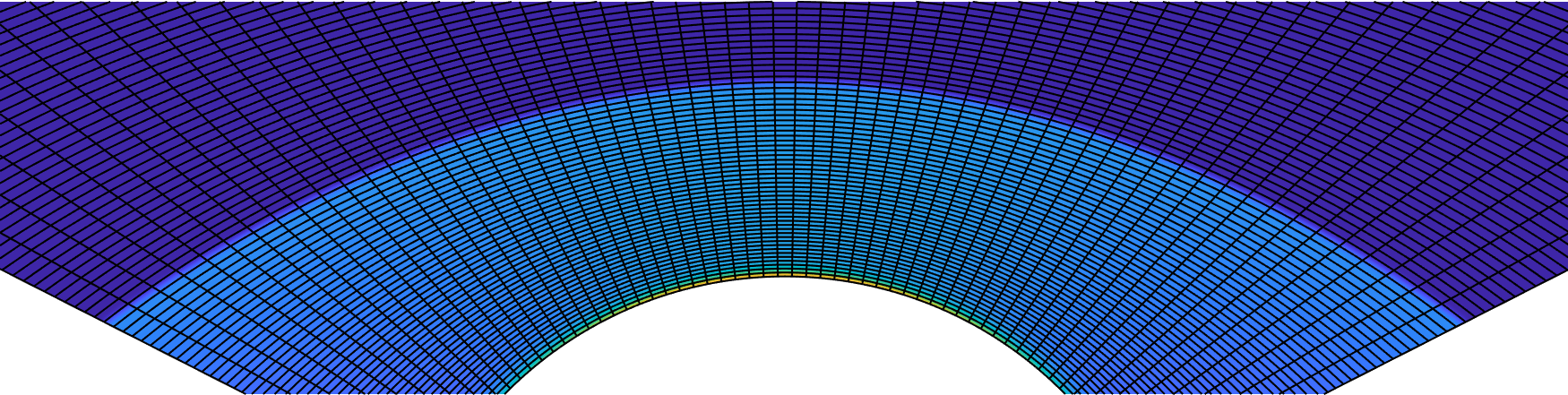}
		\caption{}
		\label{fig:cy_Uphy_2}
	\end{subfigure}
	\hfill
	\begin{subfigure}{0.49\textwidth}
		\centering
		\includegraphics[width=\textwidth,trim={0mm 0mm 0mm 0mm},clip]{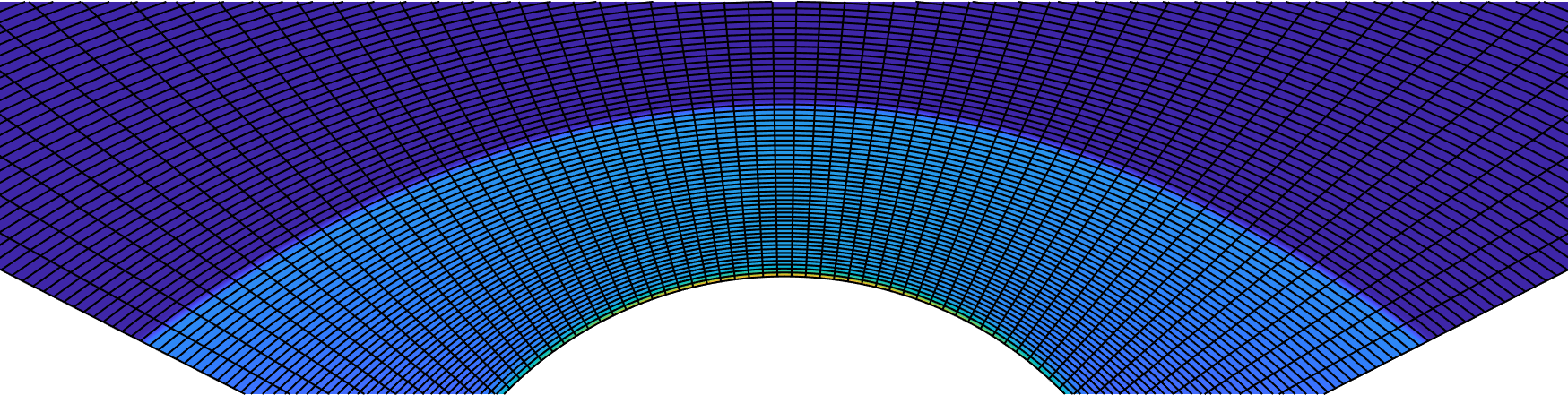}
		\caption{}
		\label{fig:cy_Uref_2}
	\end{subfigure}
	\\
	\begin{subfigure}{0.49\textwidth}
		\centering
		\includegraphics[width=\textwidth,trim={0mm 0mm 0mm 0mm},clip]{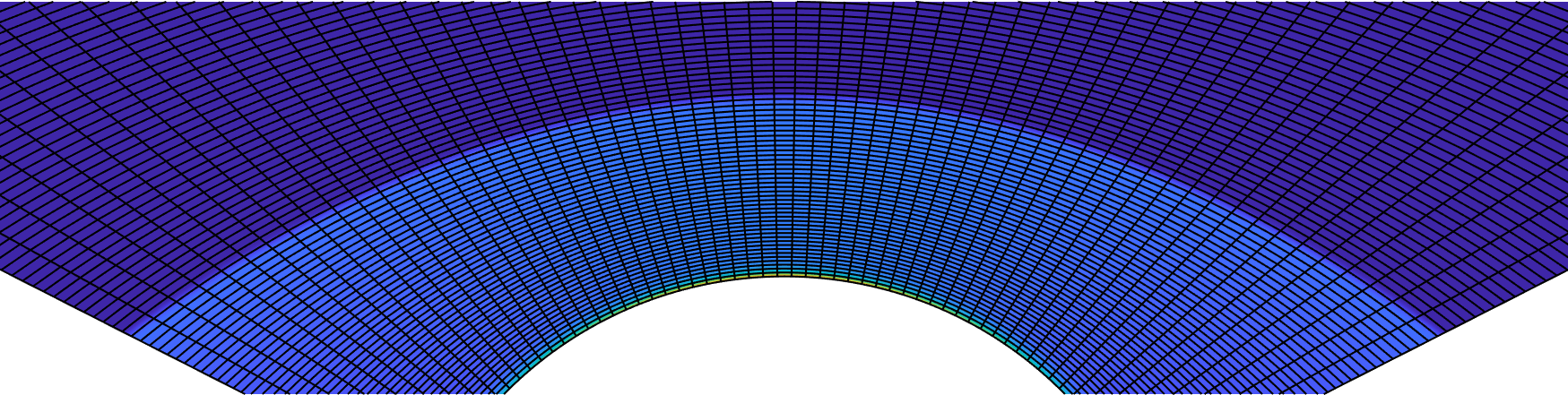}
		\caption{}
		\label{fig:cy_Uphy_3}
	\end{subfigure}
	\hfill
	\begin{subfigure}{0.49\textwidth}
		\centering
		\includegraphics[width=\textwidth,trim={0mm 0mm 0mm 0mm},clip]{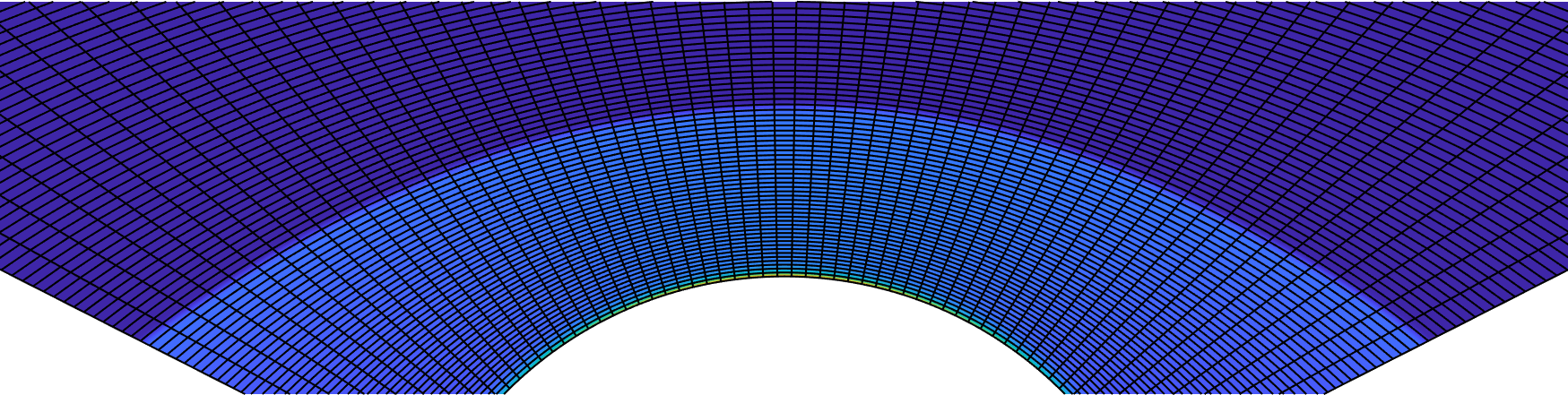}
		\caption{}
		\label{fig:cy_Uref_3}
	\end{subfigure}
	\\
	\begin{subfigure}{0.49\textwidth}
		\centering
		\includegraphics[width=\textwidth,trim={0mm 0mm 0mm 0mm},clip]{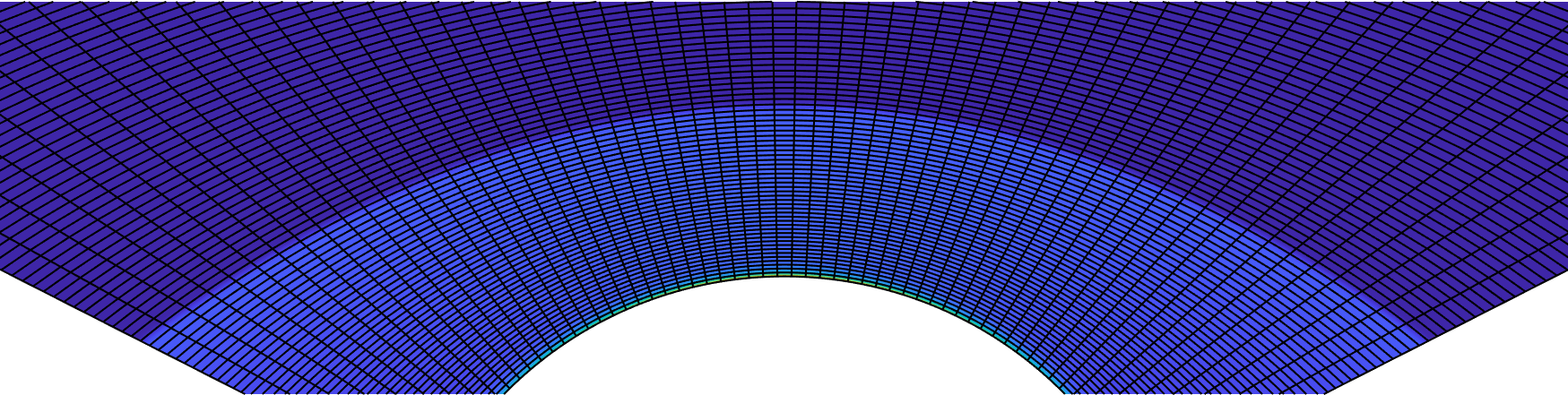}
		\caption{}
		\label{fig:cy_Uphy_4}
	\end{subfigure}
	\hfill
	\begin{subfigure}{0.49\textwidth}
		\centering
		\includegraphics[width=\textwidth,trim={0mm 0mm 0mm 0mm},clip]{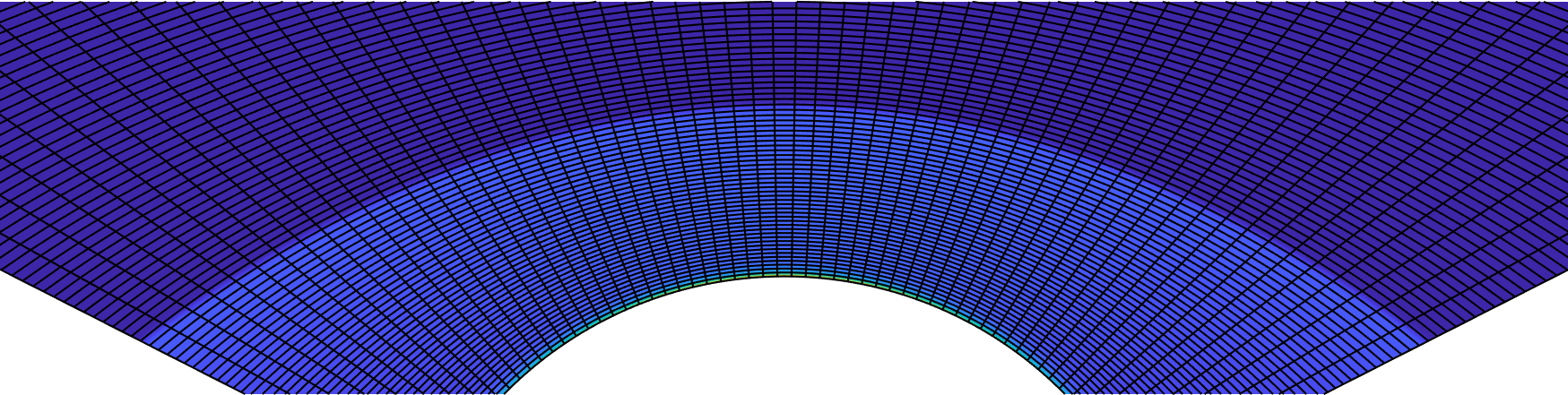}
		\caption{}
		\label{fig:cy_Uref_4}
	\end{subfigure}
	\\
	\begin{subfigure}{0.49\textwidth}
		\centering
		\includegraphics[width=\textwidth,trim={0mm 0mm 0mm 0mm},clip]{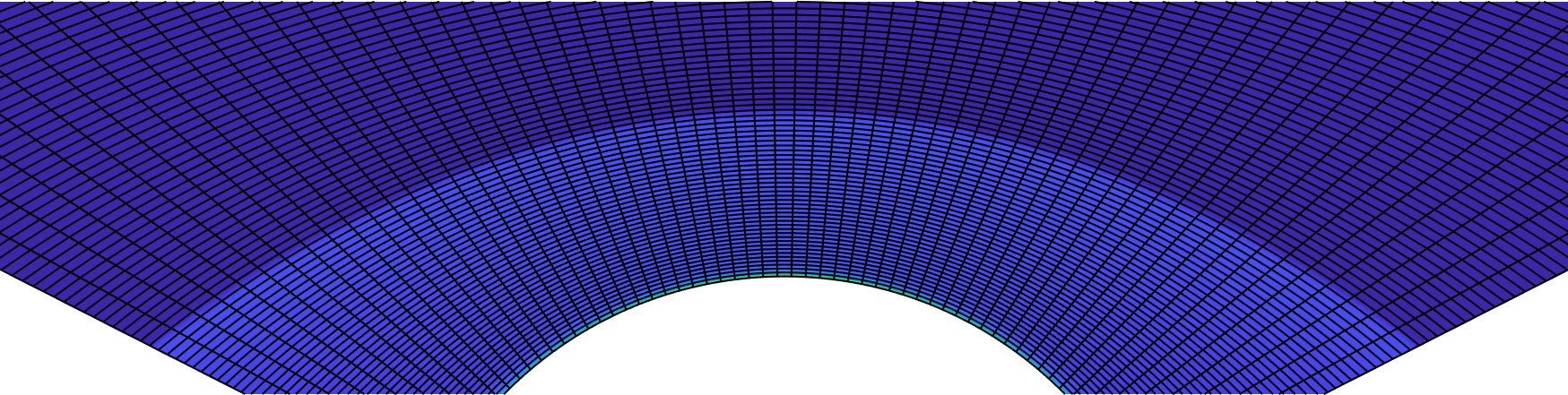}
		\caption{}
		\label{fig:cy_Uphy_5}
	\end{subfigure}
	\hfill
	\begin{subfigure}{0.49\textwidth}
		\centering
		\includegraphics[width=\textwidth,trim={0mm 0mm 0mm 0mm},clip]{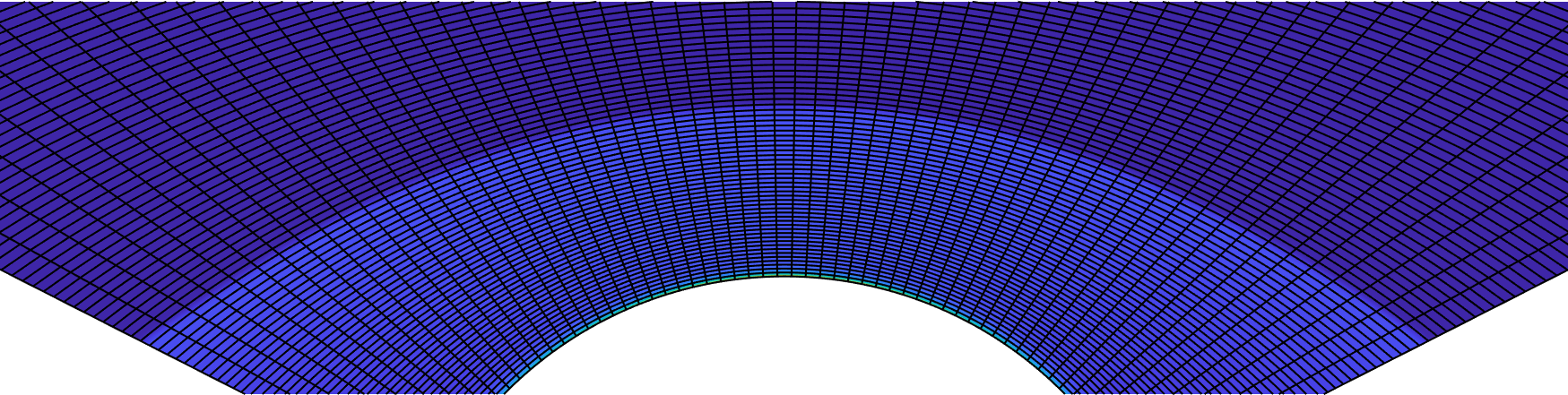}
		\caption{}
		\label{fig:cy_Uref_5}
	\end{subfigure}
	\caption{
		Hypersonic cylinder: Density snapshots before (left) and after (right) alignment for different values of Mach number \( M_\infty \in \{ 4, 5, 6, 7, 8\} \) (from top to bottom). The aligned snapshots are presented after interpolating the solutions into the reference domain.
	}
	\label{fig:cy_snap_surf}
\end{figure}

\begin{figure}[hbt!]
	\centering
	\begin{subfigure}{0.99\textwidth}
		\centering
		\includegraphics[width=\textwidth,trim={0mm 0mm 0mm 0mm},clip]{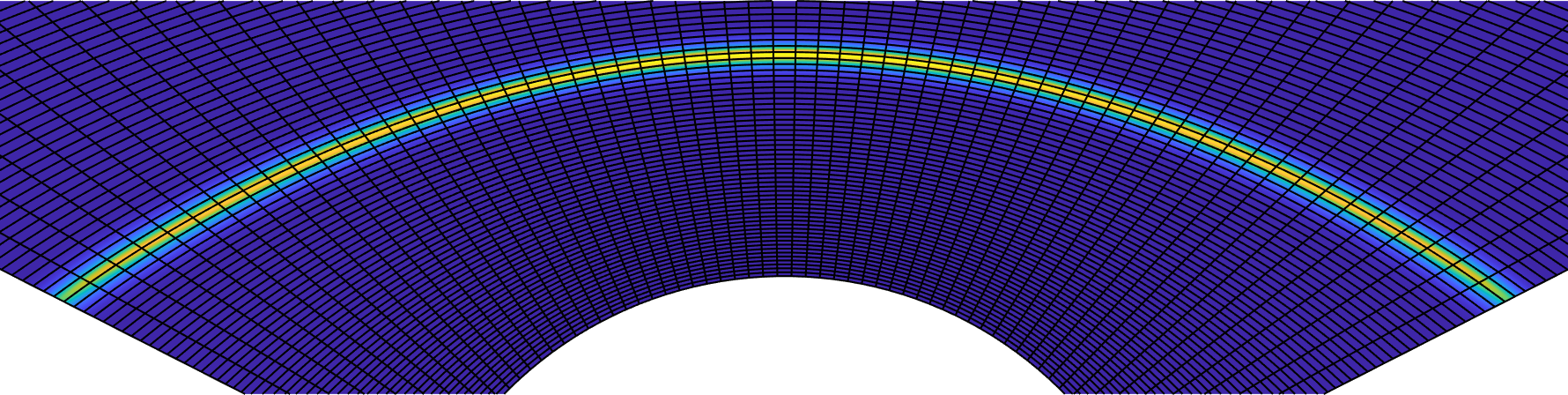}
		\caption{\( M_\infty = 4 \)}
		\label{fig:cy_av_M4}
	\end{subfigure}
	\\
	\begin{subfigure}{0.99\textwidth}
		\centering
		\includegraphics[width=\textwidth,trim={0mm 0mm 0mm 0mm},clip]{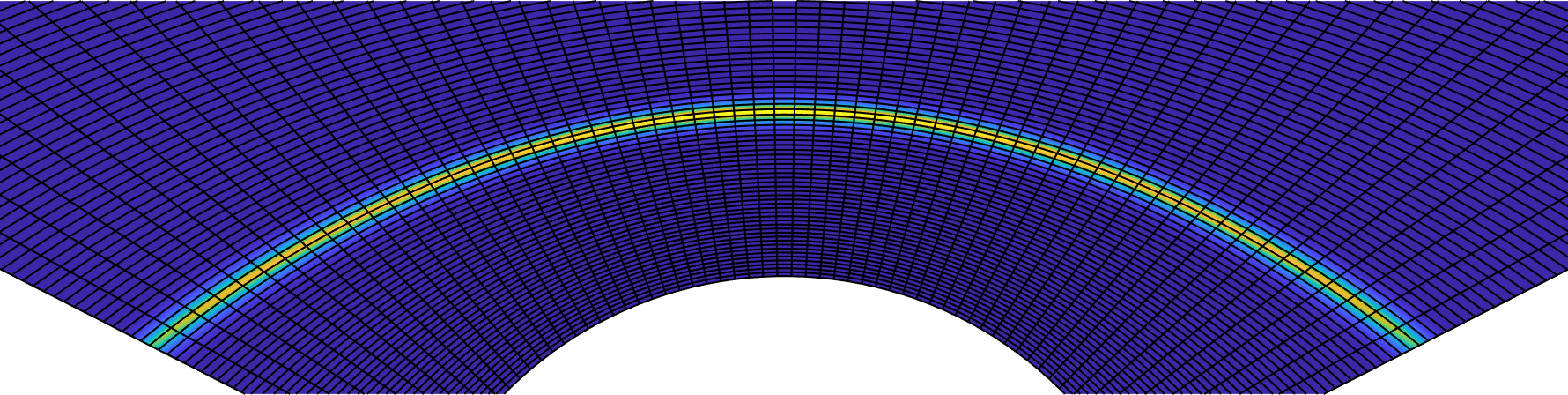}
		\caption{\( M_\infty = 8 \)}
		\label{fig:cy_av_M8}
	\end{subfigure}
	\caption{Hypersonic cylinder: Artificial viscosity field (\( \epsilon(x) \)).}
	\label{fig:cy_av}
\end{figure}

\begin{figure}[hbt!]
	\centering
	\input{figs/cylinder/control_points/cy_xc.tikz}
	\caption{Hypersonic cylinder: Placement of control points (\ref{line:cy_xc}) on the reference snapshot ($M_\infty = 7$).}
	\label{fig:cy_xc}
\end{figure}
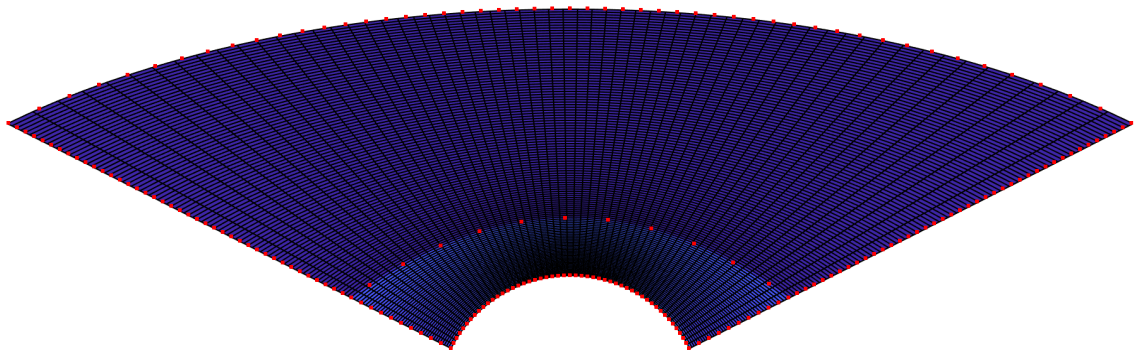

\begin{figure}[hbt!]
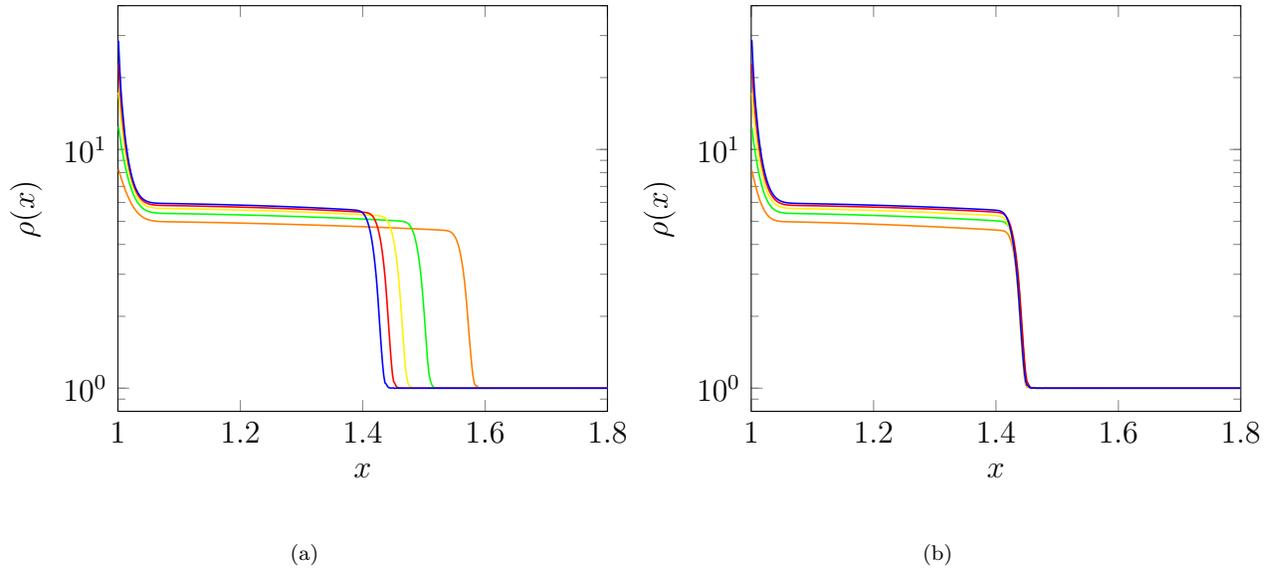

	\centering
	\begin{subfigure}{0.49\textwidth}
		\centering
		\input{figs/cylinder/slice/ref_and_phy/cy_phy_slice.tikz}
		\caption{}
		\label{fig:cy_phy_slice}
	\end{subfigure}
	\hfill
	\begin{subfigure}{0.49\textwidth}
		\centering
		\input{figs/cylinder/slice/ref_and_phy/cy_ref_slice.tikz}
		\caption{}
		\label{fig:cy_ref_slice}
	\end{subfigure}
	\caption{
		Hypersonic cylinder: The density along the stagnation line before (left) and after (right) alignment for different values of Mach number \( M_\infty \in \{ 4, 5, 6, 7, 8\} \) (from top to bottom).
	}
	\label{fig:cy_snap_slice}
\end{figure}

We solve the compressible Navier-Stokes equations coupled with the artificial viscosity equation for a Reynolds number of $Re = 1000$ and multiple values of the Mach number $M_\infty \in \{ 4, 5, 6, 7, 8 \}$ (Figure \ref{fig:cy_snap_surf}). As observed, a higher Mach number brings the bow shock closer to the cylinder. Figure \ref{fig:cy_av} illustrates the artificial viscosity state variable for two different Mach numbers, and as expected, it mirrors the bow shock location.
To perform landmark-based registration for this problem, we use the $M_\infty = 7$ solution as the reference. All other solutions are aligned to match the bow shock position of this reference. A total of $k = 11$ shock landmarks, along with additional boundary landmarks, are placed as shown in Fig. \ref{fig:cy_xc}. The radius of support for the RBF functions is set to $r = 100$. Figure \ref{fig:cy_snap_surf} also shows the density snapshots after alignment, while Fig. \ref{fig:cy_snap_slice} presents the density along the stagnation line for both aligned and unaligned solutions. Both figures demonstrate the effectivity of the landmark-based registration procedure as the bow shock location in the reference domain is essentially independent of the inflow Mach number.


\begin{figure}
	\centering
	\includegraphics[clip,width=.99\textwidth]{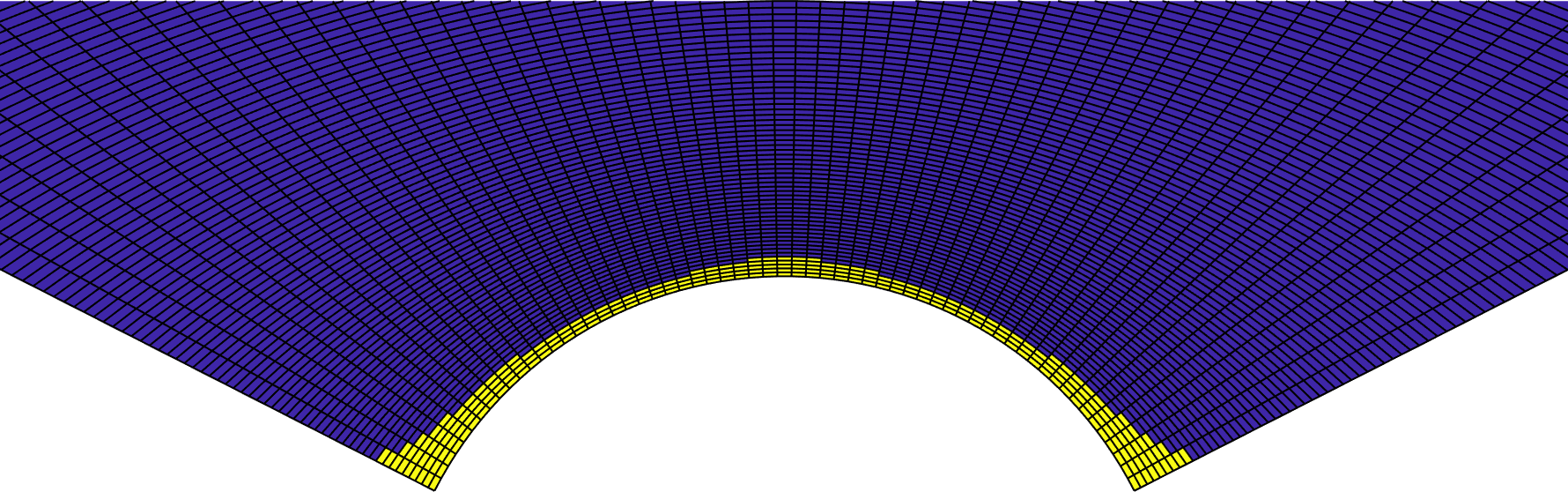}
	\caption{Hypersonic cylinder: Elements identified by the boundary layer sensor given by Eq. \eqref{eq:bl_sensor} are highlighted in yellow.}
	\label{fig:cy_bl_sensor}
\end{figure}

For this problem, we define the boundary decay function $f_w(x)$ to highlight the sharp gradient at the cylinder wall while suppressing a shock located further out. This function employs a radial logistic form, taking a value of 1 at the wall and decaying rapidly with increasing distance. Specifically, it is given by
\begin{equation}
	f_w(x) = 1 - \frac{1}{1 + e^{-k (r(x) - r_c)}},
\end{equation}
where $r(x) = \sqrt{x_1^2 + x_2^2}$ is the radial distance from the cylinder axis at point $x = (x_1, x_2)$, $r_c = 1.4$ marks the transition point of the logistic function, and $k = 10^2$ controls the steepness. This design ensures $f_w(x)$ remains near $1$ at the wall ($r = 1$), with a rapid drop-off around $r = 1.4$, preserving the wall gradient while filtering out shocks beyond this point. Figure \ref{fig:cy_bl_sensor} illustrates the elements selected by the boundary layer sensor (Eq. \eqref{eq:bl_sensor}) for the $M_\infty = 4$ solution (snapshot with thickest boundary layer). This is consistent with Fig. \ref{fig:cy_snap_slice}, showing that lower Mach numbers, at a fixed Reynolds number, yield thicker boundary layers due to weaker inertial forces and a slower transition to the no-slip condition.


For this problem, we compare the fixed mesh ROM and the implicit feature tracking ROM, while also examining the effect of an inner product matrix $\Theta$ that prioritizes boundary layer elements identified by the sensor in Eq. \eqref{eq:bl_sensor}, using two different weight values, $\omega = 1$ and $\omega = 10^4$.
For both ROMs, the reference state is defined as $Q_r = Q_\infty$, where $Q_\infty\in\Rbb^N$ is a vector containing the freestream conditions. This choice ensures that the inflow boundary conditions are naturally satisfied, a strategy commonly employed in reduced-order modeling.
Moreover, unlike the previous problem, all snapshots are used for both ROMs; however, the fixed mesh ROM utilizes unaligned snapshots, while the implicit feature tracking ROM employs aligned snapshots.

\begin{figure}[hbt!]
	\centering
	\input{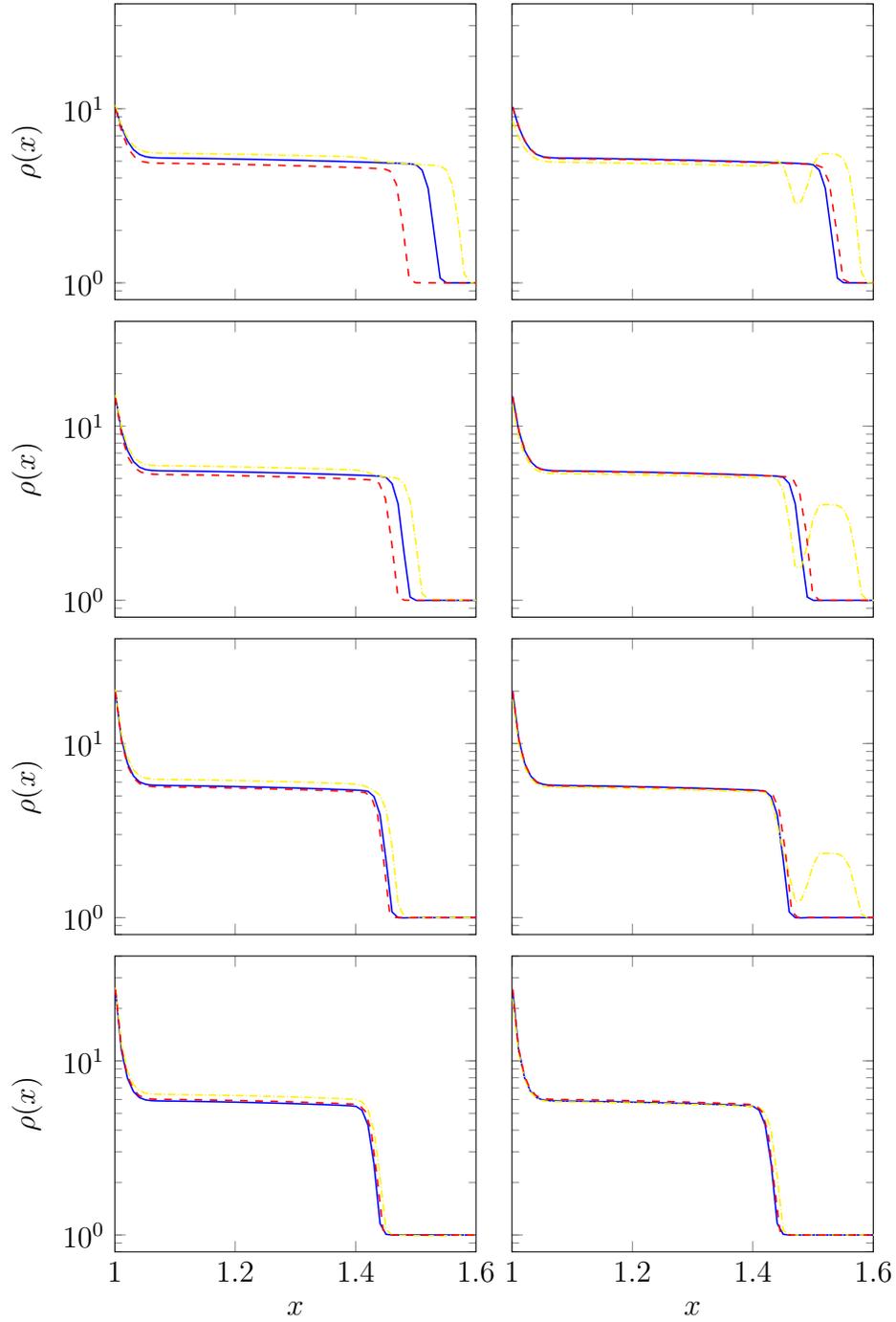}
	\caption{Hypersonic cylinder: The density along the stagnation line in the physical domain is shown for various Mach numbers, $M_\infty \in \{ 4.5, 5.5, 6.5, 7.5 \}$ (from top to bottom). Results are presented for the HDM (\ref{line:cy_hdm_45_slice_w0}), the fixed mesh ROM (\ref{line:cy_rom_45_slice_w0}), and the implicit feature tracking ROM (\ref{line:cy_romft_45_slice_w0}), for $\omega = 1$ (left) and $\omega = 10^4$ (right).
	}
	\label{fig:cy_rom_slice}
\end{figure}

Figure \ref{fig:cy_rom_slice} demonstrates that the implicit feature tracking ROM delivers superior performance compared to the fixed mesh ROM. Additionally, increasing $\omega$ in the boundary layer enhances accuracy, particularly for the implicit feature tracking ROM. In contrast, for the fixed mesh ROM, a higher $\omega$ improves accuracy near the cylinder but results in staircase-like solutions near the shock for most parameters.
It is important to emphasize that our efforts were focused solely on eliminating the convective behavior of the snapshots by aligning the shocks, without attempting to align the boundary layers. As seen in Figures \ref{fig:cy_snap_slice} and \ref{fig:cy_rom_slice}, the boundary layer varies not only in thickness but also in height, which could negatively impact the performance of our implicit feature ROM due to misalignment.

\section{Conclusions and future directions}
\label{sec:conclusions}

In this work, we present a landmark-based registration method for reduced-order models (ROMs) of convection-dominated problems, overcoming the limitations of traditional approaches that rely on extensive training data. By aligning dominant convective features through landmark-based registration and radial basis function interpolation, our method enables rapid error decay in a reference domain during the offline phase. In the online phase, it extends standard ROMs by incorporating admissible domain mappings, enhancing predictive performance.  

We assess the effectiveness of our approach on two convection-dominated viscous flow problems of increasing complexity, comparing it to the fixed-mesh ROM. Validation on the space-time Burgers' equation shows that while fixed-mesh ROMs struggle with accuracy, particularly outside the training window, the implicit feature tracking ROM consistently captures convective features, delivering accurate predictions with minimal training data. In the second problem, which is more challenging due to the presence of a boundary layer, our method continues to outperform the fixed-mesh ROM. However, because our approach focuses specifically on shock alignment, it does not account for the boundary layer when placing landmarks or aligning solutions, highlighting an avenue for further refinement.

A promising direction for future research is to integrate our method with dense image registration \cite{THIRION1998}. Dense image registration techniques are known for their sensitivity to the initial setup, including factors such as initial conditions and hyperparameters, and they typically benefit from an initial guess that is close to the optimal solution. By combining our shock alignment approach with dense registration, we could take advantage of both methods, using the shock alignment to capture sharp features while allowing dense registration to refine the alignment in smoother regions. This combination could improve overall accuracy and robustness, particularly in cases with complex features and smooth backgrounds.
Another direction for future research is to tackle more challenging problems. As problems become more complex, they typically involve a greater variety of features and potentially more distinctive points, such as corners or other sharp transitions. These features could be critical for improving alignment accuracy. In such cases, more advanced methods, as discussed in Section \ref{sec:landmark}, might offer additional benefits. These methods, which include sophisticated landmark detection and matching techniques, could potentially be integrated with our approach to enhance its performance. By combining the strengths of both, we could better handle the increased complexity and improve the alignment of even more intricate solutions.

\section*{Acknowledgments}

This material is based upon work supported by the Air Force Office of Scientific Research (AFOSR) under award numbers FA9550-20-1-0236 and FA9550-22-1-0004. The content of this publication does not necessarily reflect the position or policy of any of these supporters, and no official endorsement should be inferred.

\bibliographystyle{elsarticle-num} 
\bibliography{cas-refs}

\end{document}

%% file: preamble.pap.tex
\usepackage[margin=1in]{geometry}

\usepackage{graphicx}
\usepackage{pstool}
\usepackage{wrapfig}
\usepackage{caption}
\usepackage{subcaption}
\usepackage[section]{placeins} 

\usepackage{fancyhdr} 
\usepackage{lastpage} 
\usepackage{extramarks} 

\usepackage{xcolor} 
\usepackage{enumerate} 
\usepackage{paralist} 
\usepackage{amsmath, amsthm, amssymb, mathtools}
\usepackage{mathabx, pifont, stmaryrd} 
\usepackage[explicit]{titlesec} 
\usepackage{etoolbox} 
\usepackage{bibentry} 
\makeatletter\let\saved@bibitem\@bibitem\makeatother 
\usepackage[colorlinks, bookmarksopen, bookmarksnumbered,
            citecolor=red,urlcolor=red]{hyperref} 
\makeatletter\let\@bibitem\saved@bibitem\makeatother 

\usepackage{algpseudocode}

\algnewcommand{\algorithmicand}{\textbf{ and }}
\algnewcommand{\algorithmicor}{\textbf{ or }}
\algnewcommand{\OR}{\algorithmicor}
\algnewcommand{\AND}{\algorithmicand}
\algnewcommand{\var}{\texttt}

\newtheorem{remark}{Remark}

%% file: preamble.macros.tex
\usepackage{bm} 
\usepackage{amsfonts} 

\newcommand{\norm}[1]{\ensuremath{\left\| #1 \right\|}}

\newcommand{\pder}[2]{\ensuremath{\frac{\partial #1}{\partial #2}}}

\newcommand{\Dcal}{\ensuremath{\mathcal{D}}}

\newcommand{\Gcal}{\ensuremath{\mathcal{G}}}

\newcommand{\Xcal}{\ensuremath{\mathcal{X}}}

\newcommand{\Rbb}{\ensuremath{\mathbb{R} }}



%% file: preamble.tikz.tex
\usepackage{tikz}
\usepackage{pgfplots}
\usepackage{pgfplotstable, filecontents, booktabs}
\pgfplotsset{compat=1.9}

\usetikzlibrary{pgfplots.groupplots}
\usepgfplotslibrary{fillbetween}
\usetikzlibrary{calc,fit,matrix,arrows,automata,positioning,shapes}
\usetikzlibrary{arrows.meta}

\pgfplotsset{select coords between index/.style 2 args={
    x filter/.code={
        \ifnum\coordindex<#1\fi
        \ifnum\coordindex>#2\fi
    }
}}

\tikzset{
 invisible/.style={opacity=0},
 visible on/.style={alt={#1{}{invisible}}},
 alt/.code args={<#1>#2#3}{%
   \alt<#1>{\pgfkeysalso{#2}}{\pgfkeysalso{#3}}
 },
}


%% file: figs/methodology/bump1d_motivation/bump_1d_snapshots_xphys.tikz
\begin{tikzpicture}
\begin{axis}[
width=1.\textwidth,
xlabel={$x$},
ymax=1.2,
xmax=1,
ylabel={$y(x; c)$},
xmin=-1,
ymin=-0.2]
\addplot [black, , mark repeat=0, forget plot]
coordinates {
(-1.00000000e+00,  2.31952283e-16)
(-9.88000000e-01,  9.65003448e-16)
(-9.76000000e-01,  3.90077946e-15)
(-9.64000000e-01,  1.53202640e-14)
(-9.52000000e-01,  5.84619663e-14)
(-9.40000000e-01,  2.16756888e-13)
(-9.28000000e-01,  7.80844788e-13)
(-9.16000000e-01,  2.73305796e-12)
(-9.04000000e-01,  9.29448409e-12)
(-8.92000000e-01,  3.07110067e-11)
(-8.80000000e-01,  9.85950558e-11)
(-8.68000000e-01,  3.07544907e-10)
(-8.56000000e-01,  9.32082300e-10)
(-8.44000000e-01,  2.74468369e-09)
(-8.32000000e-01,  7.85276608e-09)
(-8.20000000e-01,  2.18295780e-08)
(-8.08000000e-01,  5.89603885e-08)
(-7.96000000e-01,  1.54727551e-07)
(-7.84000000e-01,  3.94518422e-07)
(-7.72000000e-01,  9.77370534e-07)
(-7.60000000e-01,  2.35257520e-06)
(-7.48000000e-01,  5.50199392e-06)
(-7.36000000e-01,  1.25022741e-05)
(-7.24000000e-01,  2.76026162e-05)
(-7.12000000e-01,  5.92111908e-05)
(-7.00000000e-01,  1.23409804e-04)
(-6.88000000e-01,  2.49912425e-04)
(-6.76000000e-01,  4.91720545e-04)
(-6.64000000e-01,  9.40028842e-04)
(-6.52000000e-01,  1.74604858e-03)
(-6.40000000e-01,  3.15111160e-03)
(-6.28000000e-01,  5.52539799e-03)
(-6.16000000e-01,  9.41359837e-03)
(-6.04000000e-01,  1.55826058e-02)
(-5.92000000e-01,  2.50620633e-02)
(-5.80000000e-01,  3.91638951e-02)
(-5.68000000e-01,  5.94630600e-02)
(-5.56000000e-01,  8.77204698e-02)
(-5.44000000e-01,  1.25732330e-01)
(-5.32000000e-01,  1.75099657e-01)
(-5.20000000e-01,  2.36927759e-01)
(-5.08000000e-01,  3.11486276e-01)
(-4.96000000e-01,  3.97881920e-01)
(-4.84000000e-01,  4.93812198e-01)
(-4.72000000e-01,  5.95472542e-01)
(-4.60000000e-01,  6.97676326e-01)
(-4.48000000e-01,  7.94215853e-01)
(-4.36000000e-01,  8.78446739e-01)
(-4.24000000e-01,  9.44027483e-01)
(-4.12000000e-01,  9.85703184e-01)
(-4.00000000e-01,  1.00000000e+00)
(-3.88000000e-01,  9.85703184e-01)
(-3.76000000e-01,  9.44027483e-01)
(-3.64000000e-01,  8.78446739e-01)
(-3.52000000e-01,  7.94215853e-01)
(-3.40000000e-01,  6.97676326e-01)
(-3.28000000e-01,  5.95472542e-01)
(-3.16000000e-01,  4.93812198e-01)
(-3.04000000e-01,  3.97881920e-01)
(-2.92000000e-01,  3.11486276e-01)
(-2.80000000e-01,  2.36927759e-01)
(-2.68000000e-01,  1.75099657e-01)
(-2.56000000e-01,  1.25732330e-01)
(-2.44000000e-01,  8.77204698e-02)
(-2.32000000e-01,  5.94630600e-02)
(-2.20000000e-01,  3.91638951e-02)
(-2.08000000e-01,  2.50620633e-02)
(-1.96000000e-01,  1.55826058e-02)
(-1.84000000e-01,  9.41359837e-03)
(-1.72000000e-01,  5.52539799e-03)
(-1.60000000e-01,  3.15111160e-03)
(-1.48000000e-01,  1.74604858e-03)
(-1.36000000e-01,  9.40028842e-04)
(-1.24000000e-01,  4.91720545e-04)
(-1.12000000e-01,  2.49912425e-04)
(-1.00000000e-01,  1.23409804e-04)
(-8.80000000e-02,  5.92111908e-05)
(-7.60000000e-02,  2.76026162e-05)
(-6.40000000e-02,  1.25022741e-05)
(-5.20000000e-02,  5.50199392e-06)
(-4.00000000e-02,  2.35257520e-06)
(-2.80000000e-02,  9.77370534e-07)
(-1.60000000e-02,  3.94518422e-07)
(-4.00000000e-03,  1.54727551e-07)
( 8.00000000e-03,  5.89603885e-08)
( 2.00000000e-02,  2.18295780e-08)
( 3.20000000e-02,  7.85276608e-09)
( 4.40000000e-02,  2.74468369e-09)
( 5.60000000e-02,  9.32082300e-10)
( 6.80000000e-02,  3.07544907e-10)
( 8.00000000e-02,  9.85950558e-11)
( 9.20000000e-02,  3.07110067e-11)
( 1.04000000e-01,  9.29448409e-12)
( 1.16000000e-01,  2.73305796e-12)
( 1.28000000e-01,  7.80844788e-13)
( 1.40000000e-01,  2.16756888e-13)
( 1.52000000e-01,  5.84619663e-14)
( 1.64000000e-01,  1.53202640e-14)
( 1.76000000e-01,  3.90077946e-15)
( 1.88000000e-01,  9.65003448e-16)
( 2.00000000e-01,  2.31952283e-16)
( 2.12000000e-01,  5.41702397e-17)
( 2.24000000e-01,  1.22917914e-17)
( 2.36000000e-01,  2.70995404e-18)
( 2.48000000e-01,  5.80498347e-19)
( 2.60000000e-01,  1.20818202e-19)
( 2.72000000e-01,  2.44318340e-20)
( 2.84000000e-01,  4.80034098e-21)
( 2.96000000e-01,  9.16390177e-22)
( 3.08000000e-01,  1.69973447e-22)
( 3.20000000e-01,  3.06319086e-23)
( 3.32000000e-01,  5.36363669e-24)
( 3.44000000e-01,  9.12508601e-25)
( 3.56000000e-01,  1.50836657e-25)
( 3.68000000e-01,  2.42252993e-26)
( 3.80000000e-01,  3.78027784e-27)
( 3.92000000e-01,  5.73153053e-28)
( 4.04000000e-01,  8.44325330e-29)
( 4.16000000e-01,  1.20848523e-29)
( 4.28000000e-01,  1.68060322e-30)
( 4.40000000e-01,  2.27081292e-31)
( 4.52000000e-01,  2.98119127e-32)
( 4.64000000e-01,  3.80268756e-33)
( 4.76000000e-01,  4.71285161e-34)
( 4.88000000e-01,  5.67504365e-35)
( 5.00000000e-01,  6.63967720e-36)
( 5.12000000e-01,  7.54774196e-37)
( 5.24000000e-01,  8.33641722e-38)
( 5.36000000e-01,  8.94610844e-39)
( 5.48000000e-01,  9.32784227e-40)
( 5.60000000e-01,  9.44975498e-41)
( 5.72000000e-01,  9.30148352e-42)
( 5.84000000e-01,  8.89561980e-43)
( 5.96000000e-01,  8.26594524e-44)
( 6.08000000e-01,  7.46278885e-45)
( 6.20000000e-01,  6.54639344e-46)
( 6.32000000e-01,  5.57950114e-47)
( 6.44000000e-01,  4.62041470e-48)
( 6.56000000e-01,  3.71756750e-49)
( 6.68000000e-01,  2.90622402e-50)
( 6.80000000e-01,  2.20745384e-51)
( 6.92000000e-01,  1.62909531e-52)
( 7.04000000e-01,  1.16813672e-53)
( 7.16000000e-01,  8.13829025e-55)
( 7.28000000e-01,  5.50890133e-56)
( 7.40000000e-01,  3.62317351e-57)
( 7.52000000e-01,  2.31529103e-58)
( 7.64000000e-01,  1.43752155e-59)
( 7.76000000e-01,  8.67192393e-61)
( 7.88000000e-01,  5.08286849e-62)
( 8.00000000e-01,  2.89464031e-63)
( 8.12000000e-01,  1.60166859e-64)
( 8.24000000e-01,  8.61079060e-66)
( 8.36000000e-01,  4.49785776e-67)
( 8.48000000e-01,  2.28276248e-68)
( 8.60000000e-01,  1.12566212e-69)
( 8.72000000e-01,  5.39321457e-71)
( 8.84000000e-01,  2.51061237e-72)
( 8.96000000e-01,  1.13554395e-73)
( 9.08000000e-01,  4.99022987e-75)
( 9.20000000e-01,  2.13073516e-76)
( 9.32000000e-01,  8.83956134e-78)
( 9.44000000e-01,  3.56306934e-79)
( 9.56000000e-01,  1.39543674e-80)
( 9.68000000e-01,  5.30992455e-82)
( 9.80000000e-01,  1.96317433e-83)
( 9.92000000e-01,  7.05215261e-85)};

\addplot [black, , mark repeat=0, forget plot]
coordinates {
(-1.00000000e+00,  1.60381089e-28)
(-9.88000000e-01,  1.07831271e-27)
(-9.76000000e-01,  7.04415018e-27)
(-9.64000000e-01,  4.47100143e-26)
(-9.52000000e-01,  2.75723213e-25)
(-9.40000000e-01,  1.65209178e-24)
(-9.28000000e-01,  9.61805664e-24)
(-9.16000000e-01,  5.44042476e-23)
(-9.04000000e-01,  2.98999599e-22)
(-8.92000000e-01,  1.59661690e-21)
(-8.80000000e-01,  8.28367701e-21)
(-8.68000000e-01,  4.17578289e-20)
(-8.56000000e-01,  2.04524332e-19)
(-8.44000000e-01,  9.73294738e-19)
(-8.32000000e-01,  4.50024425e-18)
(-8.20000000e-01,  2.02171585e-17)
(-8.08000000e-01,  8.82462916e-17)
(-7.96000000e-01,  3.74252863e-16)
(-7.84000000e-01,  1.54214822e-15)
(-7.72000000e-01,  6.17418194e-15)
(-7.60000000e-01,  2.40173478e-14)
(-7.48000000e-01,  9.07743172e-14)
(-7.36000000e-01,  3.33344468e-13)
(-7.24000000e-01,  1.18936686e-12)
(-7.12000000e-01,  4.12316518e-12)
(-7.00000000e-01,  1.38879439e-11)
(-6.88000000e-01,  4.54503770e-11)
(-6.76000000e-01,  1.44520464e-10)
(-6.64000000e-01,  4.46491836e-10)
(-6.52000000e-01,  1.34026296e-09)
(-6.40000000e-01,  3.90893843e-09)
(-6.28000000e-01,  1.10769441e-08)
(-6.16000000e-01,  3.04981442e-08)
(-6.04000000e-01,  8.15866656e-08)
(-5.92000000e-01,  2.12059275e-07)
(-5.80000000e-01,  5.35534780e-07)
(-5.68000000e-01,  1.31404551e-06)
(-5.56000000e-01,  3.13274795e-06)
(-5.44000000e-01,  7.25659389e-06)
(-5.32000000e-01,  1.63317417e-05)
(-5.20000000e-01,  3.57128496e-05)
(-5.08000000e-01,  7.58767686e-05)
(-4.96000000e-01,  1.56633794e-04)
(-4.84000000e-01,  3.14162560e-04)
(-4.72000000e-01,  6.12231548e-04)
(-4.60000000e-01,  1.15922917e-03)
(-4.48000000e-01,  2.13262855e-03)
(-4.36000000e-01,  3.81200493e-03)
(-4.24000000e-01,  6.62039661e-03)
(-4.12000000e-01,  1.11713814e-02)
(-4.00000000e-01,  1.83156389e-02)
(-3.88000000e-01,  2.91762575e-02)
(-3.76000000e-01,  4.51574503e-02)
(-3.64000000e-01,  6.79080972e-02)
(-3.52000000e-01,  9.92215550e-02)
(-3.40000000e-01,  1.40858421e-01)
(-3.28000000e-01,  1.94290659e-01)
(-3.16000000e-01,  2.60383431e-01)
(-3.04000000e-01,  3.39052607e-01)
(-2.92000000e-01,  4.28956399e-01)
(-2.80000000e-01,  5.27292424e-01)
(-2.68000000e-01,  6.29770381e-01)
(-2.56000000e-01,  7.30811294e-01)
(-2.44000000e-01,  8.23987433e-01)
(-2.32000000e-01,  9.02668412e-01)
(-2.20000000e-01,  9.60789439e-01)
(-2.08000000e-01,  9.93620436e-01)
(-1.96000000e-01,  9.98401279e-01)
(-1.84000000e-01,  9.74724902e-01)
(-1.72000000e-01,  9.24594515e-01)
(-1.60000000e-01,  8.52143789e-01)
(-1.48000000e-01,  7.63074204e-01)
(-1.36000000e-01,  6.63915763e-01)
(-1.24000000e-01,  5.61243736e-01)
(-1.12000000e-01,  4.60980286e-01)
(-1.00000000e-01,  3.67879441e-01)
(-8.80000000e-02,  2.85246945e-01)
(-7.60000000e-02,  2.14896234e-01)
(-6.40000000e-02,  1.57300074e-01)
(-5.20000000e-02,  1.11871991e-01)
(-4.00000000e-02,  7.73047404e-02)
(-2.80000000e-02,  5.19018938e-02)
(-1.60000000e-02,  3.38573219e-02)
(-4.00000000e-03,  2.14592391e-02)
( 8.00000000e-03,  1.32150361e-02)
( 2.00000000e-02,  7.90705405e-03)
( 3.20000000e-02,  4.59677642e-03)
( 4.40000000e-02,  2.59647636e-03)
( 5.60000000e-02,  1.42497644e-03)
( 6.80000000e-02,  7.59842031e-04)
( 8.00000000e-02,  3.93669041e-04)
( 9.20000000e-02,  1.98167089e-04)
( 1.04000000e-01,  9.69223903e-05)
( 1.16000000e-01,  4.60584185e-05)
( 1.28000000e-01,  2.12660214e-05)
( 1.40000000e-01,  9.54016287e-06)
( 1.52000000e-01,  4.15831755e-06)
( 1.64000000e-01,  1.76105059e-06)
( 1.76000000e-01,  7.24633375e-07)
( 1.88000000e-01,  2.89705774e-07)
( 2.00000000e-01,  1.12535175e-07)
( 2.12000000e-01,  4.24728835e-08)
( 2.24000000e-01,  1.55749809e-08)
( 2.36000000e-01,  5.54926694e-09)
( 2.48000000e-01,  1.92103833e-09)
( 2.60000000e-01,  6.46143177e-10)
( 2.72000000e-01,  2.11161051e-10)
( 2.84000000e-01,  6.70488271e-11)
( 2.96000000e-01,  2.06852548e-11)
( 3.08000000e-01,  6.20044532e-12)
( 3.20000000e-01,  1.80583144e-12)
( 3.32000000e-01,  5.11003477e-13)
( 3.44000000e-01,  1.40495631e-13)
( 3.56000000e-01,  3.75313429e-14)
( 3.68000000e-01,  9.74131769e-15)
( 3.80000000e-01,  2.45659537e-15)
( 3.92000000e-01,  6.01924290e-16)
( 4.04000000e-01,  1.43298759e-16)
( 4.16000000e-01,  3.31463191e-17)
( 4.28000000e-01,  7.44938645e-18)
( 4.40000000e-01,  1.62666462e-18)
( 4.52000000e-01,  3.45118195e-19)
( 4.64000000e-01,  7.11426455e-20)
( 4.76000000e-01,  1.42490032e-20)
( 4.88000000e-01,  2.77288128e-21)
( 5.00000000e-01,  5.24288566e-22)
( 5.12000000e-01,  9.63167493e-23)
( 5.24000000e-01,  1.71919669e-23)
( 5.36000000e-01,  2.98154682e-24)
( 5.48000000e-01,  5.02400450e-25)
( 5.60000000e-01,  8.22528065e-26)
( 5.72000000e-01,  1.30840968e-26)
( 5.84000000e-01,  2.02222322e-27)
( 5.96000000e-01,  3.03673402e-28)
( 6.08000000e-01,  4.43074480e-29)
( 6.20000000e-01,  6.28114815e-30)
( 6.32000000e-01,  8.65154470e-31)
( 6.44000000e-01,  1.15781879e-31)
( 6.56000000e-01,  1.50549685e-32)
( 6.68000000e-01,  1.90200400e-33)
( 6.80000000e-01,  2.33472278e-34)
( 6.92000000e-01,  2.78452765e-35)
( 7.04000000e-01,  3.22671099e-36)
( 7.16000000e-01,  3.63296254e-37)
( 7.28000000e-01,  3.97424000e-38)
( 7.40000000e-01,  4.22415241e-39)
( 7.52000000e-01,  4.36231866e-40)
( 7.64000000e-01,  4.37711054e-41)
( 7.76000000e-01,  4.26726841e-42)
( 7.88000000e-01,  4.04207833e-43)
( 8.00000000e-01,  3.72007598e-44)
( 8.12000000e-01,  3.32652822e-45)
( 8.24000000e-01,  2.89016698e-46)
( 8.36000000e-01,  2.43975926e-47)
( 8.48000000e-01,  2.00107486e-48)
( 8.60000000e-01,  1.59467437e-49)
( 8.72000000e-01,  1.23473287e-50)
( 8.84000000e-01,  9.28894363e-52)
( 8.96000000e-01,  6.78972153e-53)
( 9.08000000e-01,  4.82203010e-54)
( 9.20000000e-01,  3.32736305e-55)
( 9.32000000e-01,  2.23081102e-56)
( 9.44000000e-01,  1.45317425e-57)
( 9.56000000e-01,  9.19739679e-59)
( 9.68000000e-01,  5.65593581e-60)
( 9.80000000e-01,  3.37937463e-61)
( 9.92000000e-01,  1.96182618e-62)};

\addplot [black, , mark repeat=0, forget plot]
coordinates {
(-1.00000000e+00,  3.72007598e-44)
(-9.88000000e-01,  4.04207833e-43)
(-9.76000000e-01,  4.26726841e-42)
(-9.64000000e-01,  4.37711054e-41)
(-9.52000000e-01,  4.36231866e-40)
(-9.40000000e-01,  4.22415241e-39)
(-9.28000000e-01,  3.97424000e-38)
(-9.16000000e-01,  3.63296254e-37)
(-9.04000000e-01,  3.22671099e-36)
(-8.92000000e-01,  2.78452765e-35)
(-8.80000000e-01,  2.33472278e-34)
(-8.68000000e-01,  1.90200400e-33)
(-8.56000000e-01,  1.50549685e-32)
(-8.44000000e-01,  1.15781879e-31)
(-8.32000000e-01,  8.65154470e-31)
(-8.20000000e-01,  6.28114815e-30)
(-8.08000000e-01,  4.43074480e-29)
(-7.96000000e-01,  3.03673402e-28)
(-7.84000000e-01,  2.02222322e-27)
(-7.72000000e-01,  1.30840968e-26)
(-7.60000000e-01,  8.22528065e-26)
(-7.48000000e-01,  5.02400450e-25)
(-7.36000000e-01,  2.98154682e-24)
(-7.24000000e-01,  1.71919669e-23)
(-7.12000000e-01,  9.63167493e-23)
(-7.00000000e-01,  5.24288566e-22)
(-6.88000000e-01,  2.77288128e-21)
(-6.76000000e-01,  1.42490032e-20)
(-6.64000000e-01,  7.11426455e-20)
(-6.52000000e-01,  3.45118195e-19)
(-6.40000000e-01,  1.62666462e-18)
(-6.28000000e-01,  7.44938645e-18)
(-6.16000000e-01,  3.31463191e-17)
(-6.04000000e-01,  1.43298759e-16)
(-5.92000000e-01,  6.01924290e-16)
(-5.80000000e-01,  2.45659537e-15)
(-5.68000000e-01,  9.74131769e-15)
(-5.56000000e-01,  3.75313429e-14)
(-5.44000000e-01,  1.40495631e-13)
(-5.32000000e-01,  5.11003477e-13)
(-5.20000000e-01,  1.80583144e-12)
(-5.08000000e-01,  6.20044532e-12)
(-4.96000000e-01,  2.06852548e-11)
(-4.84000000e-01,  6.70488271e-11)
(-4.72000000e-01,  2.11161051e-10)
(-4.60000000e-01,  6.46143177e-10)
(-4.48000000e-01,  1.92103833e-09)
(-4.36000000e-01,  5.54926694e-09)
(-4.24000000e-01,  1.55749809e-08)
(-4.12000000e-01,  4.24728835e-08)
(-4.00000000e-01,  1.12535175e-07)
(-3.88000000e-01,  2.89705774e-07)
(-3.76000000e-01,  7.24633375e-07)
(-3.64000000e-01,  1.76105059e-06)
(-3.52000000e-01,  4.15831755e-06)
(-3.40000000e-01,  9.54016287e-06)
(-3.28000000e-01,  2.12660214e-05)
(-3.16000000e-01,  4.60584185e-05)
(-3.04000000e-01,  9.69223903e-05)
(-2.92000000e-01,  1.98167089e-04)
(-2.80000000e-01,  3.93669041e-04)
(-2.68000000e-01,  7.59842031e-04)
(-2.56000000e-01,  1.42497644e-03)
(-2.44000000e-01,  2.59647636e-03)
(-2.32000000e-01,  4.59677642e-03)
(-2.20000000e-01,  7.90705405e-03)
(-2.08000000e-01,  1.32150361e-02)
(-1.96000000e-01,  2.14592391e-02)
(-1.84000000e-01,  3.38573219e-02)
(-1.72000000e-01,  5.19018938e-02)
(-1.60000000e-01,  7.73047404e-02)
(-1.48000000e-01,  1.11871991e-01)
(-1.36000000e-01,  1.57300074e-01)
(-1.24000000e-01,  2.14896234e-01)
(-1.12000000e-01,  2.85246945e-01)
(-1.00000000e-01,  3.67879441e-01)
(-8.80000000e-02,  4.60980286e-01)
(-7.60000000e-02,  5.61243736e-01)
(-6.40000000e-02,  6.63915763e-01)
(-5.20000000e-02,  7.63074204e-01)
(-4.00000000e-02,  8.52143789e-01)
(-2.80000000e-02,  9.24594515e-01)
(-1.60000000e-02,  9.74724902e-01)
(-4.00000000e-03,  9.98401279e-01)
( 8.00000000e-03,  9.93620436e-01)
( 2.00000000e-02,  9.60789439e-01)
( 3.20000000e-02,  9.02668412e-01)
( 4.40000000e-02,  8.23987433e-01)
( 5.60000000e-02,  7.30811294e-01)
( 6.80000000e-02,  6.29770381e-01)
( 8.00000000e-02,  5.27292424e-01)
( 9.20000000e-02,  4.28956399e-01)
( 1.04000000e-01,  3.39052607e-01)
( 1.16000000e-01,  2.60383431e-01)
( 1.28000000e-01,  1.94290659e-01)
( 1.40000000e-01,  1.40858421e-01)
( 1.52000000e-01,  9.92215550e-02)
( 1.64000000e-01,  6.79080972e-02)
( 1.76000000e-01,  4.51574503e-02)
( 1.88000000e-01,  2.91762575e-02)
( 2.00000000e-01,  1.83156389e-02)
( 2.12000000e-01,  1.11713814e-02)
( 2.24000000e-01,  6.62039661e-03)
( 2.36000000e-01,  3.81200493e-03)
( 2.48000000e-01,  2.13262855e-03)
( 2.60000000e-01,  1.15922917e-03)
( 2.72000000e-01,  6.12231548e-04)
( 2.84000000e-01,  3.14162560e-04)
( 2.96000000e-01,  1.56633794e-04)
( 3.08000000e-01,  7.58767686e-05)
( 3.20000000e-01,  3.57128496e-05)
( 3.32000000e-01,  1.63317417e-05)
( 3.44000000e-01,  7.25659389e-06)
( 3.56000000e-01,  3.13274795e-06)
( 3.68000000e-01,  1.31404551e-06)
( 3.80000000e-01,  5.35534780e-07)
( 3.92000000e-01,  2.12059275e-07)
( 4.04000000e-01,  8.15866656e-08)
( 4.16000000e-01,  3.04981442e-08)
( 4.28000000e-01,  1.10769441e-08)
( 4.40000000e-01,  3.90893843e-09)
( 4.52000000e-01,  1.34026296e-09)
( 4.64000000e-01,  4.46491836e-10)
( 4.76000000e-01,  1.44520464e-10)
( 4.88000000e-01,  4.54503770e-11)
( 5.00000000e-01,  1.38879439e-11)
( 5.12000000e-01,  4.12316518e-12)
( 5.24000000e-01,  1.18936686e-12)
( 5.36000000e-01,  3.33344468e-13)
( 5.48000000e-01,  9.07743172e-14)
( 5.60000000e-01,  2.40173478e-14)
( 5.72000000e-01,  6.17418194e-15)
( 5.84000000e-01,  1.54214822e-15)
( 5.96000000e-01,  3.74252863e-16)
( 6.08000000e-01,  8.82462916e-17)
( 6.20000000e-01,  2.02171585e-17)
( 6.32000000e-01,  4.50024425e-18)
( 6.44000000e-01,  9.73294738e-19)
( 6.56000000e-01,  2.04524332e-19)
( 6.68000000e-01,  4.17578289e-20)
( 6.80000000e-01,  8.28367701e-21)
( 6.92000000e-01,  1.59661690e-21)
( 7.04000000e-01,  2.98999599e-22)
( 7.16000000e-01,  5.44042476e-23)
( 7.28000000e-01,  9.61805664e-24)
( 7.40000000e-01,  1.65209178e-24)
( 7.52000000e-01,  2.75723213e-25)
( 7.64000000e-01,  4.47100143e-26)
( 7.76000000e-01,  7.04415018e-27)
( 7.88000000e-01,  1.07831271e-27)
( 8.00000000e-01,  1.60381089e-28)
( 8.12000000e-01,  2.31768234e-29)
( 8.24000000e-01,  3.25422057e-30)
( 8.36000000e-01,  4.43948280e-31)
( 8.48000000e-01,  5.88450706e-32)
( 8.60000000e-01,  7.57844527e-33)
( 8.72000000e-01,  9.48292831e-34)
( 8.84000000e-01,  1.15291459e-34)
( 8.96000000e-01,  1.36189656e-35)
( 9.08000000e-01,  1.56308799e-36)
( 9.20000000e-01,  1.74307090e-37)
( 9.32000000e-01,  1.88859569e-38)
( 9.44000000e-01,  1.98817797e-39)
( 9.56000000e-01,  2.03359207e-40)
( 9.68000000e-01,  2.02099267e-41)
( 9.80000000e-01,  1.95145238e-42)
( 9.92000000e-01,  1.83081093e-43)};

\addplot [black, , mark repeat=0, forget plot]
coordinates {
(-1.00000000e+00,  2.89464031e-63)
(-9.88000000e-01,  5.08286849e-62)
(-9.76000000e-01,  8.67192393e-61)
(-9.64000000e-01,  1.43752155e-59)
(-9.52000000e-01,  2.31529103e-58)
(-9.40000000e-01,  3.62317351e-57)
(-9.28000000e-01,  5.50890133e-56)
(-9.16000000e-01,  8.13829025e-55)
(-9.04000000e-01,  1.16813672e-53)
(-8.92000000e-01,  1.62909531e-52)
(-8.80000000e-01,  2.20745384e-51)
(-8.68000000e-01,  2.90622402e-50)
(-8.56000000e-01,  3.71756750e-49)
(-8.44000000e-01,  4.62041470e-48)
(-8.32000000e-01,  5.57950114e-47)
(-8.20000000e-01,  6.54639344e-46)
(-8.08000000e-01,  7.46278885e-45)
(-7.96000000e-01,  8.26594524e-44)
(-7.84000000e-01,  8.89561980e-43)
(-7.72000000e-01,  9.30148352e-42)
(-7.60000000e-01,  9.44975498e-41)
(-7.48000000e-01,  9.32784227e-40)
(-7.36000000e-01,  8.94610844e-39)
(-7.24000000e-01,  8.33641722e-38)
(-7.12000000e-01,  7.54774196e-37)
(-7.00000000e-01,  6.63967720e-36)
(-6.88000000e-01,  5.67504365e-35)
(-6.76000000e-01,  4.71285161e-34)
(-6.64000000e-01,  3.80268756e-33)
(-6.52000000e-01,  2.98119127e-32)
(-6.40000000e-01,  2.27081292e-31)
(-6.28000000e-01,  1.68060322e-30)
(-6.16000000e-01,  1.20848523e-29)
(-6.04000000e-01,  8.44325330e-29)
(-5.92000000e-01,  5.73153053e-28)
(-5.80000000e-01,  3.78027784e-27)
(-5.68000000e-01,  2.42252993e-26)
(-5.56000000e-01,  1.50836657e-25)
(-5.44000000e-01,  9.12508601e-25)
(-5.32000000e-01,  5.36363669e-24)
(-5.20000000e-01,  3.06319086e-23)
(-5.08000000e-01,  1.69973447e-22)
(-4.96000000e-01,  9.16390177e-22)
(-4.84000000e-01,  4.80034098e-21)
(-4.72000000e-01,  2.44318340e-20)
(-4.60000000e-01,  1.20818202e-19)
(-4.48000000e-01,  5.80498347e-19)
(-4.36000000e-01,  2.70995404e-18)
(-4.24000000e-01,  1.22917914e-17)
(-4.12000000e-01,  5.41702397e-17)
(-4.00000000e-01,  2.31952283e-16)
(-3.88000000e-01,  9.65003448e-16)
(-3.76000000e-01,  3.90077946e-15)
(-3.64000000e-01,  1.53202640e-14)
(-3.52000000e-01,  5.84619663e-14)
(-3.40000000e-01,  2.16756888e-13)
(-3.28000000e-01,  7.80844788e-13)
(-3.16000000e-01,  2.73305796e-12)
(-3.04000000e-01,  9.29448409e-12)
(-2.92000000e-01,  3.07110067e-11)
(-2.80000000e-01,  9.85950558e-11)
(-2.68000000e-01,  3.07544907e-10)
(-2.56000000e-01,  9.32082300e-10)
(-2.44000000e-01,  2.74468369e-09)
(-2.32000000e-01,  7.85276608e-09)
(-2.20000000e-01,  2.18295780e-08)
(-2.08000000e-01,  5.89603885e-08)
(-1.96000000e-01,  1.54727551e-07)
(-1.84000000e-01,  3.94518422e-07)
(-1.72000000e-01,  9.77370534e-07)
(-1.60000000e-01,  2.35257520e-06)
(-1.48000000e-01,  5.50199392e-06)
(-1.36000000e-01,  1.25022741e-05)
(-1.24000000e-01,  2.76026162e-05)
(-1.12000000e-01,  5.92111908e-05)
(-1.00000000e-01,  1.23409804e-04)
(-8.80000000e-02,  2.49912425e-04)
(-7.60000000e-02,  4.91720545e-04)
(-6.40000000e-02,  9.40028842e-04)
(-5.20000000e-02,  1.74604858e-03)
(-4.00000000e-02,  3.15111160e-03)
(-2.80000000e-02,  5.52539799e-03)
(-1.60000000e-02,  9.41359837e-03)
(-4.00000000e-03,  1.55826058e-02)
( 8.00000000e-03,  2.50620633e-02)
( 2.00000000e-02,  3.91638951e-02)
( 3.20000000e-02,  5.94630600e-02)
( 4.40000000e-02,  8.77204698e-02)
( 5.60000000e-02,  1.25732330e-01)
( 6.80000000e-02,  1.75099657e-01)
( 8.00000000e-02,  2.36927759e-01)
( 9.20000000e-02,  3.11486276e-01)
( 1.04000000e-01,  3.97881920e-01)
( 1.16000000e-01,  4.93812198e-01)
( 1.28000000e-01,  5.95472542e-01)
( 1.40000000e-01,  6.97676326e-01)
( 1.52000000e-01,  7.94215853e-01)
( 1.64000000e-01,  8.78446739e-01)
( 1.76000000e-01,  9.44027483e-01)
( 1.88000000e-01,  9.85703184e-01)
( 2.00000000e-01,  1.00000000e+00)
( 2.12000000e-01,  9.85703184e-01)
( 2.24000000e-01,  9.44027483e-01)
( 2.36000000e-01,  8.78446739e-01)
( 2.48000000e-01,  7.94215853e-01)
( 2.60000000e-01,  6.97676326e-01)
( 2.72000000e-01,  5.95472542e-01)
( 2.84000000e-01,  4.93812198e-01)
( 2.96000000e-01,  3.97881920e-01)
( 3.08000000e-01,  3.11486276e-01)
( 3.20000000e-01,  2.36927759e-01)
( 3.32000000e-01,  1.75099657e-01)
( 3.44000000e-01,  1.25732330e-01)
( 3.56000000e-01,  8.77204698e-02)
( 3.68000000e-01,  5.94630600e-02)
( 3.80000000e-01,  3.91638951e-02)
( 3.92000000e-01,  2.50620633e-02)
( 4.04000000e-01,  1.55826058e-02)
( 4.16000000e-01,  9.41359837e-03)
( 4.28000000e-01,  5.52539799e-03)
( 4.40000000e-01,  3.15111160e-03)
( 4.52000000e-01,  1.74604858e-03)
( 4.64000000e-01,  9.40028842e-04)
( 4.76000000e-01,  4.91720545e-04)
( 4.88000000e-01,  2.49912425e-04)
( 5.00000000e-01,  1.23409804e-04)
( 5.12000000e-01,  5.92111908e-05)
( 5.24000000e-01,  2.76026162e-05)
( 5.36000000e-01,  1.25022741e-05)
( 5.48000000e-01,  5.50199392e-06)
( 5.60000000e-01,  2.35257520e-06)
( 5.72000000e-01,  9.77370534e-07)
( 5.84000000e-01,  3.94518422e-07)
( 5.96000000e-01,  1.54727551e-07)
( 6.08000000e-01,  5.89603885e-08)
( 6.20000000e-01,  2.18295780e-08)
( 6.32000000e-01,  7.85276608e-09)
( 6.44000000e-01,  2.74468369e-09)
( 6.56000000e-01,  9.32082300e-10)
( 6.68000000e-01,  3.07544907e-10)
( 6.80000000e-01,  9.85950558e-11)
( 6.92000000e-01,  3.07110067e-11)
( 7.04000000e-01,  9.29448409e-12)
( 7.16000000e-01,  2.73305796e-12)
( 7.28000000e-01,  7.80844788e-13)
( 7.40000000e-01,  2.16756888e-13)
( 7.52000000e-01,  5.84619663e-14)
( 7.64000000e-01,  1.53202640e-14)
( 7.76000000e-01,  3.90077946e-15)
( 7.88000000e-01,  9.65003448e-16)
( 8.00000000e-01,  2.31952283e-16)
( 8.12000000e-01,  5.41702397e-17)
( 8.24000000e-01,  1.22917914e-17)
( 8.36000000e-01,  2.70995404e-18)
( 8.48000000e-01,  5.80498347e-19)
( 8.60000000e-01,  1.20818202e-19)
( 8.72000000e-01,  2.44318340e-20)
( 8.84000000e-01,  4.80034098e-21)
( 8.96000000e-01,  9.16390177e-22)
( 9.08000000e-01,  1.69973447e-22)
( 9.20000000e-01,  3.06319086e-23)
( 9.32000000e-01,  5.36363669e-24)
( 9.44000000e-01,  9.12508601e-25)
( 9.56000000e-01,  1.50836657e-25)
( 9.68000000e-01,  2.42252993e-26)
( 9.80000000e-01,  3.78027784e-27)
( 9.92000000e-01,  5.73153053e-28)};

\addplot [black, , mark repeat=0, forget plot]
coordinates {
(-1.00000000e+00,  7.55581902e-86)
(-9.88000000e-01,  2.14415988e-84)
(-9.76000000e-01,  5.91187310e-83)
(-9.64000000e-01,  1.58374538e-81)
(-9.52000000e-01,  4.12228425e-80)
(-9.40000000e-01,  1.04251624e-78)
(-9.28000000e-01,  2.56165148e-77)
(-9.16000000e-01,  6.11574790e-76)
(-9.04000000e-01,  1.41863747e-74)
(-8.92000000e-01,  3.19731626e-73)
(-8.80000000e-01,  7.00151599e-72)
(-8.68000000e-01,  1.48967282e-70)
(-8.56000000e-01,  3.07951284e-69)
(-8.44000000e-01,  6.18536689e-68)
(-8.32000000e-01,  1.20709439e-66)
(-8.20000000e-01,  2.28880774e-65)
(-8.08000000e-01,  4.21667099e-64)
(-7.96000000e-01,  7.54783569e-63)
(-7.84000000e-01,  1.31270602e-61)
(-7.72000000e-01,  2.21822113e-60)
(-7.60000000e-01,  3.64195449e-59)
(-7.48000000e-01,  5.80973753e-58)
(-7.36000000e-01,  9.00473275e-57)
(-7.24000000e-01,  1.35605550e-55)
(-7.12000000e-01,  1.98415875e-54)
(-7.00000000e-01,  2.82077009e-53)
(-6.88000000e-01,  3.89629007e-52)
(-6.76000000e-01,  5.22910295e-51)
(-6.64000000e-01,  6.81860326e-50)
(-6.52000000e-01,  8.63885073e-49)
(-6.40000000e-01,  1.06342983e-47)
(-6.28000000e-01,  1.27190309e-46)
(-6.16000000e-01,  1.47805818e-45)
(-6.04000000e-01,  1.66886556e-44)
(-5.92000000e-01,  1.83081093e-43)
(-5.80000000e-01,  1.95145238e-42)
(-5.68000000e-01,  2.02099267e-41)
(-5.56000000e-01,  2.03359207e-40)
(-5.44000000e-01,  1.98817797e-39)
(-5.32000000e-01,  1.88859569e-38)
(-5.20000000e-01,  1.74307090e-37)
(-5.08000000e-01,  1.56308799e-36)
(-4.96000000e-01,  1.36189656e-35)
(-4.84000000e-01,  1.15291459e-34)
(-4.72000000e-01,  9.48292831e-34)
(-4.60000000e-01,  7.57844527e-33)
(-4.48000000e-01,  5.88450706e-32)
(-4.36000000e-01,  4.43948280e-31)
(-4.24000000e-01,  3.25422057e-30)
(-4.12000000e-01,  2.31768234e-29)
(-4.00000000e-01,  1.60381089e-28)
(-3.88000000e-01,  1.07831271e-27)
(-3.76000000e-01,  7.04415018e-27)
(-3.64000000e-01,  4.47100143e-26)
(-3.52000000e-01,  2.75723213e-25)
(-3.40000000e-01,  1.65209178e-24)
(-3.28000000e-01,  9.61805664e-24)
(-3.16000000e-01,  5.44042476e-23)
(-3.04000000e-01,  2.98999599e-22)
(-2.92000000e-01,  1.59661690e-21)
(-2.80000000e-01,  8.28367701e-21)
(-2.68000000e-01,  4.17578289e-20)
(-2.56000000e-01,  2.04524332e-19)
(-2.44000000e-01,  9.73294738e-19)
(-2.32000000e-01,  4.50024425e-18)
(-2.20000000e-01,  2.02171585e-17)
(-2.08000000e-01,  8.82462916e-17)
(-1.96000000e-01,  3.74252863e-16)
(-1.84000000e-01,  1.54214822e-15)
(-1.72000000e-01,  6.17418194e-15)
(-1.60000000e-01,  2.40173478e-14)
(-1.48000000e-01,  9.07743172e-14)
(-1.36000000e-01,  3.33344468e-13)
(-1.24000000e-01,  1.18936686e-12)
(-1.12000000e-01,  4.12316518e-12)
(-1.00000000e-01,  1.38879439e-11)
(-8.80000000e-02,  4.54503770e-11)
(-7.60000000e-02,  1.44520464e-10)
(-6.40000000e-02,  4.46491836e-10)
(-5.20000000e-02,  1.34026296e-09)
(-4.00000000e-02,  3.90893843e-09)
(-2.80000000e-02,  1.10769441e-08)
(-1.60000000e-02,  3.04981442e-08)
(-4.00000000e-03,  8.15866656e-08)
( 8.00000000e-03,  2.12059275e-07)
( 2.00000000e-02,  5.35534780e-07)
( 3.20000000e-02,  1.31404551e-06)
( 4.40000000e-02,  3.13274795e-06)
( 5.60000000e-02,  7.25659389e-06)
( 6.80000000e-02,  1.63317417e-05)
( 8.00000000e-02,  3.57128496e-05)
( 9.20000000e-02,  7.58767686e-05)
( 1.04000000e-01,  1.56633794e-04)
( 1.16000000e-01,  3.14162560e-04)
( 1.28000000e-01,  6.12231548e-04)
( 1.40000000e-01,  1.15922917e-03)
( 1.52000000e-01,  2.13262855e-03)
( 1.64000000e-01,  3.81200493e-03)
( 1.76000000e-01,  6.62039661e-03)
( 1.88000000e-01,  1.11713814e-02)
( 2.00000000e-01,  1.83156389e-02)
( 2.12000000e-01,  2.91762575e-02)
( 2.24000000e-01,  4.51574503e-02)
( 2.36000000e-01,  6.79080972e-02)
( 2.48000000e-01,  9.92215550e-02)
( 2.60000000e-01,  1.40858421e-01)
( 2.72000000e-01,  1.94290659e-01)
( 2.84000000e-01,  2.60383431e-01)
( 2.96000000e-01,  3.39052607e-01)
( 3.08000000e-01,  4.28956399e-01)
( 3.20000000e-01,  5.27292424e-01)
( 3.32000000e-01,  6.29770381e-01)
( 3.44000000e-01,  7.30811294e-01)
( 3.56000000e-01,  8.23987433e-01)
( 3.68000000e-01,  9.02668412e-01)
( 3.80000000e-01,  9.60789439e-01)
( 3.92000000e-01,  9.93620436e-01)
( 4.04000000e-01,  9.98401279e-01)
( 4.16000000e-01,  9.74724902e-01)
( 4.28000000e-01,  9.24594515e-01)
( 4.40000000e-01,  8.52143789e-01)
( 4.52000000e-01,  7.63074204e-01)
( 4.64000000e-01,  6.63915763e-01)
( 4.76000000e-01,  5.61243736e-01)
( 4.88000000e-01,  4.60980286e-01)
( 5.00000000e-01,  3.67879441e-01)
( 5.12000000e-01,  2.85246945e-01)
( 5.24000000e-01,  2.14896234e-01)
( 5.36000000e-01,  1.57300074e-01)
( 5.48000000e-01,  1.11871991e-01)
( 5.60000000e-01,  7.73047404e-02)
( 5.72000000e-01,  5.19018938e-02)
( 5.84000000e-01,  3.38573219e-02)
( 5.96000000e-01,  2.14592391e-02)
( 6.08000000e-01,  1.32150361e-02)
( 6.20000000e-01,  7.90705405e-03)
( 6.32000000e-01,  4.59677642e-03)
( 6.44000000e-01,  2.59647636e-03)
( 6.56000000e-01,  1.42497644e-03)
( 6.68000000e-01,  7.59842031e-04)
( 6.80000000e-01,  3.93669041e-04)
( 6.92000000e-01,  1.98167089e-04)
( 7.04000000e-01,  9.69223903e-05)
( 7.16000000e-01,  4.60584185e-05)
( 7.28000000e-01,  2.12660214e-05)
( 7.40000000e-01,  9.54016287e-06)
( 7.52000000e-01,  4.15831755e-06)
( 7.64000000e-01,  1.76105059e-06)
( 7.76000000e-01,  7.24633375e-07)
( 7.88000000e-01,  2.89705774e-07)
( 8.00000000e-01,  1.12535175e-07)
( 8.12000000e-01,  4.24728835e-08)
( 8.24000000e-01,  1.55749809e-08)
( 8.36000000e-01,  5.54926694e-09)
( 8.48000000e-01,  1.92103833e-09)
( 8.60000000e-01,  6.46143177e-10)
( 8.72000000e-01,  2.11161051e-10)
( 8.84000000e-01,  6.70488271e-11)
( 8.96000000e-01,  2.06852548e-11)
( 9.08000000e-01,  6.20044532e-12)
( 9.20000000e-01,  1.80583144e-12)
( 9.32000000e-01,  5.11003477e-13)
( 9.44000000e-01,  1.40495631e-13)
( 9.56000000e-01,  3.75313429e-14)
( 9.68000000e-01,  9.74131769e-15)
( 9.80000000e-01,  2.45659537e-15)
( 9.92000000e-01,  6.01924290e-16)};

\end{axis}
\end{tikzpicture}

%% file: figs/methodology/bump1d_motivation/bump_1d_snapshots_xref.tikz
\begin{tikzpicture}
\begin{axis}[
width=1.\textwidth,
xlabel={$\mathcal{H} (x; c, 0)$},
ymax=1.2,
xmax=1,
ylabel={$y(x; c)$},
xmin=-1,
ymin=-0.2]
\addplot [black, , mark repeat=0, forget plot]
coordinates {
(-1.00000000e+00,  2.31952283e-16)
(-9.88000000e-01,  4.83488556e-16)
(-9.76000000e-01,  1.00312685e-15)
(-9.64000000e-01,  2.06909054e-15)
(-9.52000000e-01,  4.24480282e-15)
(-9.40000000e-01,  8.65999136e-15)
(-9.28000000e-01,  1.75683152e-14)
(-9.16000000e-01,  3.54435216e-14)
(-9.04000000e-01,  7.10891079e-14)
(-8.92000000e-01,  1.41789060e-13)
(-8.80000000e-01,  2.81121741e-13)
(-8.68000000e-01,  5.54166268e-13)
(-8.56000000e-01,  1.08594563e-12)
(-8.44000000e-01,  2.11522193e-12)
(-8.32000000e-01,  4.09571323e-12)
(-8.20000000e-01,  7.88140009e-12)
(-8.08000000e-01,  1.50742913e-11)
(-7.96000000e-01,  2.86529242e-11)
(-7.84000000e-01,  5.41188125e-11)
(-7.72000000e-01,  1.01580496e-10)
(-7.60000000e-01,  1.89435417e-10)
(-7.48000000e-01,  3.50996848e-10)
(-7.36000000e-01,  6.46143177e-10)
(-7.24000000e-01,  1.18155653e-09)
(-7.12000000e-01,  2.14631071e-09)
(-7.00000000e-01,  3.87270434e-09)
(-6.88000000e-01,  6.93986029e-09)
(-6.76000000e-01,  1.23510796e-08)
(-6.64000000e-01,  2.18295780e-08)
(-6.52000000e-01,  3.83097645e-08)
(-6.40000000e-01,  6.67555497e-08)
(-6.28000000e-01,  1.15491907e-07)
(-6.16000000e-01,  1.98355964e-07)
(-6.04000000e-01,  3.38174542e-07)
(-5.92000000e-01,  5.72288589e-07)
(-5.80000000e-01,  9.61217989e-07)
(-5.68000000e-01,  1.60219464e-06)
(-5.56000000e-01,  2.65014269e-06)
(-5.44000000e-01,  4.34957914e-06)
(-5.32000000e-01,  7.08275208e-06)
(-5.20000000e-01,  1.14417141e-05)
(-5.08000000e-01,  1.83349721e-05)
(-4.96000000e-01,  2.91426194e-05)
(-4.84000000e-01,  4.59398886e-05)
(-4.72000000e-01,  7.18154426e-05)
(-4.60000000e-01,  1.11319470e-04)
(-4.48000000e-01,  1.71082679e-04)
(-4.36000000e-01,  2.60660773e-04)
(-4.24000000e-01,  3.93669041e-04)
(-4.12000000e-01,  5.89282598e-04)
(-4.00000000e-01,  8.74187670e-04)
(-3.88000000e-01,  1.28506403e-03)
(-3.76000000e-01,  1.87168346e-03)
(-3.64000000e-01,  2.70068978e-03)
(-3.52000000e-01,  3.86009173e-03)
(-3.40000000e-01,  5.46444534e-03)
(-3.28000000e-01,  7.66062295e-03)
(-3.16000000e-01,  1.06339454e-02)
(-3.04000000e-01,  1.46143135e-02)
(-2.92000000e-01,  1.98818121e-02)
(-2.80000000e-01,  2.67710587e-02)
(-2.68000000e-01,  3.56733759e-02)
(-2.56000000e-01,  4.70356911e-02)
(-2.44000000e-01,  6.13549416e-02)
(-2.32000000e-01,  7.91667326e-02)
(-2.20000000e-01,  1.01027095e-01)
(-2.08000000e-01,  1.27486467e-01)
(-1.96000000e-01,  1.59055513e-01)
(-1.84000000e-01,  1.96163093e-01)
(-1.72000000e-01,  2.39107630e-01)
(-1.60000000e-01,  2.88004180e-01)
(-1.48000000e-01,  3.42730700e-01)
(-1.36000000e-01,  4.02878058e-01)
(-1.24000000e-01,  4.67709255e-01)
(-1.12000000e-01,  5.36133741e-01)
(-1.00000000e-01,  6.06702634e-01)
(-8.80000000e-02,  6.77629783e-01)
(-7.60000000e-02,  7.46841983e-01)
(-6.40000000e-02,  8.12059277e-01)
(-5.20000000e-02,  8.70903262e-01)
(-4.00000000e-02,  9.21028004e-01)
(-2.80000000e-02,  9.60264949e-01)
(-1.60000000e-02,  9.86770463e-01)
(-4.00000000e-03,  9.99162970e-01)
( 8.00000000e-03,  9.96636112e-01)
( 2.00000000e-02,  9.79035623e-01)
( 3.20000000e-02,  9.46890569e-01)
( 4.40000000e-02,  9.01393722e-01)
( 5.60000000e-02,  8.44331099e-01)
( 6.80000000e-02,  7.77965859e-01)
( 8.00000000e-02,  7.04886848e-01)
( 9.20000000e-02,  6.27836386e-01)
( 1.04000000e-01,  5.49533480e-01)
( 1.16000000e-01,  4.72508715e-01)
( 1.28000000e-01,  3.98965811e-01)
( 1.40000000e-01,  3.30680533e-01)
( 1.52000000e-01,  2.68942635e-01)
( 1.64000000e-01,  2.14542734e-01)
( 1.76000000e-01,  1.67799277e-01)
( 1.88000000e-01,  1.28617617e-01)
( 2.00000000e-01,  9.65717914e-02)
( 2.12000000e-01,  7.09965881e-02)
( 2.24000000e-01,  5.10801525e-02)
( 2.36000000e-01,  3.59481465e-02)
( 2.48000000e-01,  2.47334112e-02)
( 2.60000000e-01,  1.66279570e-02)
( 2.72000000e-01,  1.09168277e-02)
( 2.84000000e-01,  6.99521106e-03)
( 2.96000000e-01,  4.37206447e-03)
( 3.08000000e-01,  2.66364467e-03)
( 3.20000000e-01,  1.58080933e-03)
( 3.32000000e-01,  9.13263968e-04)
( 3.44000000e-01,  5.13225706e-04)
( 3.56000000e-01,  2.80340451e-04)
( 3.68000000e-01,  1.48725502e-04)
( 3.80000000e-01,  7.65672628e-05)
( 3.92000000e-01,  3.82187366e-05)
( 4.04000000e-01,  1.84794737e-05)
( 4.16000000e-01,  8.64700425e-06)
( 4.28000000e-01,  3.91168313e-06)
( 4.40000000e-01,  1.70890867e-06)
( 4.52000000e-01,  7.20181568e-07)
( 4.64000000e-01,  2.92427344e-07)
( 4.76000000e-01,  1.14262418e-07)
( 4.88000000e-01,  4.29065303e-08)
( 5.00000000e-01,  1.54620969e-08)
( 5.12000000e-01,  5.33937709e-09)
( 5.24000000e-01,  1.76401388e-09)
( 5.36000000e-01,  5.56639627e-10)
( 5.48000000e-01,  1.67470383e-10)
( 5.60000000e-01,  4.79485797e-11)
( 5.72000000e-01,  1.30370121e-11)
( 5.84000000e-01,  3.35879083e-12)
( 5.96000000e-01,  8.18161233e-13)
( 6.08000000e-01,  1.87922789e-13)
( 6.20000000e-01,  4.05959621e-14)
( 6.32000000e-01,  8.22318165e-15)
( 6.44000000e-01,  1.55719080e-15)
( 6.56000000e-01,  2.74650060e-16)
( 6.68000000e-01,  4.49631330e-17)
( 6.80000000e-01,  6.80429921e-18)
( 6.92000000e-01,  9.47507300e-19)
( 7.04000000e-01,  1.20818202e-19)
( 7.16000000e-01,  1.40316641e-20)
( 7.28000000e-01,  1.47539635e-21)
( 7.40000000e-01,  1.39514394e-22)
( 7.52000000e-01,  1.17786621e-23)
( 7.64000000e-01,  8.80758176e-25)
( 7.76000000e-01,  5.77366965e-26)
( 7.88000000e-01,  3.28651108e-27)
( 8.00000000e-01,  1.60381089e-28)
( 8.12000000e-01,  6.61617447e-30)
( 8.24000000e-01,  2.27081292e-31)
( 8.36000000e-01,  6.35960190e-33)
( 8.48000000e-01,  1.42100987e-34)
( 8.60000000e-01,  2.46678199e-36)
( 8.72000000e-01,  3.21784027e-38)
( 8.84000000e-01,  3.02612828e-40)
( 8.96000000e-01,  1.95145238e-42)
( 9.08000000e-01,  8.05480581e-45)
( 9.20000000e-01,  1.94823938e-47)
( 9.32000000e-01,  2.43994733e-50)
( 9.44000000e-01,  1.31949509e-53)
( 9.56000000e-01,  2.32333389e-57)
( 9.68000000e-01,  8.17407472e-62)
( 9.80000000e-01,  2.11119225e-67)
( 9.92000000e-01,  2.22598655e-75)};

\addplot [black, , mark repeat=0, forget plot]
coordinates {
(-1.00000000e+00,  1.60381089e-28)
(-9.88000000e-01,  6.18515653e-28)
(-9.76000000e-01,  2.35859495e-27)
(-9.64000000e-01,  8.89583699e-27)
(-9.52000000e-01,  3.31811307e-26)
(-9.40000000e-01,  1.22383044e-25)
(-9.28000000e-01,  4.46286474e-25)
(-9.16000000e-01,  1.60875445e-24)
(-9.04000000e-01,  5.73152195e-24)
(-8.92000000e-01,  2.01819902e-23)
(-8.80000000e-01,  7.02412635e-23)
(-8.68000000e-01,  2.41615859e-22)
(-8.56000000e-01,  8.21287425e-22)
(-8.44000000e-01,  2.75807295e-21)
(-8.32000000e-01,  9.14907028e-21)
(-8.20000000e-01,  2.99809189e-20)
(-8.08000000e-01,  9.70534297e-20)
(-7.96000000e-01,  3.10321102e-19)
(-7.84000000e-01,  9.79805724e-19)
(-7.72000000e-01,  3.05434828e-18)
(-7.60000000e-01,  9.40133116e-18)
(-7.48000000e-01,  2.85713933e-17)
(-7.36000000e-01,  8.57140489e-17)
(-7.24000000e-01,  2.53768798e-16)
(-7.12000000e-01,  7.41494036e-16)
(-7.00000000e-01,  2.13823715e-15)
(-6.88000000e-01,  6.08399605e-15)
(-6.76000000e-01,  1.70764612e-14)
(-6.64000000e-01,  4.72830845e-14)
(-6.52000000e-01,  1.29146066e-13)
(-6.40000000e-01,  3.47862887e-13)
(-6.28000000e-01,  9.23916122e-13)
(-6.16000000e-01,  2.41972140e-12)
(-6.04000000e-01,  6.24763765e-12)
(-5.92000000e-01,  1.58996858e-11)
(-5.80000000e-01,  3.98835230e-11)
(-5.68000000e-01,  9.85950558e-11)
(-5.56000000e-01,  2.40144224e-10)
(-5.44000000e-01,  5.76292522e-10)
(-5.32000000e-01,  1.36234969e-09)
(-5.20000000e-01,  3.17190140e-09)
(-5.08000000e-01,  7.27314707e-09)
(-4.96000000e-01,  1.64210800e-08)
(-4.84000000e-01,  3.65002580e-08)
(-4.72000000e-01,  7.98665005e-08)
(-4.60000000e-01,  1.71991418e-07)
(-4.48000000e-01,  3.64489991e-07)
(-4.36000000e-01,  7.60010424e-07)
(-4.24000000e-01,  1.55896791e-06)
(-4.12000000e-01,  3.14544912e-06)
(-4.00000000e-01,  6.24114975e-06)
(-3.88000000e-01,  1.21767544e-05)
(-3.76000000e-01,  2.33556633e-05)
(-3.64000000e-01,  4.40337230e-05)
(-3.52000000e-01,  8.15878862e-05)
(-3.40000000e-01,  1.48539447e-04)
(-3.28000000e-01,  2.65675396e-04)
(-3.16000000e-01,  4.66741690e-04)
(-3.04000000e-01,  8.05257127e-04)
(-2.92000000e-01,  1.36409237e-03)
(-2.80000000e-01,  2.26839422e-03)
(-2.68000000e-01,  3.70232519e-03)
(-2.56000000e-01,  5.92956498e-03)
(-2.44000000e-01,  9.31700185e-03)
(-2.32000000e-01,  1.43596021e-02)
(-2.20000000e-01,  2.17035461e-02)
(-2.08000000e-01,  3.21622328e-02)
(-1.96000000e-01,  4.67190273e-02)
(-1.84000000e-01,  6.65082992e-02)
(-1.72000000e-01,  9.27669192e-02)
(-1.60000000e-01,  1.26748969e-01)
(-1.48000000e-01,  1.69600356e-01)
(-1.36000000e-01,  2.22195333e-01)
(-1.24000000e-01,  2.84945473e-01)
(-1.12000000e-01,  3.57600558e-01)
(-1.00000000e-01,  4.39068863e-01)
(-8.80000000e-02,  5.27292424e-01)
(-7.60000000e-02,  6.19211544e-01)
(-6.40000000e-02,  7.10848154e-01)
(-5.20000000e-02,  7.97524954e-01)
(-4.00000000e-02,  8.74214946e-01)
(-2.80000000e-02,  9.35993744e-01)
(-1.60000000e-02,  9.78543500e-01)
(-4.00000000e-03,  9.98639539e-01)
( 8.00000000e-03,  9.94546107e-01)
( 2.00000000e-02,  9.66255408e-01)
( 3.20000000e-02,  9.15524919e-01)
( 4.40000000e-02,  8.45701083e-01)
( 5.60000000e-02,  7.61353240e-01)
( 6.80000000e-02,  6.67773230e-01)
( 8.00000000e-02,  5.70417532e-01)
( 9.20000000e-02,  4.74375087e-01)
( 1.04000000e-01,  3.83933974e-01)
( 1.16000000e-01,  3.02296923e-01)
( 1.28000000e-01,  2.31465615e-01)
( 1.40000000e-01,  1.72283611e-01)
( 1.52000000e-01,  1.24603796e-01)
( 1.64000000e-01,  8.75324949e-02)
( 1.76000000e-01,  5.97000423e-02)
( 1.88000000e-01,  3.95147608e-02)
( 2.00000000e-01,  2.53706109e-02)
( 2.12000000e-01,  1.57939958e-02)
( 2.24000000e-01,  9.52884711e-03)
( 2.36000000e-01,  5.56888862e-03)
( 2.48000000e-01,  3.15111160e-03)
( 2.60000000e-01,  1.72547152e-03)
( 2.72000000e-01,  9.13853285e-04)
( 2.84000000e-01,  4.67885736e-04)
( 2.96000000e-01,  2.31452398e-04)
( 3.08000000e-01,  1.10560467e-04)
( 3.20000000e-01,  5.09689647e-05)
( 3.32000000e-01,  2.26633612e-05)
( 3.44000000e-01,  9.71386880e-06)
( 3.56000000e-01,  4.01084351e-06)
( 3.68000000e-01,  1.59431900e-06)
( 3.80000000e-01,  6.09718127e-07)
( 3.92000000e-01,  2.24185500e-07)
( 4.04000000e-01,  7.91954431e-08)
( 4.16000000e-01,  2.68586284e-08)
( 4.28000000e-01,  8.73887052e-09)
( 4.40000000e-01,  2.72575612e-09)
( 4.52000000e-01,  8.14331753e-10)
( 4.64000000e-01,  2.32853298e-10)
( 4.76000000e-01,  6.36740895e-11)
( 4.88000000e-01,  1.66350737e-11)
( 5.00000000e-01,  4.14912331e-12)
( 5.12000000e-01,  9.86872240e-13)
( 5.24000000e-01,  2.23678968e-13)
( 5.36000000e-01,  4.82522218e-14)
( 5.48000000e-01,  9.89935709e-15)
( 5.60000000e-01,  1.92885408e-15)
( 5.72000000e-01,  3.56682123e-16)
( 5.84000000e-01,  6.25009311e-17)
( 5.96000000e-01,  1.03690406e-17)
( 6.08000000e-01,  1.62666462e-18)
( 6.20000000e-01,  2.40919897e-19)
( 6.32000000e-01,  3.36558142e-20)
( 6.44000000e-01,  4.42858061e-21)
( 6.56000000e-01,  5.47959336e-22)
( 6.68000000e-01,  6.36609732e-23)
( 6.80000000e-01,  6.93535626e-24)
( 6.92000000e-01,  7.07371078e-25)
( 7.04000000e-01,  6.74313259e-26)
( 7.16000000e-01,  5.99704739e-27)
( 7.28000000e-01,  4.96681944e-28)
( 7.40000000e-01,  3.82348241e-29)
( 7.52000000e-01,  2.73033545e-30)
( 7.64000000e-01,  1.80486425e-31)
( 7.76000000e-01,  1.10207444e-32)
( 7.88000000e-01,  6.20265775e-34)
( 8.00000000e-01,  3.21076164e-35)
( 8.12000000e-01,  1.52467639e-36)
( 8.24000000e-01,  6.61821007e-38)
( 8.36000000e-01,  2.62066124e-39)
( 8.48000000e-01,  9.44975498e-41)
( 8.60000000e-01,  3.08574301e-42)
( 8.72000000e-01,  9.11691910e-44)
( 8.84000000e-01,  2.42488103e-45)
( 8.96000000e-01,  5.79438801e-47)
( 9.08000000e-01,  1.23612189e-48)
( 9.20000000e-01,  2.35271813e-50)
( 9.32000000e-01,  3.96849078e-52)
( 9.44000000e-01,  5.90151076e-54)
( 9.56000000e-01,  7.71293036e-56)
( 9.68000000e-01,  8.81213486e-58)
( 9.80000000e-01,  8.74910084e-60)
( 9.92000000e-01,  7.50275959e-62)};

\addplot [black, , mark repeat=0, forget plot]
coordinates {
(-1.00000000e+00,  3.72007598e-44)
(-9.88000000e-01,  4.04207833e-43)
(-9.76000000e-01,  4.26726841e-42)
(-9.64000000e-01,  4.37711054e-41)
(-9.52000000e-01,  4.36231866e-40)
(-9.40000000e-01,  4.22415241e-39)
(-9.28000000e-01,  3.97424000e-38)
(-9.16000000e-01,  3.63296254e-37)
(-9.04000000e-01,  3.22671099e-36)
(-8.92000000e-01,  2.78452765e-35)
(-8.80000000e-01,  2.33472278e-34)
(-8.68000000e-01,  1.90200400e-33)
(-8.56000000e-01,  1.50549685e-32)
(-8.44000000e-01,  1.15781879e-31)
(-8.32000000e-01,  8.65154470e-31)
(-8.20000000e-01,  6.28114815e-30)
(-8.08000000e-01,  4.43074480e-29)
(-7.96000000e-01,  3.03673402e-28)
(-7.84000000e-01,  2.02222322e-27)
(-7.72000000e-01,  1.30840968e-26)
(-7.60000000e-01,  8.22528065e-26)
(-7.48000000e-01,  5.02400450e-25)
(-7.36000000e-01,  2.98154682e-24)
(-7.24000000e-01,  1.71919669e-23)
(-7.12000000e-01,  9.63167493e-23)
(-7.00000000e-01,  5.24288566e-22)
(-6.88000000e-01,  2.77288128e-21)
(-6.76000000e-01,  1.42490032e-20)
(-6.64000000e-01,  7.11426455e-20)
(-6.52000000e-01,  3.45118195e-19)
(-6.40000000e-01,  1.62666462e-18)
(-6.28000000e-01,  7.44938645e-18)
(-6.16000000e-01,  3.31463191e-17)
(-6.04000000e-01,  1.43298759e-16)
(-5.92000000e-01,  6.01924290e-16)
(-5.80000000e-01,  2.45659537e-15)
(-5.68000000e-01,  9.74131769e-15)
(-5.56000000e-01,  3.75313429e-14)
(-5.44000000e-01,  1.40495631e-13)
(-5.32000000e-01,  5.11003477e-13)
(-5.20000000e-01,  1.80583144e-12)
(-5.08000000e-01,  6.20044532e-12)
(-4.96000000e-01,  2.06852548e-11)
(-4.84000000e-01,  6.70488271e-11)
(-4.72000000e-01,  2.11161051e-10)
(-4.60000000e-01,  6.46143177e-10)
(-4.48000000e-01,  1.92103833e-09)
(-4.36000000e-01,  5.54926694e-09)
(-4.24000000e-01,  1.55749809e-08)
(-4.12000000e-01,  4.24728835e-08)
(-4.00000000e-01,  1.12535175e-07)
(-3.88000000e-01,  2.89705774e-07)
(-3.76000000e-01,  7.24633375e-07)
(-3.64000000e-01,  1.76105059e-06)
(-3.52000000e-01,  4.15831755e-06)
(-3.40000000e-01,  9.54016287e-06)
(-3.28000000e-01,  2.12660214e-05)
(-3.16000000e-01,  4.60584185e-05)
(-3.04000000e-01,  9.69223903e-05)
(-2.92000000e-01,  1.98167089e-04)
(-2.80000000e-01,  3.93669041e-04)
(-2.68000000e-01,  7.59842031e-04)
(-2.56000000e-01,  1.42497644e-03)
(-2.44000000e-01,  2.59647636e-03)
(-2.32000000e-01,  4.59677642e-03)
(-2.20000000e-01,  7.90705405e-03)
(-2.08000000e-01,  1.32150361e-02)
(-1.96000000e-01,  2.14592391e-02)
(-1.84000000e-01,  3.38573219e-02)
(-1.72000000e-01,  5.19018938e-02)
(-1.60000000e-01,  7.73047404e-02)
(-1.48000000e-01,  1.11871991e-01)
(-1.36000000e-01,  1.57300074e-01)
(-1.24000000e-01,  2.14896234e-01)
(-1.12000000e-01,  2.85246945e-01)
(-1.00000000e-01,  3.67879441e-01)
(-8.80000000e-02,  4.60980286e-01)
(-7.60000000e-02,  5.61243736e-01)
(-6.40000000e-02,  6.63915763e-01)
(-5.20000000e-02,  7.63074204e-01)
(-4.00000000e-02,  8.52143789e-01)
(-2.80000000e-02,  9.24594515e-01)
(-1.60000000e-02,  9.74724902e-01)
(-4.00000000e-03,  9.98401279e-01)
( 8.00000000e-03,  9.93620436e-01)
( 2.00000000e-02,  9.60789439e-01)
( 3.20000000e-02,  9.02668412e-01)
( 4.40000000e-02,  8.23987433e-01)
( 5.60000000e-02,  7.30811294e-01)
( 6.80000000e-02,  6.29770381e-01)
( 8.00000000e-02,  5.27292424e-01)
( 9.20000000e-02,  4.28956399e-01)
( 1.04000000e-01,  3.39052607e-01)
( 1.16000000e-01,  2.60383431e-01)
( 1.28000000e-01,  1.94290659e-01)
( 1.40000000e-01,  1.40858421e-01)
( 1.52000000e-01,  9.92215550e-02)
( 1.64000000e-01,  6.79080972e-02)
( 1.76000000e-01,  4.51574503e-02)
( 1.88000000e-01,  2.91762575e-02)
( 2.00000000e-01,  1.83156389e-02)
( 2.12000000e-01,  1.11713814e-02)
( 2.24000000e-01,  6.62039661e-03)
( 2.36000000e-01,  3.81200493e-03)
( 2.48000000e-01,  2.13262855e-03)
( 2.60000000e-01,  1.15922917e-03)
( 2.72000000e-01,  6.12231548e-04)
( 2.84000000e-01,  3.14162560e-04)
( 2.96000000e-01,  1.56633794e-04)
( 3.08000000e-01,  7.58767686e-05)
( 3.20000000e-01,  3.57128496e-05)
( 3.32000000e-01,  1.63317417e-05)
( 3.44000000e-01,  7.25659389e-06)
( 3.56000000e-01,  3.13274795e-06)
( 3.68000000e-01,  1.31404551e-06)
( 3.80000000e-01,  5.35534780e-07)
( 3.92000000e-01,  2.12059275e-07)
( 4.04000000e-01,  8.15866656e-08)
( 4.16000000e-01,  3.04981442e-08)
( 4.28000000e-01,  1.10769441e-08)
( 4.40000000e-01,  3.90893843e-09)
( 4.52000000e-01,  1.34026296e-09)
( 4.64000000e-01,  4.46491836e-10)
( 4.76000000e-01,  1.44520464e-10)
( 4.88000000e-01,  4.54503770e-11)
( 5.00000000e-01,  1.38879439e-11)
( 5.12000000e-01,  4.12316518e-12)
( 5.24000000e-01,  1.18936686e-12)
( 5.36000000e-01,  3.33344468e-13)
( 5.48000000e-01,  9.07743172e-14)
( 5.60000000e-01,  2.40173478e-14)
( 5.72000000e-01,  6.17418194e-15)
( 5.84000000e-01,  1.54214822e-15)
( 5.96000000e-01,  3.74252863e-16)
( 6.08000000e-01,  8.82462916e-17)
( 6.20000000e-01,  2.02171585e-17)
( 6.32000000e-01,  4.50024425e-18)
( 6.44000000e-01,  9.73294738e-19)
( 6.56000000e-01,  2.04524332e-19)
( 6.68000000e-01,  4.17578289e-20)
( 6.80000000e-01,  8.28367701e-21)
( 6.92000000e-01,  1.59661690e-21)
( 7.04000000e-01,  2.98999599e-22)
( 7.16000000e-01,  5.44042476e-23)
( 7.28000000e-01,  9.61805664e-24)
( 7.40000000e-01,  1.65209178e-24)
( 7.52000000e-01,  2.75723213e-25)
( 7.64000000e-01,  4.47100143e-26)
( 7.76000000e-01,  7.04415018e-27)
( 7.88000000e-01,  1.07831271e-27)
( 8.00000000e-01,  1.60381089e-28)
( 8.12000000e-01,  2.31768234e-29)
( 8.24000000e-01,  3.25422057e-30)
( 8.36000000e-01,  4.43948280e-31)
( 8.48000000e-01,  5.88450706e-32)
( 8.60000000e-01,  7.57844527e-33)
( 8.72000000e-01,  9.48292831e-34)
( 8.84000000e-01,  1.15291459e-34)
( 8.96000000e-01,  1.36189656e-35)
( 9.08000000e-01,  1.56308799e-36)
( 9.20000000e-01,  1.74307090e-37)
( 9.32000000e-01,  1.88859569e-38)
( 9.44000000e-01,  1.98817797e-39)
( 9.56000000e-01,  2.03359207e-40)
( 9.68000000e-01,  2.02099267e-41)
( 9.80000000e-01,  1.95145238e-42)
( 9.92000000e-01,  1.83081093e-43)};

\addplot [black, , mark repeat=0, forget plot]
coordinates {
(-1.00000000e+00,  2.89464031e-63)
(-9.88000000e-01,  3.73145906e-61)
(-9.76000000e-01,  4.13967169e-59)
(-9.64000000e-01,  3.97704641e-57)
(-9.52000000e-01,  3.32736305e-55)
(-9.40000000e-01,  2.43562230e-53)
(-9.28000000e-01,  1.56751824e-51)
(-9.16000000e-01,  8.92372724e-50)
(-9.04000000e-01,  4.51759095e-48)
(-8.92000000e-01,  2.03546543e-46)
(-8.80000000e-01,  8.22194312e-45)
(-8.68000000e-01,  2.98246406e-43)
(-8.56000000e-01,  9.76452709e-42)
(-8.44000000e-01,  2.88963693e-40)
(-8.32000000e-01,  7.77280163e-39)
(-8.20000000e-01,  1.90219674e-37)
(-8.08000000e-01,  4.24951970e-36)
(-7.96000000e-01,  8.69715059e-35)
(-7.84000000e-01,  1.63413988e-33)
(-7.72000000e-01,  2.82500218e-32)
(-7.60000000e-01,  4.50382519e-31)
(-7.48000000e-01,  6.63747459e-30)
(-7.36000000e-01,  9.06243705e-29)
(-7.24000000e-01,  1.14866776e-27)
(-7.12000000e-01,  1.35422431e-26)
(-7.00000000e-01,  1.48773850e-25)
(-6.88000000e-01,  1.52563604e-24)
(-6.76000000e-01,  1.46287374e-23)
(-6.64000000e-01,  1.31383099e-22)
(-6.52000000e-01,  1.10708668e-21)
(-6.40000000e-01,  8.76600942e-21)
(-6.28000000e-01,  6.53003513e-20)
(-6.16000000e-01,  4.58272030e-19)
(-6.04000000e-01,  3.03497074e-18)
(-5.92000000e-01,  1.89890128e-17)
(-5.80000000e-01,  1.12360783e-16)
(-5.68000000e-01,  6.29782381e-16)
(-5.56000000e-01,  3.34641401e-15)
(-5.44000000e-01,  1.68794016e-14)
(-5.32000000e-01,  8.09035178e-14)
(-5.20000000e-01,  3.68853368e-13)
(-5.08000000e-01,  1.60130332e-12)
(-4.96000000e-01,  6.62587530e-12)
(-4.84000000e-01,  2.61551183e-11)
(-4.72000000e-01,  9.85950558e-11)
(-4.60000000e-01,  3.55186440e-10)
(-4.48000000e-01,  1.22404279e-09)
(-4.36000000e-01,  4.03835117e-09)
(-4.24000000e-01,  1.27649348e-08)
(-4.12000000e-01,  3.86909134e-08)
(-4.00000000e-01,  1.12535175e-07)
(-3.88000000e-01,  3.14314077e-07)
(-3.76000000e-01,  8.43646527e-07)
(-3.64000000e-01,  2.17760833e-06)
(-3.52000000e-01,  5.40888681e-06)
(-3.40000000e-01,  1.29367505e-05)
(-3.28000000e-01,  2.98133371e-05)
(-3.16000000e-01,  6.62421058e-05)
(-3.04000000e-01,  1.41989455e-04)
(-2.92000000e-01,  2.93784838e-04)
(-2.80000000e-01,  5.87084217e-04)
(-2.68000000e-01,  1.13372928e-03)
(-2.56000000e-01,  2.11685109e-03)
(-2.44000000e-01,  3.82358961e-03)
(-2.32000000e-01,  6.68458427e-03)
(-2.20000000e-01,  1.13166044e-02)
(-2.08000000e-01,  1.85612586e-02)
(-1.96000000e-01,  2.95089627e-02)
(-1.84000000e-01,  4.54942288e-02)
(-1.72000000e-01,  6.80471943e-02)
(-1.60000000e-01,  9.87885632e-02)
(-1.48000000e-01,  1.39261813e-01)
(-1.36000000e-01,  1.90707743e-01)
(-1.24000000e-01,  2.53800918e-01)
(-1.12000000e-01,  3.28382624e-01)
(-1.00000000e-01,  4.13236616e-01)
(-8.80000000e-02,  5.05958312e-01)
(-7.60000000e-02,  6.02961943e-01)
(-6.40000000e-02,  6.99652800e-01)
(-5.20000000e-02,  7.90765111e-01)
(-4.00000000e-02,  8.70835205e-01)
(-2.80000000e-02,  9.34751211e-01)
(-1.60000000e-02,  9.78301018e-01)
(-4.00000000e-03,  9.98635672e-01)
( 8.00000000e-03,  9.94576915e-01)
( 2.00000000e-02,  9.66723248e-01)
( 3.20000000e-02,  9.17342204e-01)
( 4.40000000e-02,  8.50072801e-01)
( 5.60000000e-02,  7.69491508e-01)
( 6.80000000e-02,  6.80611661e-01)
( 8.00000000e-02,  5.88389701e-01)
( 9.20000000e-02,  4.97302596e-01)
( 1.04000000e-01,  4.11039824e-01)
( 1.16000000e-01,  3.32330040e-01)
( 1.28000000e-01,  2.62900106e-01)
( 1.40000000e-01,  2.03543924e-01)
( 1.52000000e-01,  1.54269506e-01)
( 1.64000000e-01,  1.14488856e-01)
( 1.76000000e-01,  8.32169664e-02)
( 1.88000000e-01,  5.92554824e-02)
( 2.00000000e-01,  4.13441124e-02)
( 2.12000000e-01,  2.82726202e-02)
( 2.24000000e-01,  1.89531947e-02)
( 2.36000000e-01,  1.24582696e-02)
( 2.48000000e-01,  8.03129989e-03)
( 2.60000000e-01,  5.07875424e-03)
( 2.72000000e-01,  3.15111160e-03)
( 2.84000000e-01,  1.91862958e-03)
( 2.96000000e-01,  1.14664709e-03)
( 3.08000000e-01,  6.72760459e-04)
( 3.20000000e-01,  3.87589908e-04)
( 3.32000000e-01,  2.19301504e-04)
( 3.44000000e-01,  1.21886808e-04)
( 3.56000000e-01,  6.65560935e-05)
( 3.68000000e-01,  3.57128496e-05)
( 3.80000000e-01,  1.88334730e-05)
( 3.92000000e-01,  9.76327715e-06)
( 4.04000000e-01,  4.97599673e-06)
( 4.16000000e-01,  2.49384153e-06)
( 4.28000000e-01,  1.22921483e-06)
( 4.40000000e-01,  5.95970225e-07)
( 4.52000000e-01,  2.84278316e-07)
( 4.64000000e-01,  1.33421917e-07)
( 4.76000000e-01,  6.16273414e-08)
( 4.88000000e-01,  2.80171986e-08)
( 5.00000000e-01,  1.25386669e-08)
( 5.12000000e-01,  5.52509295e-09)
( 5.24000000e-01,  2.39722881e-09)
( 5.36000000e-01,  1.02438216e-09)
( 5.48000000e-01,  4.31175838e-10)
( 5.60000000e-01,  1.78775786e-10)
( 5.72000000e-01,  7.30359308e-11)
( 5.84000000e-01,  2.94020704e-11)
( 5.96000000e-01,  1.16641368e-11)
( 6.08000000e-01,  4.56117662e-12)
( 6.20000000e-01,  1.75830194e-12)
( 6.32000000e-01,  6.68191449e-13)
( 6.44000000e-01,  2.50388984e-13)
( 6.56000000e-01,  9.25349577e-14)
( 6.68000000e-01,  3.37247891e-14)
( 6.80000000e-01,  1.21234377e-14)
( 6.92000000e-01,  4.29974924e-15)
( 7.04000000e-01,  1.50459311e-15)
( 7.16000000e-01,  5.19439901e-16)
( 7.28000000e-01,  1.76972348e-16)
( 7.40000000e-01,  5.95141626e-17)
( 7.52000000e-01,  1.97553711e-17)
( 7.64000000e-01,  6.47255667e-18)
( 7.76000000e-01,  2.09363177e-18)
( 7.88000000e-01,  6.68736608e-19)
( 8.00000000e-01,  2.10947468e-19)
( 8.12000000e-01,  6.57077001e-20)
( 8.24000000e-01,  2.02127810e-20)
( 8.36000000e-01,  6.14204544e-21)
( 8.48000000e-01,  1.84397100e-21)
( 8.60000000e-01,  5.46985204e-22)
( 8.72000000e-01,  1.60302868e-22)
( 8.84000000e-01,  4.64165191e-23)
( 8.96000000e-01,  1.32820952e-23)
( 9.08000000e-01,  3.75668054e-24)
( 9.20000000e-01,  1.05036322e-24)
( 9.32000000e-01,  2.90329262e-25)
( 9.44000000e-01,  7.93295196e-26)
( 9.56000000e-01,  2.14277530e-26)
( 9.68000000e-01,  5.72258743e-27)
( 9.80000000e-01,  1.51128808e-27)
( 9.92000000e-01,  3.94358527e-28)};

\addplot [black, , mark repeat=0, forget plot]
coordinates {
(-1.00000000e+00,  7.55581902e-86)
(-9.88000000e-01,  2.80385592e-72)
(-9.76000000e-01,  2.28880774e-65)
(-9.64000000e-01,  3.12058405e-60)
(-9.52000000e-01,  4.81467295e-56)
(-9.40000000e-01,  1.80396818e-52)
(-9.28000000e-01,  2.45711193e-49)
(-9.16000000e-01,  1.54995778e-46)
(-9.04000000e-01,  5.29842349e-44)
(-8.92000000e-01,  1.09721255e-41)
(-8.80000000e-01,  1.49108640e-39)
(-8.68000000e-01,  1.41443116e-37)
(-8.56000000e-01,  9.81612352e-36)
(-8.44000000e-01,  5.18228807e-34)
(-8.32000000e-01,  2.14407871e-32)
(-8.20000000e-01,  7.13859819e-31)
(-8.08000000e-01,  1.95249187e-29)
(-7.96000000e-01,  4.46706676e-28)
(-7.84000000e-01,  8.68419298e-27)
(-7.72000000e-01,  1.45414731e-25)
(-7.60000000e-01,  2.11981968e-24)
(-7.48000000e-01,  2.72074446e-23)
(-7.36000000e-01,  3.10020660e-22)
(-7.24000000e-01,  3.16146848e-21)
(-7.12000000e-01,  2.90652986e-20)
(-7.00000000e-01,  2.42441140e-19)
(-6.88000000e-01,  1.84535747e-18)
(-6.76000000e-01,  1.28844229e-17)
(-6.64000000e-01,  8.29042770e-17)
(-6.52000000e-01,  4.93713311e-16)
(-6.40000000e-01,  2.73256162e-15)
(-6.28000000e-01,  1.41055544e-14)
(-6.16000000e-01,  6.81301789e-14)
(-6.04000000e-01,  3.08918836e-13)
(-5.92000000e-01,  1.31843003e-12)
(-5.80000000e-01,  5.31093959e-12)
(-5.68000000e-01,  2.02425187e-11)
(-5.56000000e-01,  7.31602277e-11)
(-5.44000000e-01,  2.51285178e-10)
(-5.32000000e-01,  8.21887109e-10)
(-5.20000000e-01,  2.56454429e-09)
(-5.08000000e-01,  7.64744296e-09)
(-4.96000000e-01,  2.18295780e-08)
(-4.84000000e-01,  5.97409316e-08)
(-4.72000000e-01,  1.56975151e-07)
(-4.60000000e-01,  3.96569865e-07)
(-4.48000000e-01,  9.64500424e-07)
(-4.36000000e-01,  2.26106306e-06)
(-4.24000000e-01,  5.11513097e-06)
(-4.12000000e-01,  1.11793260e-05)
(-4.00000000e-01,  2.36290355e-05)
(-3.88000000e-01,  4.83482631e-05)
(-3.76000000e-01,  9.58583000e-05)
(-3.64000000e-01,  1.84324425e-04)
(-3.52000000e-01,  3.44047112e-04)
(-3.40000000e-01,  6.23866773e-04)
(-3.28000000e-01,  1.09987107e-03)
(-3.16000000e-01,  1.88666631e-03)
(-3.04000000e-01,  3.15111160e-03)
(-2.92000000e-01,  5.12794641e-03)
(-2.80000000e-01,  8.13619858e-03)
(-2.68000000e-01,  1.25941903e-02)
(-2.56000000e-01,  1.90306010e-02)
(-2.44000000e-01,  2.80879426e-02)
(-2.32000000e-01,  4.05148728e-02)
(-2.20000000e-01,  5.71437329e-02)
(-2.08000000e-01,  7.88507687e-02)
(-1.96000000e-01,  1.06498047e-01)
(-1.84000000e-01,  1.40858421e-01)
(-1.72000000e-01,  1.82527896e-01)
(-1.60000000e-01,  2.31831869e-01)
(-1.48000000e-01,  2.88735483e-01)
(-1.36000000e-01,  3.52768195e-01)
(-1.24000000e-01,  4.22973997e-01)
(-1.12000000e-01,  4.97897546e-01)
(-1.00000000e-01,  5.75611870e-01)
(-8.80000000e-02,  6.53790940e-01)
(-7.60000000e-02,  7.29824419e-01)
(-6.40000000e-02,  8.00966297e-01)
(-5.20000000e-02,  8.64505608e-01)
(-4.00000000e-02,  9.17943693e-01)
(-2.80000000e-02,  9.59161142e-01)
(-1.60000000e-02,  9.86558775e-01)
(-4.00000000e-03,  9.99159621e-01)
( 8.00000000e-03,  9.96662843e-01)
( 2.00000000e-02,  9.79446001e-01)
( 3.20000000e-02,  9.48517763e-01)
( 4.40000000e-02,  9.05427713e-01)
( 5.60000000e-02,  8.52143789e-01)
( 6.80000000e-02,  7.90910202e-01)
( 8.00000000e-02,  7.24099353e-01)
( 9.20000000e-02,  6.54070413e-01)
( 1.04000000e-01,  5.83045149e-01)
( 1.16000000e-01,  5.13008523e-01)
( 1.28000000e-01,  4.45638352e-01)
( 1.40000000e-01,  3.82265050e-01)
( 1.52000000e-01,  3.23859639e-01)
( 1.64000000e-01,  2.71046110e-01)
( 1.76000000e-01,  2.24132792e-01)
( 1.88000000e-01,  1.83156855e-01)
( 2.00000000e-01,  1.47936105e-01)
( 2.12000000e-01,  1.18122913e-01)
( 2.24000000e-01,  9.32560356e-02)
( 2.36000000e-01,  7.28072539e-02)
( 2.48000000e-01,  5.62208720e-02)
( 2.60000000e-01,  4.29451748e-02)
( 2.72000000e-01,  3.24557809e-02)
( 2.84000000e-01,  2.42714706e-02)
( 2.96000000e-01,  1.79634758e-02)
( 3.08000000e-01,  1.31594374e-02)
( 3.20000000e-01,  9.54328413e-03)
( 3.32000000e-01,  6.85222368e-03)
( 3.44000000e-01,  4.87188108e-03)
( 3.56000000e-01,  3.43044457e-03)
( 3.68000000e-01,  2.39248170e-03)
( 3.80000000e-01,  1.65289776e-03)
( 3.92000000e-01,  1.13134683e-03)
( 4.04000000e-01,  7.67272278e-04)
( 4.16000000e-01,  5.15653200e-04)
( 4.28000000e-01,  3.43456040e-04)
( 4.40000000e-01,  2.26746056e-04)
( 4.52000000e-01,  1.48392454e-04)
( 4.64000000e-01,  9.62793775e-05)
( 4.76000000e-01,  6.19366218e-05)
( 4.88000000e-01,  3.95094731e-05)
( 5.00000000e-01,  2.49943309e-05)
( 5.12000000e-01,  1.56823770e-05)
( 5.24000000e-01,  9.76000867e-06)
( 5.36000000e-01,  6.02554201e-06)
( 5.48000000e-01,  3.69060339e-06)
( 5.60000000e-01,  2.24282693e-06)
( 5.72000000e-01,  1.35245333e-06)
( 5.84000000e-01,  8.09306295e-07)
( 5.96000000e-01,  4.80637551e-07)
( 6.08000000e-01,  2.83317063e-07)
( 6.20000000e-01,  1.65767003e-07)
( 6.32000000e-01,  9.62809558e-08)
( 6.44000000e-01,  5.55196230e-08)
( 6.56000000e-01,  3.17857061e-08)
( 6.68000000e-01,  1.80687348e-08)
( 6.80000000e-01,  1.01998337e-08)
( 6.92000000e-01,  5.71796089e-09)
( 7.04000000e-01,  3.18339377e-09)
( 7.16000000e-01,  1.76039343e-09)
( 7.28000000e-01,  9.66951983e-10)
( 7.40000000e-01,  5.27583946e-10)
( 7.52000000e-01,  2.85987199e-10)
( 7.64000000e-01,  1.54012271e-10)
( 7.76000000e-01,  8.24048137e-11)
( 7.88000000e-01,  4.38136318e-11)
( 8.00000000e-01,  2.31463568e-11)
( 8.12000000e-01,  1.21522536e-11)
( 8.24000000e-01,  6.34094852e-12)
( 8.36000000e-01,  3.28809880e-12)
( 8.48000000e-01,  1.69491440e-12)
( 8.60000000e-01,  8.68353392e-13)
( 8.72000000e-01,  4.42255616e-13)
( 8.84000000e-01,  2.23927247e-13)
( 8.96000000e-01,  1.12707271e-13)
( 9.08000000e-01,  5.64092494e-14)
( 9.20000000e-01,  2.80663127e-14)
( 9.32000000e-01,  1.38871509e-14)
( 9.44000000e-01,  6.83218735e-15)
( 9.56000000e-01,  3.34267017e-15)
( 9.68000000e-01,  1.62648740e-15)
( 9.80000000e-01,  7.87039243e-16)
( 9.92000000e-01,  3.78792550e-16)};

\end{axis}
\end{tikzpicture}

%% file: figs/methodology/bump1d_motivation/bump_1d_rom.tikz
\begin{tikzpicture}
\begin{axis}[
width=1.\textwidth,
xlabel={$x$},
ymax=1.2,
xmax=1,
ylabel={$y(x; c)$},
xmin=-1,
ymin=-0.2]
\addplot [blue, , mark repeat=0]
coordinates {
(-1.00000000e+00,  4.02006022e-74)
(-9.88000000e-01,  8.97382304e-73)
(-9.76000000e-01,  1.94632232e-71)
(-9.64000000e-01,  4.10151538e-70)
(-9.52000000e-01,  8.39781385e-69)
(-9.40000000e-01,  1.67063072e-67)
(-9.28000000e-01,  3.22914094e-66)
(-9.16000000e-01,  6.06437232e-65)
(-9.04000000e-01,  1.10656526e-63)
(-8.92000000e-01,  1.96182618e-62)
(-8.80000000e-01,  3.37937463e-61)
(-8.68000000e-01,  5.65593581e-60)
(-8.56000000e-01,  9.19739679e-59)
(-8.44000000e-01,  1.45317425e-57)
(-8.32000000e-01,  2.23081102e-56)
(-8.20000000e-01,  3.32736305e-55)
(-8.08000000e-01,  4.82203010e-54)
(-7.96000000e-01,  6.78972153e-53)
(-7.84000000e-01,  9.28894363e-52)
(-7.72000000e-01,  1.23473287e-50)
(-7.60000000e-01,  1.59467437e-49)
(-7.48000000e-01,  2.00107486e-48)
(-7.36000000e-01,  2.43975926e-47)
(-7.24000000e-01,  2.89016698e-46)
(-7.12000000e-01,  3.32652822e-45)
(-7.00000000e-01,  3.72007598e-44)
(-6.88000000e-01,  4.04207833e-43)
(-6.76000000e-01,  4.26726841e-42)
(-6.64000000e-01,  4.37711054e-41)
(-6.52000000e-01,  4.36231866e-40)
(-6.40000000e-01,  4.22415241e-39)
(-6.28000000e-01,  3.97424000e-38)
(-6.16000000e-01,  3.63296254e-37)
(-6.04000000e-01,  3.22671099e-36)
(-5.92000000e-01,  2.78452765e-35)
(-5.80000000e-01,  2.33472278e-34)
(-5.68000000e-01,  1.90200400e-33)
(-5.56000000e-01,  1.50549685e-32)
(-5.44000000e-01,  1.15781879e-31)
(-5.32000000e-01,  8.65154470e-31)
(-5.20000000e-01,  6.28114815e-30)
(-5.08000000e-01,  4.43074480e-29)
(-4.96000000e-01,  3.03673402e-28)
(-4.84000000e-01,  2.02222322e-27)
(-4.72000000e-01,  1.30840968e-26)
(-4.60000000e-01,  8.22528065e-26)
(-4.48000000e-01,  5.02400450e-25)
(-4.36000000e-01,  2.98154682e-24)
(-4.24000000e-01,  1.71919669e-23)
(-4.12000000e-01,  9.63167493e-23)
(-4.00000000e-01,  5.24288566e-22)
(-3.88000000e-01,  2.77288128e-21)
(-3.76000000e-01,  1.42490032e-20)
(-3.64000000e-01,  7.11426455e-20)
(-3.52000000e-01,  3.45118195e-19)
(-3.40000000e-01,  1.62666462e-18)
(-3.28000000e-01,  7.44938645e-18)
(-3.16000000e-01,  3.31463191e-17)
(-3.04000000e-01,  1.43298759e-16)
(-2.92000000e-01,  6.01924290e-16)
(-2.80000000e-01,  2.45659537e-15)
(-2.68000000e-01,  9.74131769e-15)
(-2.56000000e-01,  3.75313429e-14)
(-2.44000000e-01,  1.40495631e-13)
(-2.32000000e-01,  5.11003477e-13)
(-2.20000000e-01,  1.80583144e-12)
(-2.08000000e-01,  6.20044532e-12)
(-1.96000000e-01,  2.06852548e-11)
(-1.84000000e-01,  6.70488271e-11)
(-1.72000000e-01,  2.11161051e-10)
(-1.60000000e-01,  6.46143177e-10)
(-1.48000000e-01,  1.92103833e-09)
(-1.36000000e-01,  5.54926694e-09)
(-1.24000000e-01,  1.55749809e-08)
(-1.12000000e-01,  4.24728835e-08)
(-1.00000000e-01,  1.12535175e-07)
(-8.80000000e-02,  2.89705774e-07)
(-7.60000000e-02,  7.24633375e-07)
(-6.40000000e-02,  1.76105059e-06)
(-5.20000000e-02,  4.15831755e-06)
(-4.00000000e-02,  9.54016287e-06)
(-2.80000000e-02,  2.12660214e-05)
(-1.60000000e-02,  4.60584185e-05)
(-4.00000000e-03,  9.69223903e-05)
( 8.00000000e-03,  1.98167089e-04)
( 2.00000000e-02,  3.93669041e-04)
( 3.20000000e-02,  7.59842031e-04)
( 4.40000000e-02,  1.42497644e-03)
( 5.60000000e-02,  2.59647636e-03)
( 6.80000000e-02,  4.59677642e-03)
( 8.00000000e-02,  7.90705405e-03)
( 9.20000000e-02,  1.32150361e-02)
( 1.04000000e-01,  2.14592391e-02)
( 1.16000000e-01,  3.38573219e-02)
( 1.28000000e-01,  5.19018938e-02)
( 1.40000000e-01,  7.73047404e-02)
( 1.52000000e-01,  1.11871991e-01)
( 1.64000000e-01,  1.57300074e-01)
( 1.76000000e-01,  2.14896234e-01)
( 1.88000000e-01,  2.85246945e-01)
( 2.00000000e-01,  3.67879441e-01)
( 2.12000000e-01,  4.60980286e-01)
( 2.24000000e-01,  5.61243736e-01)
( 2.36000000e-01,  6.63915763e-01)
( 2.48000000e-01,  7.63074204e-01)
( 2.60000000e-01,  8.52143789e-01)
( 2.72000000e-01,  9.24594515e-01)
( 2.84000000e-01,  9.74724902e-01)
( 2.96000000e-01,  9.98401279e-01)
( 3.08000000e-01,  9.93620436e-01)
( 3.20000000e-01,  9.60789439e-01)
( 3.32000000e-01,  9.02668412e-01)
( 3.44000000e-01,  8.23987433e-01)
( 3.56000000e-01,  7.30811294e-01)
( 3.68000000e-01,  6.29770381e-01)
( 3.80000000e-01,  5.27292424e-01)
( 3.92000000e-01,  4.28956399e-01)
( 4.04000000e-01,  3.39052607e-01)
( 4.16000000e-01,  2.60383431e-01)
( 4.28000000e-01,  1.94290659e-01)
( 4.40000000e-01,  1.40858421e-01)
( 4.52000000e-01,  9.92215550e-02)
( 4.64000000e-01,  6.79080972e-02)
( 4.76000000e-01,  4.51574503e-02)
( 4.88000000e-01,  2.91762575e-02)
( 5.00000000e-01,  1.83156389e-02)
( 5.12000000e-01,  1.11713814e-02)
( 5.24000000e-01,  6.62039661e-03)
( 5.36000000e-01,  3.81200493e-03)
( 5.48000000e-01,  2.13262855e-03)
( 5.60000000e-01,  1.15922917e-03)
( 5.72000000e-01,  6.12231548e-04)
( 5.84000000e-01,  3.14162560e-04)
( 5.96000000e-01,  1.56633794e-04)
( 6.08000000e-01,  7.58767686e-05)
( 6.20000000e-01,  3.57128496e-05)
( 6.32000000e-01,  1.63317417e-05)
( 6.44000000e-01,  7.25659389e-06)
( 6.56000000e-01,  3.13274795e-06)
( 6.68000000e-01,  1.31404551e-06)
( 6.80000000e-01,  5.35534780e-07)
( 6.92000000e-01,  2.12059275e-07)
( 7.04000000e-01,  8.15866656e-08)
( 7.16000000e-01,  3.04981442e-08)
( 7.28000000e-01,  1.10769441e-08)
( 7.40000000e-01,  3.90893843e-09)
( 7.52000000e-01,  1.34026296e-09)
( 7.64000000e-01,  4.46491836e-10)
( 7.76000000e-01,  1.44520464e-10)
( 7.88000000e-01,  4.54503770e-11)
( 8.00000000e-01,  1.38879439e-11)
( 8.12000000e-01,  4.12316518e-12)
( 8.24000000e-01,  1.18936686e-12)
( 8.36000000e-01,  3.33344468e-13)
( 8.48000000e-01,  9.07743172e-14)
( 8.60000000e-01,  2.40173478e-14)
( 8.72000000e-01,  6.17418194e-15)
( 8.84000000e-01,  1.54214822e-15)
( 8.96000000e-01,  3.74252863e-16)
( 9.08000000e-01,  8.82462916e-17)
( 9.20000000e-01,  2.02171585e-17)
( 9.32000000e-01,  4.50024425e-18)
( 9.44000000e-01,  9.73294738e-19)
( 9.56000000e-01,  2.04524332e-19)
( 9.68000000e-01,  4.17578289e-20)
( 9.80000000e-01,  8.28367701e-21)
( 9.92000000e-01,  1.59661690e-21)};\label{line:bump_1d_rom_ref}

\addplot [green, dashdotted, , mark repeat=0]
coordinates {
(-1.00000000e+00, -9.30503331e-20)
(-9.88000000e-01, -2.51780584e-16)
(-9.76000000e-01, -4.45668292e-18)
(-9.64000000e-01, -1.75035681e-17)
(-9.52000000e-01, -6.67934319e-17)
(-9.40000000e-01, -2.47647100e-16)
(-9.28000000e-01, -8.92123656e-16)
(-9.16000000e-01, -3.12254842e-15)
(-9.04000000e-01, -1.06190491e-14)
(-8.92000000e-01, -3.50876589e-14)
(-8.80000000e-01, -1.12645923e-13)
(-8.68000000e-01, -3.51373398e-13)
(-8.56000000e-01, -1.06491416e-12)
(-8.44000000e-01, -3.13583093e-12)
(-8.32000000e-01, -8.97187055e-12)
(-8.20000000e-01, -2.49405299e-11)
(-8.08000000e-01, -6.73628838e-11)
(-7.96000000e-01, -1.76777906e-10)
(-7.84000000e-01, -4.50741571e-10)
(-7.72000000e-01, -1.11665640e-09)
(-7.60000000e-01, -2.68784253e-09)
(-7.48000000e-01, -6.28608707e-09)
(-7.36000000e-01, -1.42839812e-08)
(-7.24000000e-01, -3.15362787e-08)
(-7.12000000e-01, -6.76493935e-08)
(-7.00000000e-01, -1.40996923e-07)
(-6.88000000e-01, -2.85527271e-07)
(-6.76000000e-01, -5.61794824e-07)
(-6.64000000e-01, -1.07398933e-06)
(-6.52000000e-01, -1.99486795e-06)
(-6.40000000e-01, -3.60014563e-06)
(-6.28000000e-01, -6.31273199e-06)
(-6.16000000e-01, -1.07548762e-05)
(-6.04000000e-01, -1.78025948e-05)
(-5.92000000e-01, -2.86318559e-05)
(-5.80000000e-01, -4.47405702e-05)
(-5.68000000e-01, -6.79258915e-05)
(-5.56000000e-01, -1.00194648e-04)
(-5.44000000e-01, -1.43588146e-04)
(-5.32000000e-01, -1.99912808e-04)
(-5.20000000e-01, -2.70385410e-04)
(-5.08000000e-01, -3.55223918e-04)
(-4.96000000e-01, -4.53237374e-04)
(-4.84000000e-01, -5.61484019e-04)
(-4.72000000e-01, -6.75068693e-04)
(-4.60000000e-01, -7.87133494e-04)
(-4.48000000e-01, -8.89059832e-04)
(-4.36000000e-01, -9.70852186e-04)
(-4.24000000e-01, -1.02162746e-03)
(-4.12000000e-01, -1.03010580e-03)
(-4.00000000e-01, -9.85004165e-04)
(-3.88000000e-01, -8.75280235e-04)
(-3.76000000e-01, -6.90255292e-04)
(-3.64000000e-01, -4.19739831e-04)
(-3.52000000e-01, -5.43632097e-05)
(-3.40000000e-01,  4.13666943e-04)
(-3.28000000e-01,  9.89204120e-04)
(-3.16000000e-01,  1.67216796e-03)
(-3.04000000e-01,  2.45508046e-03)
(-2.92000000e-01,  3.32050651e-03)
(-2.80000000e-01,  4.23891924e-03)
(-2.68000000e-01,  5.16753245e-03)
(-2.56000000e-01,  6.05050047e-03)
(-2.44000000e-01,  6.82059217e-03)
(-2.32000000e-01,  7.40207071e-03)
(-2.20000000e-01,  7.71417289e-03)
(-2.08000000e-01,  7.67441675e-03)
(-1.96000000e-01,  7.20108017e-03)
(-1.84000000e-01,  6.21461645e-03)
(-1.72000000e-01,  4.63842610e-03)
(-1.60000000e-01,  2.40010092e-03)
(-1.48000000e-01, -5.65260565e-04)
(-1.36000000e-01, -4.30804917e-03)
(-1.24000000e-01, -8.85243449e-03)
(-1.12000000e-01, -1.41801187e-02)
(-1.00000000e-01, -2.02112552e-02)
(-8.80000000e-02, -2.67856640e-02)
(-7.60000000e-02, -3.36482879e-02)
(-6.40000000e-02, -4.04425889e-02)
(-5.20000000e-02, -4.67141217e-02)
(-4.00000000e-02, -5.19240908e-02)
(-2.80000000e-02, -5.54699680e-02)
(-1.60000000e-02, -5.67082009e-02)
(-4.00000000e-03, -5.49736659e-02)
( 8.00000000e-03, -4.95924201e-02)
( 2.00000000e-02, -3.98884266e-02)
( 3.20000000e-02, -2.51902795e-02)
( 4.40000000e-02, -4.84879797e-03)
( 5.60000000e-02,  2.17214516e-02)
( 6.80000000e-02,  5.49659750e-02)
( 8.00000000e-02,  9.50764497e-02)
( 9.20000000e-02,  1.41848273e-01)
( 1.04000000e-01,  1.94532899e-01)
( 1.16000000e-01,  2.51720389e-01)
( 1.28000000e-01,  3.11293110e-01)
( 1.40000000e-01,  3.70486981e-01)
( 1.52000000e-01,  4.26080841e-01)
( 1.64000000e-01,  4.74709529e-01)
( 1.76000000e-01,  5.13267229e-01)
( 1.88000000e-01,  5.39341828e-01)
( 2.00000000e-01,  5.51605530e-01)
( 2.12000000e-01,  5.50086872e-01)
( 2.24000000e-01,  5.36265572e-01)
( 2.36000000e-01,  5.12961123e-01)
( 2.48000000e-01,  4.84022055e-01)
( 2.60000000e-01,  4.53857331e-01)
( 2.72000000e-01,  4.26877377e-01)
( 2.84000000e-01,  4.06925524e-01)
( 2.96000000e-01,  3.96780221e-01)
( 3.08000000e-01,  3.97796185e-01)
( 3.20000000e-01,  4.09732149e-01)
( 3.32000000e-01,  4.30787436e-01)
( 3.44000000e-01,  4.57842346e-01)
( 3.56000000e-01,  4.86870638e-01)
( 3.68000000e-01,  5.13468696e-01)
( 3.80000000e-01,  5.33428323e-01)
( 3.92000000e-01,  5.43271932e-01)
( 4.04000000e-01,  5.40673002e-01)
( 4.16000000e-01,  5.24702159e-01)
( 4.28000000e-01,  4.95868140e-01)
( 4.40000000e-01,  4.55958027e-01)
( 4.52000000e-01,  4.07715300e-01)
( 4.64000000e-01,  3.54419856e-01)
( 4.76000000e-01,  2.99445687e-01)
( 4.88000000e-01,  2.45867575e-01)
( 5.00000000e-01,  1.96170199e-01)
( 5.12000000e-01,  1.52086952e-01)
( 5.24000000e-01,  1.14568376e-01)
( 5.36000000e-01,  8.38577526e-02)
( 5.48000000e-01,  5.96378881e-02)
( 5.60000000e-01,  4.12096306e-02)
( 5.72000000e-01,  2.76675477e-02)
( 5.84000000e-01,  1.80483244e-02)
( 5.96000000e-01,  1.14392260e-02)
( 6.08000000e-01,  7.04448842e-03)
( 6.20000000e-01,  4.21497556e-03)
( 6.32000000e-01,  2.45037904e-03)
( 6.44000000e-01,  1.38408886e-03)
( 6.56000000e-01,  7.59603759e-04)
( 6.68000000e-01,  4.05044393e-04)
( 6.80000000e-01,  2.09850740e-04)
( 6.92000000e-01,  1.05635704e-04)
( 7.04000000e-01,  5.16658157e-05)
( 7.16000000e-01,  2.45520738e-05)
( 7.28000000e-01,  1.13361451e-05)
( 7.40000000e-01,  5.08551490e-06)
( 7.52000000e-01,  2.21664829e-06)
( 7.64000000e-01,  9.38752202e-07)
( 7.76000000e-01,  3.86275770e-07)
( 7.88000000e-01,  1.54431641e-07)
( 8.00000000e-01,  5.99884199e-08)
( 8.12000000e-01,  2.26407536e-08)
( 8.24000000e-01,  8.30245738e-09)
( 8.36000000e-01,  2.95811293e-09)
( 8.48000000e-01,  1.02403585e-09)
( 8.60000000e-01,  3.44435491e-10)
( 8.72000000e-01,  1.12562297e-10)
( 8.84000000e-01,  3.57412977e-11)
( 8.96000000e-01,  1.10265590e-11)
( 9.08000000e-01,  3.30523250e-12)
( 9.20000000e-01,  9.62623238e-13)
( 9.32000000e-01,  2.72397419e-13)
( 9.44000000e-01,  7.48931250e-14)
( 9.56000000e-01,  2.00065976e-14)
( 9.68000000e-01,  5.19274312e-15)
( 9.80000000e-01,  1.30952189e-15)
( 9.92000000e-01,  3.20864006e-16)};\label{line:bump_1d_rom_r}

\addplot [red, dashed, , mark repeat=0]
coordinates {
(-1.00000000e+00,  1.26453838e-15)
(-9.80843200e-01, -7.15719777e-32)
(-9.61772800e-01, -3.70271825e-31)
(-9.42788800e-01, -1.63647239e-31)
(-9.23891200e-01, -1.00848165e-30)
(-9.05080000e-01, -1.51985410e-30)
(-8.86355200e-01, -5.18816942e-31)
(-8.67716800e-01, -1.09388241e-30)
(-8.49164800e-01, -1.42356749e-29)
(-8.30699200e-01, -1.03543740e-29)
(-8.12320000e-01, -1.02737029e-29)
(-7.94027200e-01, -8.98103927e-29)
(-7.75820800e-01, -7.98404473e-29)
(-7.57700800e-01, -3.32783416e-28)
(-7.39667200e-01, -4.40826794e-28)
(-7.21720000e-01, -2.03222080e-27)
(-7.03859200e-01, -3.64298293e-27)
(-6.86084800e-01, -3.57125824e-27)
(-6.68396800e-01, -1.03996900e-26)
(-6.50795200e-01, -1.00014176e-26)
(-6.33280000e-01,  8.57435832e-26)
(-6.15851200e-01,  4.89167441e-25)
(-5.98508800e-01,  2.93733079e-24)
(-5.81252800e-01,  1.71928544e-23)
(-5.64083200e-01,  9.61121595e-23)
(-5.47000000e-01,  5.23151394e-22)
(-5.30003200e-01,  2.77090667e-21)
(-5.13092800e-01,  1.42463765e-20)
(-4.96268800e-01,  7.11398859e-20)
(-4.79531200e-01,  3.45122920e-19)
(-4.62880000e-01,  1.62664553e-18)
(-4.46315200e-01,  7.44935080e-18)
(-4.29836800e-01,  3.31463106e-17)
(-4.13444800e-01,  1.43298764e-16)
(-3.97139200e-01,  6.01924274e-16)
(-3.80920000e-01,  2.45659519e-15)
(-3.64787200e-01,  9.74131778e-15)
(-3.48740800e-01,  3.75313424e-14)
(-3.32780800e-01,  1.40495631e-13)
(-3.16907200e-01,  5.11003476e-13)
(-3.01120000e-01,  1.80583144e-12)
(-2.85419200e-01,  6.20044532e-12)
(-2.69804800e-01,  2.06852548e-11)
(-2.54276800e-01,  6.70488271e-11)
(-2.38835200e-01,  2.11161051e-10)
(-2.23480000e-01,  6.46143177e-10)
(-2.08211200e-01,  1.92103833e-09)
(-1.93028800e-01,  5.54926694e-09)
(-1.77932800e-01,  1.55749809e-08)
(-1.62923200e-01,  4.24728835e-08)
(-1.48000000e-01,  1.12535175e-07)
(-1.33163200e-01,  2.89705774e-07)
(-1.18412800e-01,  7.24633375e-07)
(-1.03748800e-01,  1.76105059e-06)
(-8.91712000e-02,  4.15831755e-06)
(-7.46800000e-02,  9.54016287e-06)
(-6.02752000e-02,  2.12660214e-05)
(-4.59568000e-02,  4.60584185e-05)
(-3.17248000e-02,  9.69223903e-05)
(-1.75792000e-02,  1.98167089e-04)
(-3.52000000e-03,  3.93669041e-04)
( 1.04528000e-02,  7.59842031e-04)
( 2.43392000e-02,  1.42497644e-03)
( 3.81392000e-02,  2.59647636e-03)
( 5.18528000e-02,  4.59677642e-03)
( 6.54800000e-02,  7.90705405e-03)
( 7.90208000e-02,  1.32150361e-02)
( 9.24752000e-02,  2.14592391e-02)
( 1.05843200e-01,  3.38573219e-02)
( 1.19124800e-01,  5.19018938e-02)
( 1.32320000e-01,  7.73047404e-02)
( 1.45428800e-01,  1.11871991e-01)
( 1.58451200e-01,  1.57300074e-01)
( 1.71387200e-01,  2.14896234e-01)
( 1.84236800e-01,  2.85246945e-01)
( 1.97000000e-01,  3.67879441e-01)
( 2.09676800e-01,  4.60980286e-01)
( 2.22267200e-01,  5.61243736e-01)
( 2.34771200e-01,  6.63915763e-01)
( 2.47188800e-01,  7.63074204e-01)
( 2.59520000e-01,  8.52143789e-01)
( 2.71764800e-01,  9.24594515e-01)
( 2.83923200e-01,  9.74724902e-01)
( 2.95995200e-01,  9.98401279e-01)
( 3.07980800e-01,  9.93620436e-01)
( 3.19880000e-01,  9.60789439e-01)
( 3.31692800e-01,  9.02668412e-01)
( 3.43419200e-01,  8.23987433e-01)
( 3.55059200e-01,  7.30811294e-01)
( 3.66612800e-01,  6.29770381e-01)
( 3.78080000e-01,  5.27292424e-01)
( 3.89460800e-01,  4.28956399e-01)
( 4.00755200e-01,  3.39052607e-01)
( 4.11963200e-01,  2.60383431e-01)
( 4.23084800e-01,  1.94290659e-01)
( 4.34120000e-01,  1.40858421e-01)
( 4.45068800e-01,  9.92215550e-02)
( 4.55931200e-01,  6.79080972e-02)
( 4.66707200e-01,  4.51574503e-02)
( 4.77396800e-01,  2.91762575e-02)
( 4.88000000e-01,  1.83156389e-02)
( 4.98516800e-01,  1.11713814e-02)
( 5.08947200e-01,  6.62039661e-03)
( 5.19291200e-01,  3.81200493e-03)
( 5.29548800e-01,  2.13262855e-03)
( 5.39720000e-01,  1.15922917e-03)
( 5.49804800e-01,  6.12231548e-04)
( 5.59803200e-01,  3.14162560e-04)
( 5.69715200e-01,  1.56633794e-04)
( 5.79540800e-01,  7.58767686e-05)
( 5.89280000e-01,  3.57128496e-05)
( 5.98932800e-01,  1.63317417e-05)
( 6.08499200e-01,  7.25659389e-06)
( 6.17979200e-01,  3.13274795e-06)
( 6.27372800e-01,  1.31404551e-06)
( 6.36680000e-01,  5.35534780e-07)
( 6.45900800e-01,  2.12059275e-07)
( 6.55035200e-01,  8.15866656e-08)
( 6.64083200e-01,  3.04981442e-08)
( 6.73044800e-01,  1.10769441e-08)
( 6.81920000e-01,  3.90893843e-09)
( 6.90708800e-01,  1.34026296e-09)
( 6.99411200e-01,  4.46491836e-10)
( 7.08027200e-01,  1.44520464e-10)
( 7.16556800e-01,  4.54503771e-11)
( 7.25000000e-01,  1.38879439e-11)
( 7.33356800e-01,  4.12316519e-12)
( 7.41627200e-01,  1.18936687e-12)
( 7.49811200e-01,  3.33344472e-13)
( 7.57908800e-01,  9.07743196e-14)
( 7.65920000e-01,  2.40173493e-14)
( 7.73844800e-01,  6.17418287e-15)
( 7.81683200e-01,  1.54214876e-15)
( 7.89435200e-01,  3.74253181e-16)
( 7.97100800e-01,  8.82464741e-17)
( 8.04680000e-01,  2.02172729e-17)
( 8.12172800e-01,  4.50030849e-18)
( 8.19579200e-01,  9.73332563e-19)
( 8.26899200e-01,  2.04545640e-19)
( 8.34132800e-01,  4.17700188e-20)
( 8.41280000e-01,  8.29042275e-21)
( 8.48340800e-01,  1.60032462e-21)
( 8.55315200e-01,  3.01048148e-22)
( 8.62203200e-01,  5.56380423e-23)
( 8.69004800e-01,  1.02442395e-23)
( 8.75720000e-01,  1.99890827e-24)
( 8.82348800e-01,  4.69672742e-25)
( 8.88891200e-01,  1.43624904e-25)
( 8.95347200e-01,  6.04284686e-26)
( 9.01716800e-01,  2.84387360e-26)
( 9.08000000e-01,  1.57151306e-26)
( 9.14196800e-01,  8.21872101e-27)
( 9.20307200e-01,  4.06284771e-27)
( 9.26331200e-01,  2.31088749e-27)
( 9.32268800e-01,  1.11176896e-27)
( 9.38120000e-01,  5.76736259e-28)
( 9.43884800e-01,  2.98512646e-28)
( 9.49563200e-01,  1.50129293e-28)
( 9.55155200e-01,  7.17750495e-29)
( 9.60660800e-01,  3.82619657e-29)
( 9.66080000e-01,  1.89693519e-29)
( 9.71412800e-01,  9.14711779e-30)
( 9.76659200e-01,  4.33434425e-30)
( 9.81819200e-01,  2.22755491e-30)
( 9.86892800e-01,  1.07137407e-30)
( 9.91880000e-01,  4.96824528e-31)
( 9.96780800e-01,  2.54343252e-31)};\label{line:bump_1d_rom_rft}

\end{axis}
\end{tikzpicture}

%% file: figs/methodology/bump1d_motivation/bump_1d_sig_val.tikz
\begin{tikzpicture}
\begin{semilogyaxis}[
width=1.\textwidth,
xlabel=index,
ymax=1,
xmax=5,
ylabel=normalized singular value,
xmin=1,
ymode=log,
ymin=0.001]
\addplot [blue, mark options={solid, thin}, mark=o, mark size=3, mark repeat=1]
coordinates {
( 1.00000000e+00,  1.00000000e+00)
( 2.00000000e+00,  9.53028404e-02)
( 3.00000000e+00,  2.51527567e-02)
( 4.00000000e+00,  4.54846794e-03)
( 5.00000000e+00,  4.40169942e-03)};\label{line:bump_1d_sig_val_r}

\addplot [red, mark options={solid, thin}, mark=square, mark size=3, mark repeat=1]
coordinates {
( 1.00000000e+00,  1.00000000e+00)
( 2.00000000e+00,  9.58808470e-01)
( 3.00000000e+00,  8.99714485e-01)
( 4.00000000e+00,  8.36726608e-01)
( 5.00000000e+00,  7.87610740e-01)};\label{line:bump_1d_sig_val_rft}

\end{semilogyaxis}
\end{tikzpicture}

%% file: figs/methodology/density_sampling/density_sampl.tikz
\begin{tikzpicture}
\begin{axis}[
axis equal image,
width=1.\textwidth,
ytick={-.3,0,.3},
ymax=0.3,
xmax=3,
xmin=-3,
ymin=-0.3]
\addplot []
graphics [xmin=-3,xmax=3,ymin=-0.3,ymax=0.3] { 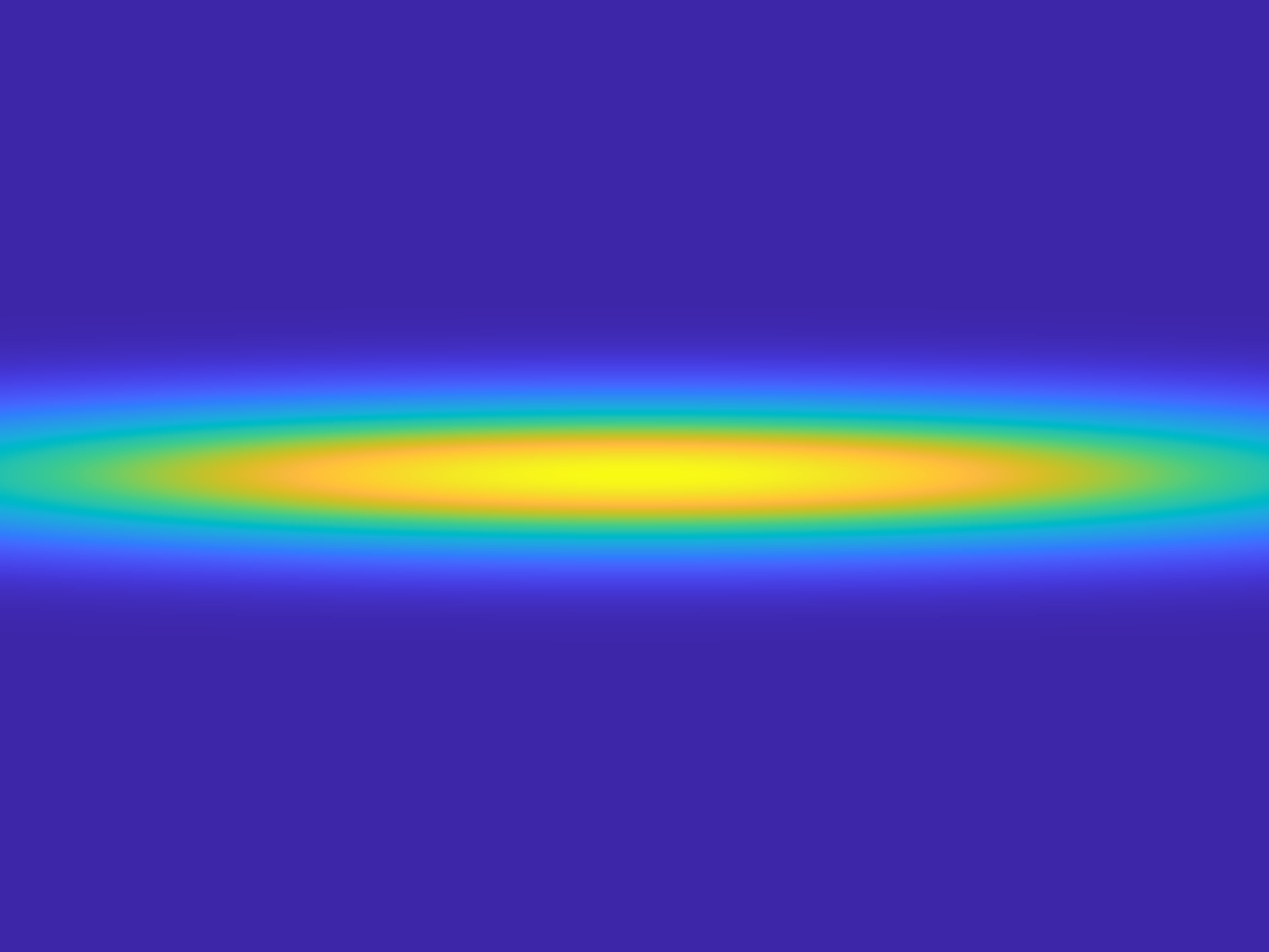};

\addplot [black, only marks, mark options={solid, thin}, mark=x, mark size=2, mark repeat=0]
coordinates {
( 2.48025514e+00,  2.64718492e-02)
(-6.46637883e-01,  3.10955780e-02)
(-2.79332352e+00, -1.22511281e-02)
(-6.14136253e-02, -1.08827599e-02)
( 1.07821606e+00,  3.10196008e-02)
(-9.57685640e-01,  1.70535502e-02)
( 6.96268057e-01, -5.34223022e-03)
( 1.98497177e+00,  1.70528182e-02)
( 1.52237457e+00, -2.39108306e-02)
( 4.12941965e-01, -6.12187179e-03)
(-1.13270975e+00,  5.70662710e-03)
(-2.36008338e+00,  9.23796162e-02)
(-2.49338493e+00, -2.00434702e-02)
(-1.90891783e+00, -4.72394167e-02)
( 7.32330789e-01, -2.98095238e-02)
( 2.66872314e+00, -1.82718151e-03)
(-1.54985228e+00, -1.92175709e-02)
(-8.81048573e-01,  6.42388080e-02)
( 1.39033431e+00,  2.95491926e-02)
(-1.86626991e+00,  3.73550867e-02)
(-7.92502056e-02, -1.28282823e-02)
(-1.16190317e+00,  1.70173108e-03)
( 3.00938057e-01,  2.44950172e-02)
( 2.87849027e+00, -1.22260054e-02)
(-1.45161182e+00, -1.82560308e-02)
(-1.08733019e+00, -1.51666481e-02)
(-6.83461572e-02,  1.57050122e-02)
( 7.44360529e-01,  3.58271082e-02)
( 1.83293655e+00,  1.53443031e-02)
(-1.56040794e+00,  7.73023866e-02)
( 1.27616683e+00,  9.43248310e-05)
(-2.64228679e+00,  3.63943808e-02)
(-2.57132721e+00,  4.32996849e-03)
(-1.24809552e+00, -1.36697660e-02)
( 1.18863312e+00,  3.33055827e-02)
( 2.13313684e+00,  2.89529074e-02)
(-1.85445783e+00, -1.43494014e-02)
(-2.27633032e+00,  1.79014969e-02)
(-6.92285254e-01,  1.65972765e-02)
(-1.25735601e+00,  2.34181769e-02)
(-9.36737975e-01,  1.68138667e-02)
(-1.43563201e+00,  1.88712501e-02)
(-4.48444079e-01, -3.74562226e-02)
(-1.92740288e+00, -1.54228622e-02)
( 5.91142013e-01, -5.81514873e-03)
( 1.19932710e+00,  2.77061517e-02)
( 1.67281345e+00, -1.53094162e-02)
(-3.59489166e-01,  5.42854835e-03)
(-2.95970811e+00,  2.04340975e-02)
(-4.53905761e-01, -7.81672679e-03)
( 2.79558688e-01, -1.48542316e-02)
(-2.99173818e-01, -8.25490127e-03)
(-1.46135405e+00,  2.26921474e-02)
( 8.73311250e-01, -4.10735501e-03)
( 2.68296663e-01,  2.94622961e-02)
(-2.61845177e+00, -1.90840008e-02)
(-2.43707984e+00,  5.08088077e-03)
( 2.16683887e+00, -3.02933329e-03)
(-1.54328785e+00, -1.15195374e-02)
( 1.53046259e+00, -2.45208910e-02)
( 1.07191383e+00, -9.64596182e-04)
(-2.99650501e-02, -7.04783556e-02)
( 2.27408343e+00,  9.77823232e-02)
( 1.66080416e-01, -4.09532296e-03)
(-1.74356950e+00,  1.04582683e-02)
(-2.80805391e+00,  2.29426838e-02)
(-2.70280453e+00, -2.08600216e-03)
(-1.30879886e+00,  7.71933561e-03)
(-5.30391910e-03,  7.16021115e-03)
( 2.24356444e+00, -4.59411335e-02)
( 3.89877424e-01,  2.80623650e-02)
( 7.25751864e-01,  1.47419530e-02)
(-1.19216361e+00, -4.08932331e-02)
(-1.97590878e-01,  2.96396813e-02)
(-9.12724834e-01, -1.07946704e-02)
(-1.93735480e+00,  3.25616124e-02)
(-1.27290393e+00, -1.70954922e-02)
( 1.58374247e+00,  6.36408078e-02)
(-1.93129828e+00, -2.80730173e-02)
( 1.31314042e-01, -3.28302051e-02)
(-1.89483541e+00,  1.94422701e-02)
(-2.57128312e+00, -5.15026882e-02)
(-1.00744300e+00, -4.05306368e-02)
( 1.56855532e+00,  2.62139998e-02)
( 6.20807903e-01,  5.22049316e-03)
( 1.15519192e+00,  1.13339670e-02)
( 6.47195444e-01,  4.82508099e-02)
(-2.23266972e+00,  9.90802140e-03)
(-2.20750227e+00,  4.45449079e-02)
(-2.29504289e+00,  2.81435846e-02)
(-1.76595287e+00, -2.23456738e-02)
(-1.62628049e+00,  2.83881241e-02)
(-2.08892685e+00,  5.63863933e-02)
(-2.45100761e+00, -1.89369160e-02)
(-4.22613244e-01, -9.65215302e-03)
(-2.34709237e+00,  2.63532747e-02)
(-1.89339799e+00, -5.75938315e-02)
(-8.95339531e-01,  3.71071417e-02)
(-4.64126081e-01, -2.80787364e-02)
( 1.45527219e+00, -1.51330433e-02)
( 2.61438524e+00, -8.42273323e-03)
(-2.20436528e-01, -5.75673589e-02)
( 2.45304487e-01,  3.59467796e-02)
( 1.85522311e+00,  4.97237744e-02)
(-6.06715486e-01, -1.69813227e-02)
( 2.11358334e+00,  1.12732351e-03)
( 2.70536649e+00, -1.12071690e-02)
( 2.20049938e+00,  2.62377469e-02)
( 6.29943851e-01, -2.25509137e-02)
(-2.84919009e+00, -1.57775492e-02)
( 1.40537815e+00,  1.42051746e-02)
(-1.65737757e+00, -2.52872385e-02)
( 8.40699289e-01, -6.38766224e-02)
(-6.96788394e-01,  2.54693005e-02)
( 4.52969158e-01,  6.01034095e-03)
(-1.50822624e+00, -9.67224591e-03)
(-9.79804414e-01,  3.24763721e-02)
( 5.01199046e-01,  1.03585030e-02)
(-5.19436266e-01, -5.64535863e-02)
( 2.53399199e+00,  5.41908441e-02)
(-1.35046814e-01,  2.47432825e-02)
(-2.40542831e+00, -2.04723962e-03)
(-1.92610512e+00, -3.22088644e-02)
(-2.39079667e+00, -2.18290495e-02)
( 9.75923006e-01, -4.36996882e-02)
( 1.26677131e+00,  2.49145834e-02)
( 2.26847147e-02, -2.08113226e-03)
( 1.00272180e+00,  1.72879229e-02)
( 2.58024376e+00, -2.01960062e-02)
(-2.67124926e+00, -3.92677241e-02)
( 2.26679443e+00,  1.64865929e-02)
( 2.03381758e+00, -1.33478878e-02)
( 2.33266408e-02,  2.93619333e-02)
(-2.16765218e+00, -4.88541324e-03)
( 1.03935592e+00,  3.28559809e-02)
(-6.59840923e-01, -4.19411407e-04)
(-1.04374697e+00, -8.71507983e-03)
( 2.30643027e+00,  4.41711342e-02)
( 1.04865880e+00, -1.22982352e-02)
(-1.52263131e+00, -3.14573558e-02)
( 2.79322758e-01,  1.23840309e-02)
(-6.11214723e-01,  3.07344384e-03)
(-1.01590760e+00, -1.39996403e-02)
( 2.38719411e+00,  1.86723721e-02)
(-1.78754942e+00, -9.22130686e-03)
( 1.38832240e+00, -2.79937916e-02)
( 1.70539559e+00,  4.11143716e-02)
(-6.60416057e-01,  1.81809461e-02)
( 1.01425554e+00,  4.22648562e-05)
(-2.35189984e+00,  3.39935162e-03)
( 3.07986108e-02, -4.57156751e-02)
( 4.70929850e-02,  1.71218251e-02)
( 3.23322395e-01,  3.60131060e-02)
(-5.78228023e-01, -2.19648124e-02)
(-2.15846784e+00, -4.79739612e-02)
(-2.55546254e+00, -2.12233146e-02)
(-2.14509567e+00, -4.63848002e-02)
(-7.43933304e-01,  4.75600739e-03)
(-2.58985668e+00, -1.27345845e-02)
( 2.76935144e+00,  5.24828968e-02)
( 2.13826145e+00, -1.95132900e-02)
(-3.45612642e-01,  8.08710956e-02)
( 1.94558895e-01,  4.32994693e-02)
(-5.76859733e-01,  9.71325997e-03)
( 3.16392799e-01, -4.50377190e-02)
(-1.94512963e+00, -2.79258046e-02)
( 7.51531081e-01,  8.61243507e-03)
(-1.27543640e+00,  3.31821349e-04)
( 1.16642346e+00, -4.86430873e-02)
( 2.43865960e+00, -2.14630848e-02)
( 2.59738737e+00, -2.18664927e-02)
(-2.08831752e+00, -2.05782315e-02)
(-2.21331176e+00, -1.29918564e-02)
(-1.32674633e+00,  4.93234444e-02)
( 1.26245211e+00, -3.76280110e-02)
(-1.85408136e+00, -1.14940076e-02)
( 1.95944387e+00,  3.53742187e-02)
(-8.14278783e-01,  6.46998700e-03)
(-5.16681013e-01, -3.81727147e-02)
( 2.60987452e+00, -4.10309087e-03)
(-4.08116335e-01,  8.53339749e-03)
( 2.28202716e+00, -5.77962700e-03)
( 2.06692438e+00,  2.30650194e-02)
(-2.44206644e+00, -7.30215045e-03)
( 2.49015547e+00,  2.85483478e-02)
(-2.45339895e+00,  1.88074063e-02)
(-2.85820520e+00,  2.14865221e-02)
( 2.87391280e+00, -4.33464203e-02)
(-9.22435700e-01, -1.62749171e-02)
(-1.47311494e+00, -3.51691590e-02)
(-5.61760072e-01, -2.27617723e-02)
(-1.99865566e+00, -6.23815854e-02)
( 1.86176874e+00,  1.13996523e-02)
( 2.57690683e+00, -1.80970133e-02)
( 1.49105374e+00,  8.65988662e-03)
( 2.35360313e+00, -5.14793227e-02)
(-1.64959271e+00, -2.99972162e-02)
( 1.75548998e+00, -3.41917608e-02)
( 1.96325451e+00,  3.51723238e-02)
(-1.45286241e-01, -2.01849547e-02)
( 2.04651798e+00, -2.90987511e-02)
( 1.54730490e+00, -2.21742473e-02)
(-1.34192810e+00,  2.44647230e-02)
(-2.10920986e+00,  2.39631855e-02)
(-1.17689026e+00, -3.34108645e-03)
( 1.32205924e+00, -3.19229664e-03)
( 9.75096591e-01,  4.66266368e-03)
( 2.77196260e+00,  8.04080750e-03)
( 1.17788678e+00,  3.94323158e-03)
( 2.34021740e+00, -3.39595515e-02)
(-2.31630815e+00, -3.78154573e-02)
(-2.69612101e+00, -6.75963324e-03)
( 6.18937923e-01,  1.99758192e-02)
(-2.78745924e+00,  2.76296701e-03)
(-2.35172378e+00, -8.02487421e-03)
(-2.88752346e-01, -3.33143627e-02)
( 1.44543182e+00,  1.35890481e-03)
(-7.89894950e-01,  8.83635590e-02)
( 9.03045812e-01,  4.53259065e-02)
( 1.73534798e-01,  3.88701121e-02)
( 1.70912755e+00, -4.58336996e-02)
( 4.24097103e-01,  1.43659272e-02)
(-3.22706330e-01, -6.86751833e-03)
(-2.12560740e+00,  1.70087231e-02)
(-8.09804248e-01, -3.81700762e-02)
( 4.98797809e-01,  3.96507853e-02)
(-2.18395567e+00,  5.77782639e-02)
(-1.57278725e+00, -5.12704175e-02)
( 9.03183747e-01,  3.25411315e-03)
( 1.79022155e+00,  3.40210473e-04)
( 2.52474065e+00, -3.54414109e-04)
(-1.21587424e+00, -2.07162929e-02)
(-1.13114785e+00,  3.87686342e-02)
(-5.87468286e-01, -4.09638392e-02)
(-2.36663322e+00,  1.87655209e-02)
( 3.04865250e-01,  7.41804383e-02)
( 2.83939053e-02, -1.90083040e-02)
( 4.51101675e-01,  2.12435865e-02)
( 1.16923927e+00, -1.77156345e-02)
(-2.50010433e+00,  2.22127419e-03)
( 1.92597887e+00,  2.75795240e-02)
( 2.37572777e+00,  3.07501142e-03)
( 1.40263743e+00, -1.25655126e-02)
(-2.05915267e+00, -3.47931432e-02)
(-5.46786552e-01,  4.16021788e-02)
(-4.82893776e-01, -2.83743294e-02)
(-1.47595290e+00, -1.37563105e-02)
( 5.05207347e-01, -2.52842435e-02)
(-1.68603529e+00,  4.44648861e-03)
(-1.80099525e+00,  2.14679844e-02)
( 1.58872534e+00, -1.57861295e-02)
( 1.66195743e+00,  3.80683831e-04)
( 1.89261094e+00, -9.97284342e-02)
(-2.47518710e+00, -4.78545298e-02)
( 2.15688064e+00, -3.19043378e-02)
( 1.17414556e+00,  2.55808356e-02)
(-1.32102648e+00, -1.05985381e-02)
( 2.26580483e+00, -6.17989595e-03)
(-2.00927926e-01, -3.79102575e-04)
(-2.91544155e+00,  2.45760689e-02)
( 5.30198675e-01, -1.33130320e-02)
( 2.31741067e+00, -5.06117270e-02)
( 2.01242730e+00,  6.71426696e-02)
( 1.49093946e+00,  1.35682299e-02)
( 1.23664557e+00,  6.62719513e-02)
( 8.29462257e-01, -3.13988376e-02)
( 4.94588572e-01, -1.57988437e-02)
(-2.85583605e+00, -1.77083920e-03)
(-9.61457015e-01, -4.25300781e-02)
(-6.04420114e-01,  3.95299137e-02)
( 9.97958931e-01, -1.13867821e-02)
( 1.63324039e+00, -4.92914001e-03)
(-4.98391470e-01, -2.39701594e-02)
( 2.07782183e+00, -1.83873011e-02)
( 8.78076393e-01, -2.40853029e-02)
(-1.70592741e-01,  8.59833897e-03)
( 6.52552244e-01, -4.92501586e-02)
( 2.39865611e+00, -5.64269752e-02)
(-5.17380365e-01,  5.49252770e-04)
(-1.56073516e+00,  1.95309033e-02)
(-2.60856117e+00, -3.28177657e-02)
( 1.96857362e+00,  1.48707546e-03)
(-1.64016169e+00,  6.96634192e-03)
( 2.39363842e+00, -2.92721996e-02)
( 1.44315954e+00, -8.69497454e-03)
( 1.19547645e+00,  1.42714425e-02)
(-1.94897887e-01, -7.17327943e-02)
( 9.14525013e-01,  6.20976767e-03)
( 2.88637381e-02, -2.40012830e-03)
( 2.61585332e+00, -2.21436634e-02)
( 5.04617099e-01, -1.94096845e-02)
(-6.62477946e-01,  5.50609707e-02)
( 2.53628058e+00,  1.28521839e-02)
(-9.72932608e-01,  4.41446856e-02)
(-1.56678781e-01, -6.92536312e-04)
( 2.95553891e+00, -8.83125959e-03)
( 7.12465880e-02,  1.03923693e-02)
( 5.26923758e-01, -7.14473549e-02)
( 1.09984749e+00,  2.17118957e-02)
(-5.62320346e-01,  2.59815890e-02)
(-2.23452563e+00, -6.61603337e-02)
(-4.91025173e-01, -2.30286422e-03)
(-1.49816468e+00,  6.22693231e-02)
(-2.15002212e+00, -1.24211484e-02)
(-1.29015017e-01,  1.74806800e-02)
( 2.43198364e+00,  2.80388497e-02)
( 1.11412282e+00, -6.75614300e-03)
(-2.67841311e+00,  5.40499054e-03)
(-1.87570267e+00,  6.33791718e-03)
( 5.97182223e-02,  2.29807130e-02)
(-6.11968881e-01,  3.93137379e-02)};\label{line:density_sampl_sampl}

\addplot [red, only marks, mark options={solid, thin}, mark=square*, mark size=2, mark repeat=0]
coordinates {
(-5.28753632e-01,  6.65183740e-04)
( 1.42586637e+00,  4.04070365e-03)
(-2.59269459e+00, -1.09013502e-03)
(-9.99966406e-01,  4.62330015e-03)
(-2.04544539e+00, -6.78484361e-03)
(-1.49053644e+00,  1.67393507e-03)
( 7.10800424e-01, -2.51451575e-03)
( 2.35401518e+00, -6.94789252e-04)
( 4.41430260e-02,  1.97254216e-03)};\label{line:density_sampl_C}

\end{axis}
\end{tikzpicture}

%% file: figs/methodology/parabola_correspondence/parabola_correspondence.tikz
\begin{tikzpicture}
\begin{axis}[
axis lines=none,
width=.45\textwidth,
ymax=1.2,
xmax=1,
xmin=-1,
ymin=-0.6]
\addplot [blue, , mark repeat=0]
coordinates {
(-1.00000000e+00,  1.10000000e+00)
(-9.90000000e-01,  1.07214000e+00)
(-9.80000000e-01,  1.04456000e+00)
(-9.70000000e-01,  1.01726000e+00)
(-9.60000000e-01,  9.90240000e-01)
(-9.50000000e-01,  9.63500000e-01)
(-9.40000000e-01,  9.37040000e-01)
(-9.30000000e-01,  9.10860000e-01)
(-9.20000000e-01,  8.84960000e-01)
(-9.10000000e-01,  8.59340000e-01)
(-9.00000000e-01,  8.34000000e-01)
(-8.90000000e-01,  8.08940000e-01)
(-8.80000000e-01,  7.84160000e-01)
(-8.70000000e-01,  7.59660000e-01)
(-8.60000000e-01,  7.35440000e-01)
(-8.50000000e-01,  7.11500000e-01)
(-8.40000000e-01,  6.87840000e-01)
(-8.30000000e-01,  6.64460000e-01)
(-8.20000000e-01,  6.41360000e-01)
(-8.10000000e-01,  6.18540000e-01)
(-8.00000000e-01,  5.96000000e-01)
(-7.90000000e-01,  5.73740000e-01)
(-7.80000000e-01,  5.51760000e-01)
(-7.70000000e-01,  5.30060000e-01)
(-7.60000000e-01,  5.08640000e-01)
(-7.50000000e-01,  4.87500000e-01)
(-7.40000000e-01,  4.66640000e-01)
(-7.30000000e-01,  4.46060000e-01)
(-7.20000000e-01,  4.25760000e-01)
(-7.10000000e-01,  4.05740000e-01)
(-7.00000000e-01,  3.86000000e-01)
(-6.90000000e-01,  3.66540000e-01)
(-6.80000000e-01,  3.47360000e-01)
(-6.70000000e-01,  3.28460000e-01)
(-6.60000000e-01,  3.09840000e-01)
(-6.50000000e-01,  2.91500000e-01)
(-6.40000000e-01,  2.73440000e-01)
(-6.30000000e-01,  2.55660000e-01)
(-6.20000000e-01,  2.38160000e-01)
(-6.10000000e-01,  2.20940000e-01)
(-6.00000000e-01,  2.04000000e-01)
(-5.90000000e-01,  1.87340000e-01)
(-5.80000000e-01,  1.70960000e-01)
(-5.70000000e-01,  1.54860000e-01)
(-5.60000000e-01,  1.39040000e-01)
(-5.50000000e-01,  1.23500000e-01)
(-5.40000000e-01,  1.08240000e-01)
(-5.30000000e-01,  9.32600000e-02)
(-5.20000000e-01,  7.85600000e-02)
(-5.10000000e-01,  6.41400000e-02)
(-5.00000000e-01,  5.00000000e-02)
(-4.90000000e-01,  3.61400000e-02)
(-4.80000000e-01,  2.25600000e-02)
(-4.70000000e-01,  9.26000000e-03)
(-4.60000000e-01, -3.76000000e-03)
(-4.50000000e-01, -1.65000000e-02)
(-4.40000000e-01, -2.89600000e-02)
(-4.30000000e-01, -4.11400000e-02)
(-4.20000000e-01, -5.30400000e-02)
(-4.10000000e-01, -6.46600000e-02)
(-4.00000000e-01, -7.60000000e-02)
(-3.90000000e-01, -8.70600000e-02)
(-3.80000000e-01, -9.78400000e-02)
(-3.70000000e-01, -1.08340000e-01)
(-3.60000000e-01, -1.18560000e-01)
(-3.50000000e-01, -1.28500000e-01)
(-3.40000000e-01, -1.38160000e-01)
(-3.30000000e-01, -1.47540000e-01)
(-3.20000000e-01, -1.56640000e-01)
(-3.10000000e-01, -1.65460000e-01)
(-3.00000000e-01, -1.74000000e-01)
(-2.90000000e-01, -1.82260000e-01)
(-2.80000000e-01, -1.90240000e-01)
(-2.70000000e-01, -1.97940000e-01)
(-2.60000000e-01, -2.05360000e-01)
(-2.50000000e-01, -2.12500000e-01)
(-2.40000000e-01, -2.19360000e-01)
(-2.30000000e-01, -2.25940000e-01)
(-2.20000000e-01, -2.32240000e-01)
(-2.10000000e-01, -2.38260000e-01)
(-2.00000000e-01, -2.44000000e-01)
(-1.90000000e-01, -2.49460000e-01)
(-1.80000000e-01, -2.54640000e-01)
(-1.70000000e-01, -2.59540000e-01)
(-1.60000000e-01, -2.64160000e-01)
(-1.50000000e-01, -2.68500000e-01)
(-1.40000000e-01, -2.72560000e-01)
(-1.30000000e-01, -2.76340000e-01)
(-1.20000000e-01, -2.79840000e-01)
(-1.10000000e-01, -2.83060000e-01)
(-1.00000000e-01, -2.86000000e-01)
(-9.00000000e-02, -2.88660000e-01)
(-8.00000000e-02, -2.91040000e-01)
(-7.00000000e-02, -2.93140000e-01)
(-6.00000000e-02, -2.94960000e-01)
(-5.00000000e-02, -2.96500000e-01)
(-4.00000000e-02, -2.97760000e-01)
(-3.00000000e-02, -2.98740000e-01)
(-2.00000000e-02, -2.99440000e-01)
(-1.00000000e-02, -2.99860000e-01)
( 0.00000000e+00, -3.00000000e-01)
( 1.00000000e-02, -2.99860000e-01)
( 2.00000000e-02, -2.99440000e-01)
( 3.00000000e-02, -2.98740000e-01)
( 4.00000000e-02, -2.97760000e-01)
( 5.00000000e-02, -2.96500000e-01)
( 6.00000000e-02, -2.94960000e-01)
( 7.00000000e-02, -2.93140000e-01)
( 8.00000000e-02, -2.91040000e-01)
( 9.00000000e-02, -2.88660000e-01)
( 1.00000000e-01, -2.86000000e-01)
( 1.10000000e-01, -2.83060000e-01)
( 1.20000000e-01, -2.79840000e-01)
( 1.30000000e-01, -2.76340000e-01)
( 1.40000000e-01, -2.72560000e-01)
( 1.50000000e-01, -2.68500000e-01)
( 1.60000000e-01, -2.64160000e-01)
( 1.70000000e-01, -2.59540000e-01)
( 1.80000000e-01, -2.54640000e-01)
( 1.90000000e-01, -2.49460000e-01)
( 2.00000000e-01, -2.44000000e-01)
( 2.10000000e-01, -2.38260000e-01)
( 2.20000000e-01, -2.32240000e-01)
( 2.30000000e-01, -2.25940000e-01)
( 2.40000000e-01, -2.19360000e-01)
( 2.50000000e-01, -2.12500000e-01)
( 2.60000000e-01, -2.05360000e-01)
( 2.70000000e-01, -1.97940000e-01)
( 2.80000000e-01, -1.90240000e-01)
( 2.90000000e-01, -1.82260000e-01)
( 3.00000000e-01, -1.74000000e-01)
( 3.10000000e-01, -1.65460000e-01)
( 3.20000000e-01, -1.56640000e-01)
( 3.30000000e-01, -1.47540000e-01)
( 3.40000000e-01, -1.38160000e-01)
( 3.50000000e-01, -1.28500000e-01)
( 3.60000000e-01, -1.18560000e-01)
( 3.70000000e-01, -1.08340000e-01)
( 3.80000000e-01, -9.78400000e-02)
( 3.90000000e-01, -8.70600000e-02)
( 4.00000000e-01, -7.60000000e-02)
( 4.10000000e-01, -6.46600000e-02)
( 4.20000000e-01, -5.30400000e-02)
( 4.30000000e-01, -4.11400000e-02)
( 4.40000000e-01, -2.89600000e-02)
( 4.50000000e-01, -1.65000000e-02)
( 4.60000000e-01, -3.76000000e-03)
( 4.70000000e-01,  9.26000000e-03)
( 4.80000000e-01,  2.25600000e-02)
( 4.90000000e-01,  3.61400000e-02)
( 5.00000000e-01,  5.00000000e-02)
( 5.10000000e-01,  6.41400000e-02)
( 5.20000000e-01,  7.85600000e-02)
( 5.30000000e-01,  9.32600000e-02)
( 5.40000000e-01,  1.08240000e-01)
( 5.50000000e-01,  1.23500000e-01)
( 5.60000000e-01,  1.39040000e-01)
( 5.70000000e-01,  1.54860000e-01)
( 5.80000000e-01,  1.70960000e-01)
( 5.90000000e-01,  1.87340000e-01)
( 6.00000000e-01,  2.04000000e-01)
( 6.10000000e-01,  2.20940000e-01)
( 6.20000000e-01,  2.38160000e-01)
( 6.30000000e-01,  2.55660000e-01)
( 6.40000000e-01,  2.73440000e-01)
( 6.50000000e-01,  2.91500000e-01)
( 6.60000000e-01,  3.09840000e-01)
( 6.70000000e-01,  3.28460000e-01)
( 6.80000000e-01,  3.47360000e-01)
( 6.90000000e-01,  3.66540000e-01)
( 7.00000000e-01,  3.86000000e-01)
( 7.10000000e-01,  4.05740000e-01)
( 7.20000000e-01,  4.25760000e-01)
( 7.30000000e-01,  4.46060000e-01)
( 7.40000000e-01,  4.66640000e-01)
( 7.50000000e-01,  4.87500000e-01)
( 7.60000000e-01,  5.08640000e-01)
( 7.70000000e-01,  5.30060000e-01)
( 7.80000000e-01,  5.51760000e-01)
( 7.90000000e-01,  5.73740000e-01)
( 8.00000000e-01,  5.96000000e-01)
( 8.10000000e-01,  6.18540000e-01)
( 8.20000000e-01,  6.41360000e-01)
( 8.30000000e-01,  6.64460000e-01)
( 8.40000000e-01,  6.87840000e-01)
( 8.50000000e-01,  7.11500000e-01)
( 8.60000000e-01,  7.35440000e-01)
( 8.70000000e-01,  7.59660000e-01)
( 8.80000000e-01,  7.84160000e-01)
( 8.90000000e-01,  8.08940000e-01)
( 9.00000000e-01,  8.34000000e-01)
( 9.10000000e-01,  8.59340000e-01)
( 9.20000000e-01,  8.84960000e-01)
( 9.30000000e-01,  9.10860000e-01)
( 9.40000000e-01,  9.37040000e-01)
( 9.50000000e-01,  9.63500000e-01)
( 9.60000000e-01,  9.90240000e-01)
( 9.70000000e-01,  1.01726000e+00)
( 9.80000000e-01,  1.04456000e+00)
( 9.90000000e-01,  1.07214000e+00)
( 1.00000000e+00,  1.10000000e+00)};\label{line:par_corres_y1}

\addplot [green, , mark repeat=0]
coordinates {
(-1.00000000e+00,  6.00000000e-01)
(-9.90000000e-01,  5.80100000e-01)
(-9.80000000e-01,  5.60400000e-01)
(-9.70000000e-01,  5.40900000e-01)
(-9.60000000e-01,  5.21600000e-01)
(-9.50000000e-01,  5.02500000e-01)
(-9.40000000e-01,  4.83600000e-01)
(-9.30000000e-01,  4.64900000e-01)
(-9.20000000e-01,  4.46400000e-01)
(-9.10000000e-01,  4.28100000e-01)
(-9.00000000e-01,  4.10000000e-01)
(-8.90000000e-01,  3.92100000e-01)
(-8.80000000e-01,  3.74400000e-01)
(-8.70000000e-01,  3.56900000e-01)
(-8.60000000e-01,  3.39600000e-01)
(-8.50000000e-01,  3.22500000e-01)
(-8.40000000e-01,  3.05600000e-01)
(-8.30000000e-01,  2.88900000e-01)
(-8.20000000e-01,  2.72400000e-01)
(-8.10000000e-01,  2.56100000e-01)
(-8.00000000e-01,  2.40000000e-01)
(-7.90000000e-01,  2.24100000e-01)
(-7.80000000e-01,  2.08400000e-01)
(-7.70000000e-01,  1.92900000e-01)
(-7.60000000e-01,  1.77600000e-01)
(-7.50000000e-01,  1.62500000e-01)
(-7.40000000e-01,  1.47600000e-01)
(-7.30000000e-01,  1.32900000e-01)
(-7.20000000e-01,  1.18400000e-01)
(-7.10000000e-01,  1.04100000e-01)
(-7.00000000e-01,  9.00000000e-02)
(-6.90000000e-01,  7.61000000e-02)
(-6.80000000e-01,  6.24000000e-02)
(-6.70000000e-01,  4.89000000e-02)
(-6.60000000e-01,  3.56000000e-02)
(-6.50000000e-01,  2.25000000e-02)
(-6.40000000e-01,  9.60000000e-03)
(-6.30000000e-01, -3.10000000e-03)
(-6.20000000e-01, -1.56000000e-02)
(-6.10000000e-01, -2.79000000e-02)
(-6.00000000e-01, -4.00000000e-02)
(-5.90000000e-01, -5.19000000e-02)
(-5.80000000e-01, -6.36000000e-02)
(-5.70000000e-01, -7.51000000e-02)
(-5.60000000e-01, -8.64000000e-02)
(-5.50000000e-01, -9.75000000e-02)
(-5.40000000e-01, -1.08400000e-01)
(-5.30000000e-01, -1.19100000e-01)
(-5.20000000e-01, -1.29600000e-01)
(-5.10000000e-01, -1.39900000e-01)
(-5.00000000e-01, -1.50000000e-01)
(-4.90000000e-01, -1.59900000e-01)
(-4.80000000e-01, -1.69600000e-01)
(-4.70000000e-01, -1.79100000e-01)
(-4.60000000e-01, -1.88400000e-01)
(-4.50000000e-01, -1.97500000e-01)
(-4.40000000e-01, -2.06400000e-01)
(-4.30000000e-01, -2.15100000e-01)
(-4.20000000e-01, -2.23600000e-01)
(-4.10000000e-01, -2.31900000e-01)
(-4.00000000e-01, -2.40000000e-01)
(-3.90000000e-01, -2.47900000e-01)
(-3.80000000e-01, -2.55600000e-01)
(-3.70000000e-01, -2.63100000e-01)
(-3.60000000e-01, -2.70400000e-01)
(-3.50000000e-01, -2.77500000e-01)
(-3.40000000e-01, -2.84400000e-01)
(-3.30000000e-01, -2.91100000e-01)
(-3.20000000e-01, -2.97600000e-01)
(-3.10000000e-01, -3.03900000e-01)
(-3.00000000e-01, -3.10000000e-01)
(-2.90000000e-01, -3.15900000e-01)
(-2.80000000e-01, -3.21600000e-01)
(-2.70000000e-01, -3.27100000e-01)
(-2.60000000e-01, -3.32400000e-01)
(-2.50000000e-01, -3.37500000e-01)
(-2.40000000e-01, -3.42400000e-01)
(-2.30000000e-01, -3.47100000e-01)
(-2.20000000e-01, -3.51600000e-01)
(-2.10000000e-01, -3.55900000e-01)
(-2.00000000e-01, -3.60000000e-01)
(-1.90000000e-01, -3.63900000e-01)
(-1.80000000e-01, -3.67600000e-01)
(-1.70000000e-01, -3.71100000e-01)
(-1.60000000e-01, -3.74400000e-01)
(-1.50000000e-01, -3.77500000e-01)
(-1.40000000e-01, -3.80400000e-01)
(-1.30000000e-01, -3.83100000e-01)
(-1.20000000e-01, -3.85600000e-01)
(-1.10000000e-01, -3.87900000e-01)
(-1.00000000e-01, -3.90000000e-01)
(-9.00000000e-02, -3.91900000e-01)
(-8.00000000e-02, -3.93600000e-01)
(-7.00000000e-02, -3.95100000e-01)
(-6.00000000e-02, -3.96400000e-01)
(-5.00000000e-02, -3.97500000e-01)
(-4.00000000e-02, -3.98400000e-01)
(-3.00000000e-02, -3.99100000e-01)
(-2.00000000e-02, -3.99600000e-01)
(-1.00000000e-02, -3.99900000e-01)
( 0.00000000e+00, -4.00000000e-01)
( 1.00000000e-02, -3.99900000e-01)
( 2.00000000e-02, -3.99600000e-01)
( 3.00000000e-02, -3.99100000e-01)
( 4.00000000e-02, -3.98400000e-01)
( 5.00000000e-02, -3.97500000e-01)
( 6.00000000e-02, -3.96400000e-01)
( 7.00000000e-02, -3.95100000e-01)
( 8.00000000e-02, -3.93600000e-01)
( 9.00000000e-02, -3.91900000e-01)
( 1.00000000e-01, -3.90000000e-01)
( 1.10000000e-01, -3.87900000e-01)
( 1.20000000e-01, -3.85600000e-01)
( 1.30000000e-01, -3.83100000e-01)
( 1.40000000e-01, -3.80400000e-01)
( 1.50000000e-01, -3.77500000e-01)
( 1.60000000e-01, -3.74400000e-01)
( 1.70000000e-01, -3.71100000e-01)
( 1.80000000e-01, -3.67600000e-01)
( 1.90000000e-01, -3.63900000e-01)
( 2.00000000e-01, -3.60000000e-01)
( 2.10000000e-01, -3.55900000e-01)
( 2.20000000e-01, -3.51600000e-01)
( 2.30000000e-01, -3.47100000e-01)
( 2.40000000e-01, -3.42400000e-01)
( 2.50000000e-01, -3.37500000e-01)
( 2.60000000e-01, -3.32400000e-01)
( 2.70000000e-01, -3.27100000e-01)
( 2.80000000e-01, -3.21600000e-01)
( 2.90000000e-01, -3.15900000e-01)
( 3.00000000e-01, -3.10000000e-01)
( 3.10000000e-01, -3.03900000e-01)
( 3.20000000e-01, -2.97600000e-01)
( 3.30000000e-01, -2.91100000e-01)
( 3.40000000e-01, -2.84400000e-01)
( 3.50000000e-01, -2.77500000e-01)
( 3.60000000e-01, -2.70400000e-01)
( 3.70000000e-01, -2.63100000e-01)
( 3.80000000e-01, -2.55600000e-01)
( 3.90000000e-01, -2.47900000e-01)
( 4.00000000e-01, -2.40000000e-01)
( 4.10000000e-01, -2.31900000e-01)
( 4.20000000e-01, -2.23600000e-01)
( 4.30000000e-01, -2.15100000e-01)
( 4.40000000e-01, -2.06400000e-01)
( 4.50000000e-01, -1.97500000e-01)
( 4.60000000e-01, -1.88400000e-01)
( 4.70000000e-01, -1.79100000e-01)
( 4.80000000e-01, -1.69600000e-01)
( 4.90000000e-01, -1.59900000e-01)
( 5.00000000e-01, -1.50000000e-01)
( 5.10000000e-01, -1.39900000e-01)
( 5.20000000e-01, -1.29600000e-01)
( 5.30000000e-01, -1.19100000e-01)
( 5.40000000e-01, -1.08400000e-01)
( 5.50000000e-01, -9.75000000e-02)
( 5.60000000e-01, -8.64000000e-02)
( 5.70000000e-01, -7.51000000e-02)
( 5.80000000e-01, -6.36000000e-02)
( 5.90000000e-01, -5.19000000e-02)
( 6.00000000e-01, -4.00000000e-02)
( 6.10000000e-01, -2.79000000e-02)
( 6.20000000e-01, -1.56000000e-02)
( 6.30000000e-01, -3.10000000e-03)
( 6.40000000e-01,  9.60000000e-03)
( 6.50000000e-01,  2.25000000e-02)
( 6.60000000e-01,  3.56000000e-02)
( 6.70000000e-01,  4.89000000e-02)
( 6.80000000e-01,  6.24000000e-02)
( 6.90000000e-01,  7.61000000e-02)
( 7.00000000e-01,  9.00000000e-02)
( 7.10000000e-01,  1.04100000e-01)
( 7.20000000e-01,  1.18400000e-01)
( 7.30000000e-01,  1.32900000e-01)
( 7.40000000e-01,  1.47600000e-01)
( 7.50000000e-01,  1.62500000e-01)
( 7.60000000e-01,  1.77600000e-01)
( 7.70000000e-01,  1.92900000e-01)
( 7.80000000e-01,  2.08400000e-01)
( 7.90000000e-01,  2.24100000e-01)
( 8.00000000e-01,  2.40000000e-01)
( 8.10000000e-01,  2.56100000e-01)
( 8.20000000e-01,  2.72400000e-01)
( 8.30000000e-01,  2.88900000e-01)
( 8.40000000e-01,  3.05600000e-01)
( 8.50000000e-01,  3.22500000e-01)
( 8.60000000e-01,  3.39600000e-01)
( 8.70000000e-01,  3.56900000e-01)
( 8.80000000e-01,  3.74400000e-01)
( 8.90000000e-01,  3.92100000e-01)
( 9.00000000e-01,  4.10000000e-01)
( 9.10000000e-01,  4.28100000e-01)
( 9.20000000e-01,  4.46400000e-01)
( 9.30000000e-01,  4.64900000e-01)
( 9.40000000e-01,  4.83600000e-01)
( 9.50000000e-01,  5.02500000e-01)
( 9.60000000e-01,  5.21600000e-01)
( 9.70000000e-01,  5.40900000e-01)
( 9.80000000e-01,  5.60400000e-01)
( 9.90000000e-01,  5.80100000e-01)
( 1.00000000e+00,  6.00000000e-01)};\label{line:par_corres_y2}

\addplot [red, only marks, mark options={solid, thin}, mark=square*, mark size=2, mark repeat=0]
coordinates {
( 4.98800000e-01,  6.66726933e-02)
( 5.00000000e-04, -2.57999655e-01)
(-8.48300000e-01,  7.18211553e-01)
( 8.48600000e-01,  7.18881789e-01)
(-4.98000000e-01,  6.55933120e-02)};\label{line:par_corres_C1}

\addplot [red, only marks, mark options={solid, thin}, mark=square*, mark size=2, mark repeat=0]
coordinates {
( 4.67800000e-01, -1.66949080e-01)
( 9.00000000e-04, -3.77396473e-01)
(-8.36700000e-01,  3.08966740e-01)
( 8.37200000e-01,  3.09749307e-01)
(-4.66400000e-01, -1.68201840e-01)};\label{line:par_corres_C2}

\addplot [black, dotted, thick, , mark repeat=0]
coordinates {
( 4.98800000e-01,  6.66726933e-02)
( 4.67800000e-01, -1.66949080e-01)};\label{line:par_corres_corres}

\addplot [black, dotted, thick, , mark repeat=0]
coordinates {
( 5.00000000e-04, -2.57999655e-01)
( 9.00000000e-04, -3.77396473e-01)};\label{line:par_corres_corres}

\addplot [black, dotted, thick, , mark repeat=0]
coordinates {
(-8.48300000e-01,  7.18211553e-01)
(-8.36700000e-01,  3.08966740e-01)};\label{line:par_corres_corres}

\addplot [black, dotted, thick, , mark repeat=0]
coordinates {
( 8.48600000e-01,  7.18881789e-01)
( 8.37200000e-01,  3.09749307e-01)};\label{line:par_corres_corres}

\addplot [black, dotted, thick, , mark repeat=0]
coordinates {
(-4.98000000e-01,  6.55933120e-02)
(-4.66400000e-01, -1.68201840e-01)};\label{line:par_corres_corres}

\end{axis}
\end{tikzpicture}

%% file: figs/methodology/arc/arc0.tikz
\begin{tikzpicture}
\begin{axis}[
axis equal image,
xmin=-1,
width=1.\textwidth,
ymin=-1,
ymax=1,
xmax=1]
\addplot []
graphics [xmin=-1,xmax=1,ymin=-1,ymax=1] { 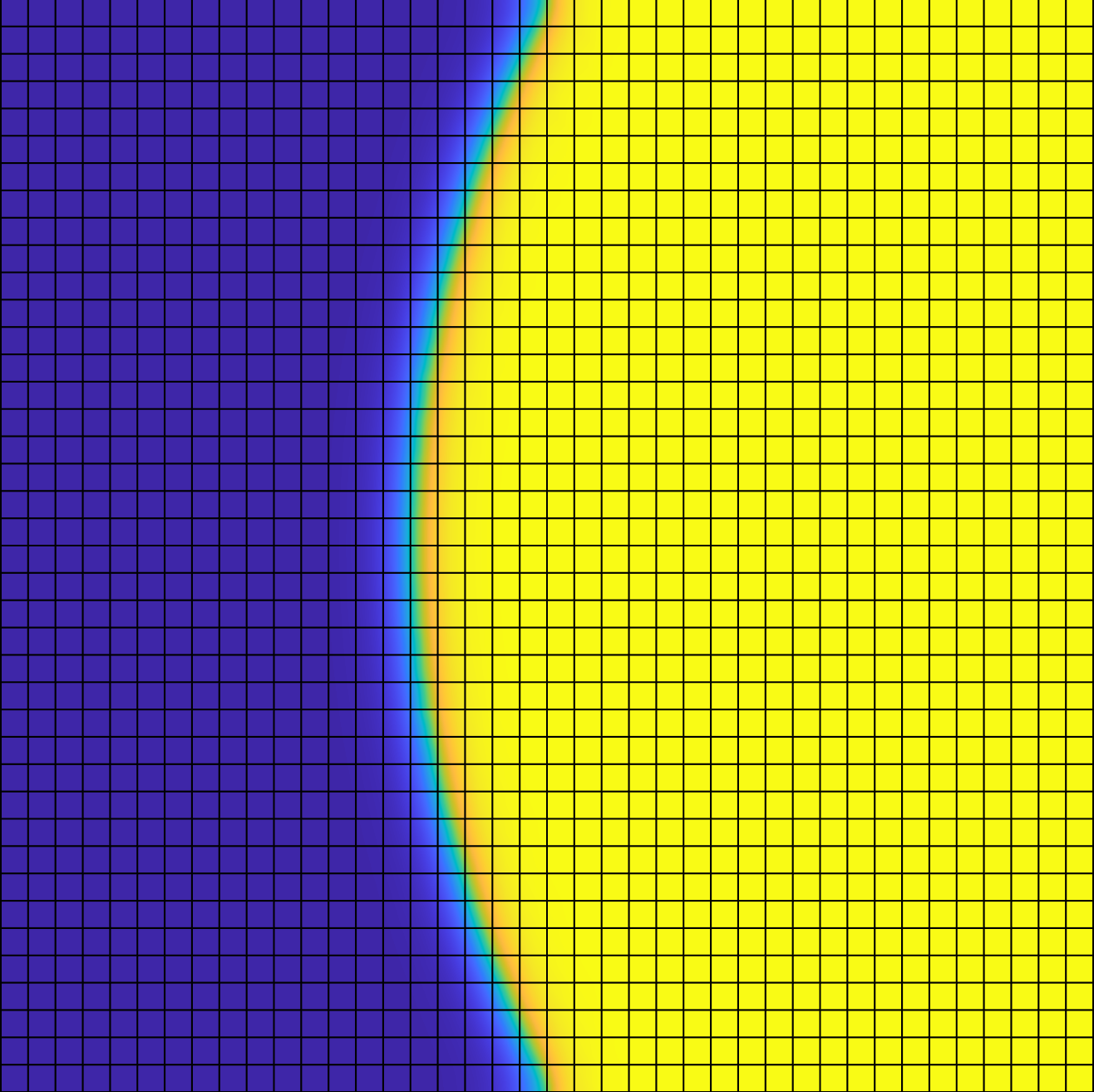};

\end{axis}
\end{tikzpicture}

%% file: figs/methodology/arc/arc1.tikz
\begin{tikzpicture}
\begin{axis}[
axis equal image,
xmin=-1,
width=1.\textwidth,
ymin=-1,
ymax=1,
xmax=1]
\addplot []
graphics [xmin=-1,xmax=1,ymin=-1,ymax=1] { 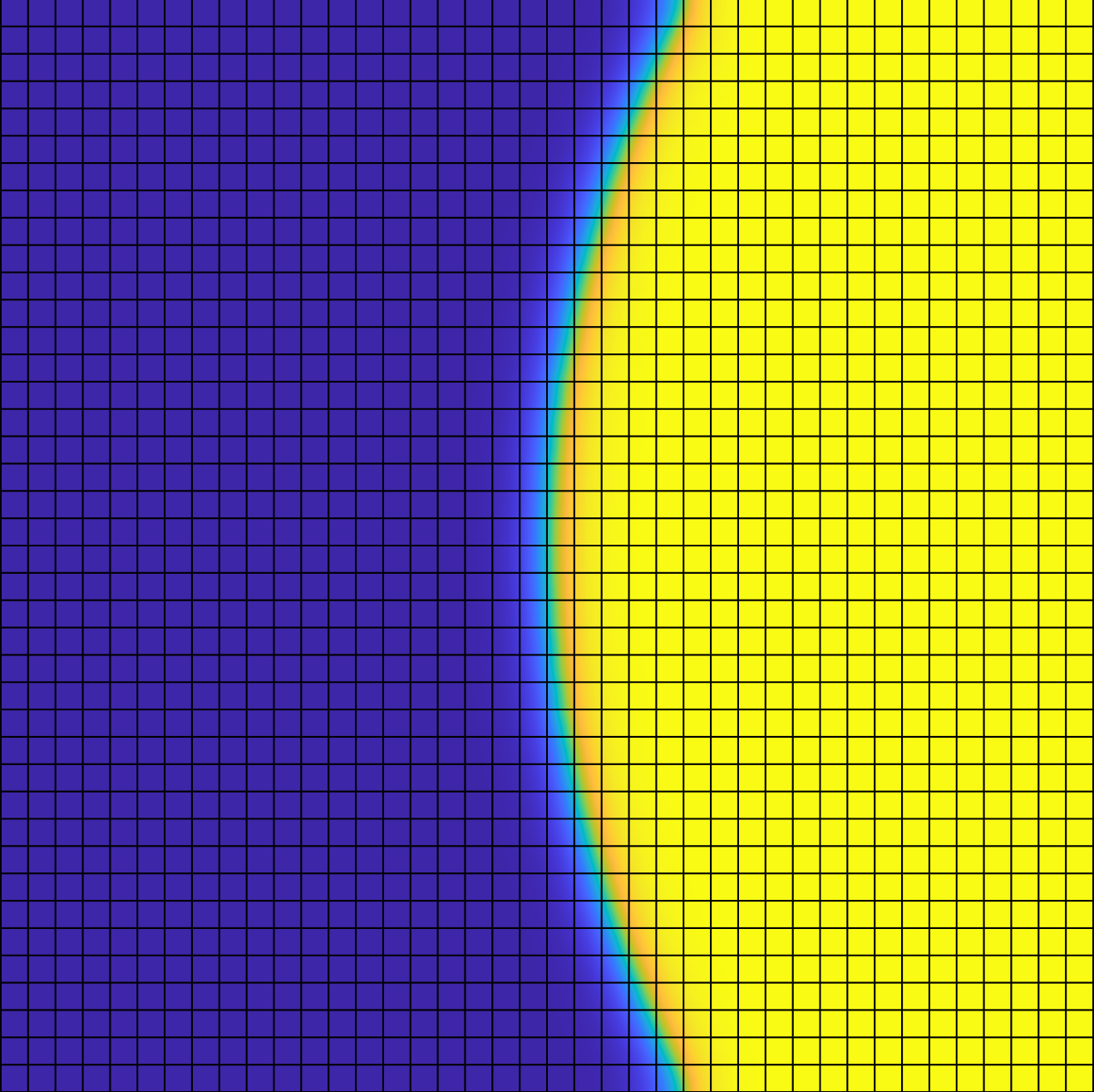};

\end{axis}
\end{tikzpicture}

%% file: figs/methodology/arc/arc0_cntl_pnts.tikz
\begin{tikzpicture}
\begin{axis}[
axis equal image,
xmin=-1,
width=1.\textwidth,
ymin=-1,
ymax=1,
xmax=1]
\addplot []
graphics [xmin=-1,xmax=1,ymin=-1,ymax=1] { arc0.png};

\addplot [red, only marks, mark options={solid, thin}, mark=square*, mark size=2, mark repeat=0, forget plot]
coordinates {
(-1.00000000e+00, -1.00000000e+00)
(-1.00000000e+00,  1.00000000e+00)
( 1.00000000e+00, -1.00000000e+00)
( 1.00000000e+00,  1.00000000e+00)};

\addplot [red, only marks, mark options={solid, thin}, mark=square*, mark size=2, mark repeat=0, forget plot]
coordinates {
( 1.00000000e+00, -8.00000000e-01)
( 1.00000000e+00, -6.00000000e-01)
( 1.00000000e+00, -4.00000000e-01)
( 1.00000000e+00, -2.00000000e-01)
( 1.00000000e+00, -2.22044605e-16)
( 1.00000000e+00,  2.00000000e-01)
( 1.00000000e+00,  4.00000000e-01)
( 1.00000000e+00,  6.00000000e-01)
( 1.00000000e+00,  8.00000000e-01)};

\addplot [red, only marks, mark options={solid, thin}, mark=square*, mark size=2, mark repeat=0, forget plot]
coordinates {
(-1.00000000e+00, -8.00000000e-01)
(-1.00000000e+00, -6.00000000e-01)
(-1.00000000e+00, -4.00000000e-01)
(-1.00000000e+00, -2.00000000e-01)
(-1.00000000e+00, -2.22044605e-16)
(-1.00000000e+00,  2.00000000e-01)
(-1.00000000e+00,  4.00000000e-01)
(-1.00000000e+00,  6.00000000e-01)
(-1.00000000e+00,  8.00000000e-01)};

\addplot [red, only marks, mark options={solid, thin}, mark=square*, mark size=2, mark repeat=0, forget plot]
coordinates {
(-8.00000000e-01, -1.00000000e+00)
(-6.00000000e-01, -1.00000000e+00)
(-4.00000000e-01, -1.00000000e+00)
(-2.00000000e-01, -1.00000000e+00)
(-2.22044605e-16, -1.00000000e+00)
( 2.00000000e-01, -1.00000000e+00)
( 4.00000000e-01, -1.00000000e+00)
( 6.00000000e-01, -1.00000000e+00)
( 8.00000000e-01, -1.00000000e+00)};

\addplot [red, only marks, mark options={solid, thin}, mark=square*, mark size=2, mark repeat=0]
coordinates {
(-8.00000000e-01,  1.00000000e+00)
(-6.00000000e-01,  1.00000000e+00)
(-4.00000000e-01,  1.00000000e+00)
(-2.00000000e-01,  1.00000000e+00)
(-2.22044605e-16,  1.00000000e+00)
( 2.00000000e-01,  1.00000000e+00)
( 4.00000000e-01,  1.00000000e+00)
( 6.00000000e-01,  1.00000000e+00)
( 8.00000000e-01,  1.00000000e+00)};\label{line:arc0_bnd}

\addplot [black, only marks, mark options={solid, thin}, mark=square*, mark size=2, mark repeat=0]
coordinates {
( 0.00000000e+00, -1.00000000e+00)
( 0.00000000e+00,  1.00000000e+00)};\label{line:arc0_end}

\addplot [violet, only marks, mark options={solid, thin}, mark=square*, mark size=2, mark repeat=0]
coordinates {
(-1.09375000e-01, -7.50000000e-01)
(-1.87500000e-01, -5.00000000e-01)
(-2.34375000e-01, -2.50000000e-01)
(-2.50000000e-01,  0.00000000e+00)
(-2.34375000e-01,  2.50000000e-01)
(-1.87500000e-01,  5.00000000e-01)
(-1.09375000e-01,  7.50000000e-01)};\label{line:arc0_arc}

\end{axis}
\end{tikzpicture}

%% file: figs/methodology/arc/arc0_match_cntl_pnts.tikz
\begin{tikzpicture}
\begin{axis}[
axis equal image,
xmin=-1,
width=1.\textwidth,
ymin=-1,
ymax=1,
xmax=1]
\addplot []
graphics [xmin=-1,xmax=1,ymin=-1,ymax=1] { 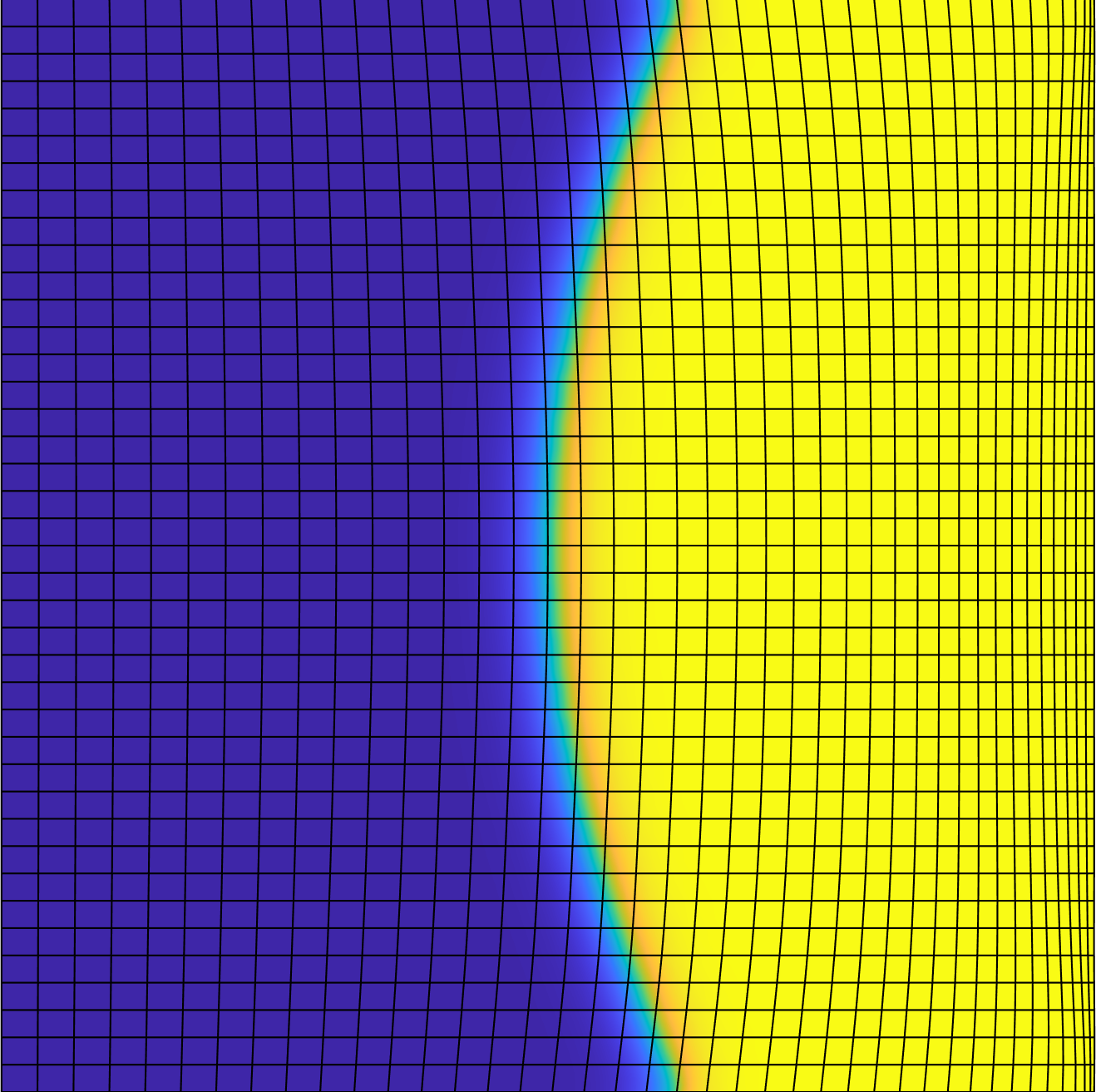};

\addplot [red, only marks, mark options={solid, thin}, mark=square*, mark size=2, mark repeat=0, forget plot]
coordinates {
(-1.00000000e+00, -1.00000000e+00)
(-1.00000000e+00,  1.00000000e+00)
( 1.00000000e+00, -1.00000000e+00)
( 1.00000000e+00,  1.00000000e+00)};

\addplot [red, only marks, mark options={solid, thin}, mark=square*, mark size=2, mark repeat=0, forget plot]
coordinates {
( 1.00000000e+00, -8.00000000e-01)
( 1.00000000e+00, -6.00000000e-01)
( 1.00000000e+00, -4.00000000e-01)
( 1.00000000e+00, -2.00000000e-01)
( 1.00000000e+00, -2.22044605e-16)
( 1.00000000e+00,  2.00000000e-01)
( 1.00000000e+00,  4.00000000e-01)
( 1.00000000e+00,  6.00000000e-01)
( 1.00000000e+00,  8.00000000e-01)};

\addplot [red, only marks, mark options={solid, thin}, mark=square*, mark size=2, mark repeat=0, forget plot]
coordinates {
(-1.00000000e+00, -8.00000000e-01)
(-1.00000000e+00, -6.00000000e-01)
(-1.00000000e+00, -4.00000000e-01)
(-1.00000000e+00, -2.00000000e-01)
(-1.00000000e+00, -2.22044605e-16)
(-1.00000000e+00,  2.00000000e-01)
(-1.00000000e+00,  4.00000000e-01)
(-1.00000000e+00,  6.00000000e-01)
(-1.00000000e+00,  8.00000000e-01)};

\addplot [red, only marks, mark options={solid, thin}, mark=square*, mark size=2, mark repeat=0, forget plot]
coordinates {
(-7.10000000e-01, -1.00000000e+00)
(-4.40000000e-01, -1.00000000e+00)
(-1.90000000e-01, -1.00000000e+00)
( 4.00000000e-02, -1.00000000e+00)
( 2.50000000e-01, -1.00000000e+00)
( 4.40000000e-01, -1.00000000e+00)
( 6.10000000e-01, -1.00000000e+00)
( 7.60000000e-01, -1.00000000e+00)
( 8.90000000e-01, -1.00000000e+00)};

\addplot [red, only marks, mark options={solid, thin}, mark=square*, mark size=2, mark repeat=0]
coordinates {
(-7.10000000e-01,  1.00000000e+00)
(-4.40000000e-01,  1.00000000e+00)
(-1.90000000e-01,  1.00000000e+00)
( 4.00000000e-02,  1.00000000e+00)
( 2.50000000e-01,  1.00000000e+00)
( 4.40000000e-01,  1.00000000e+00)
( 6.10000000e-01,  1.00000000e+00)
( 7.60000000e-01,  1.00000000e+00)
( 8.90000000e-01,  1.00000000e+00)};\label{line:arc0_match_bnd}

\addplot [black, only marks, mark options={solid, thin}, mark=square*, mark size=2, mark repeat=0]
coordinates {
( 2.50000000e-01, -1.00000000e+00)
( 2.50000000e-01,  1.00000000e+00)};\label{line:arc0_match_end}

\addplot [violet, only marks, mark options={solid, thin}, mark=square*, mark size=2, mark repeat=0]
coordinates {
( 1.40625000e-01, -7.50000000e-01)
( 6.25000000e-02, -5.00000000e-01)
( 1.56250000e-02, -2.50000000e-01)
( 0.00000000e+00,  0.00000000e+00)
( 1.56250000e-02,  2.50000000e-01)
( 6.25000000e-02,  5.00000000e-01)
( 1.40625000e-01,  7.50000000e-01)};\label{line:arc0_match_arc}

\end{axis}
\end{tikzpicture}

%% file: figs/burgers/control_points/burgers_xc.tikz
\begin{tikzpicture}
\begin{axis}[
axis equal image,
axis lines=none,
width=.8\textwidth,
xtick={-1,-.5,0,.5,1},
ymax=1,
xmax=1,
xmin=-1,
ymin=0]
\addplot []
graphics [xmin=-1,xmax=1,ymin=0,ymax=1] { 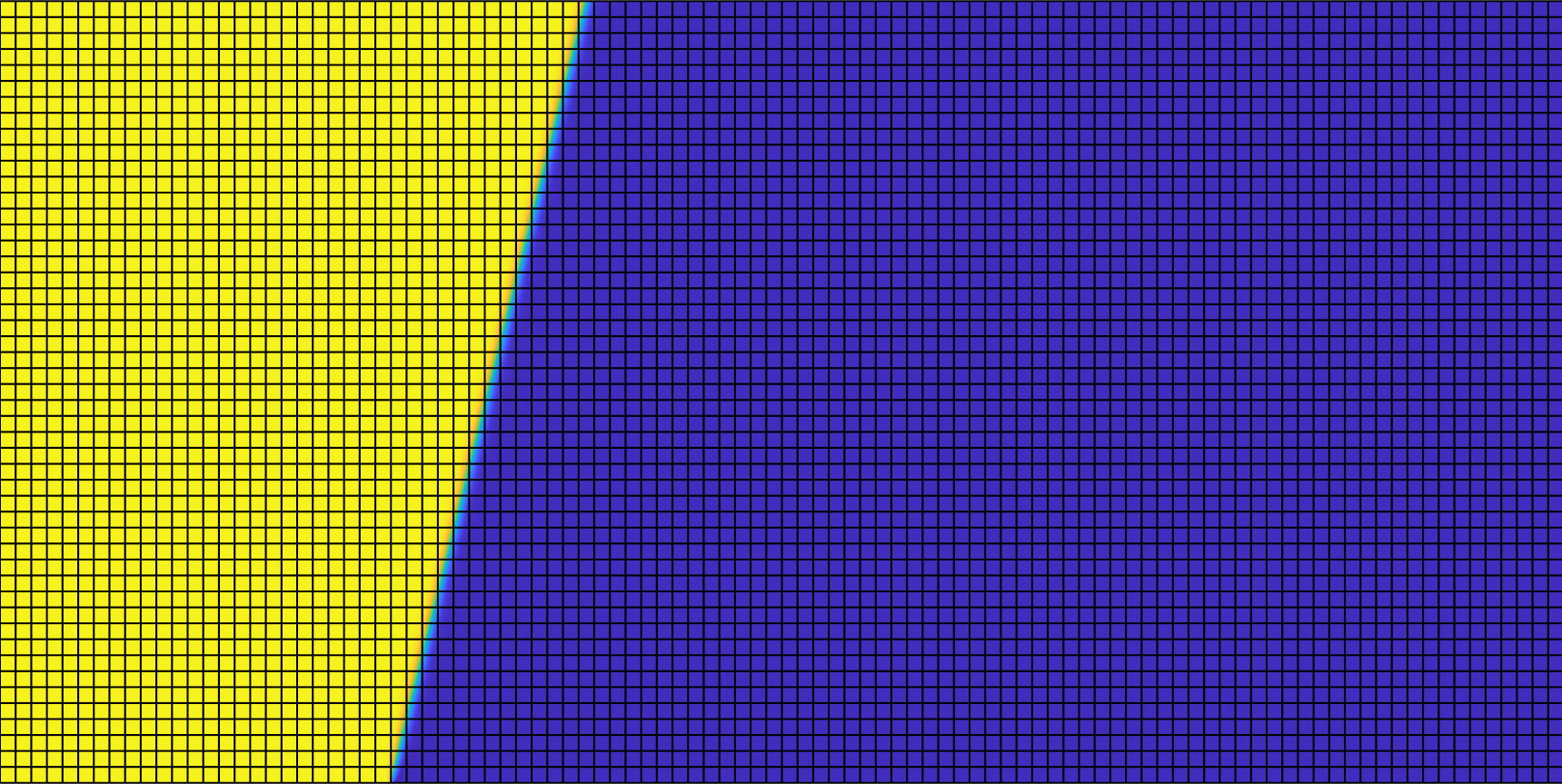};

\addplot [red, only marks, mark options={solid, thin}, mark=square*, mark size=0.5, mark repeat=0]
coordinates {
(-1.00000000e+00,  0.00000000e+00)
(-1.00000000e+00,  1.00000000e+00)
( 1.00000000e+00,  0.00000000e+00)
( 1.00000000e+00,  1.00000000e+00)
(-1.00000000e+00,  2.00000000e-02)
(-1.00000000e+00,  4.00000000e-02)
(-1.00000000e+00,  6.00000000e-02)
(-1.00000000e+00,  8.00000000e-02)
(-1.00000000e+00,  1.00000000e-01)
(-1.00000000e+00,  1.20000000e-01)
(-1.00000000e+00,  1.40000000e-01)
(-1.00000000e+00,  1.60000000e-01)
(-1.00000000e+00,  1.80000000e-01)
(-1.00000000e+00,  2.00000000e-01)
(-1.00000000e+00,  2.20000000e-01)
(-1.00000000e+00,  2.40000000e-01)
(-1.00000000e+00,  2.60000000e-01)
(-1.00000000e+00,  2.80000000e-01)
(-1.00000000e+00,  3.00000000e-01)
(-1.00000000e+00,  3.20000000e-01)
(-1.00000000e+00,  3.40000000e-01)
(-1.00000000e+00,  3.60000000e-01)
(-1.00000000e+00,  3.80000000e-01)
(-1.00000000e+00,  4.00000000e-01)
(-1.00000000e+00,  4.20000000e-01)
(-1.00000000e+00,  4.40000000e-01)
(-1.00000000e+00,  4.60000000e-01)
(-1.00000000e+00,  4.80000000e-01)
(-1.00000000e+00,  5.00000000e-01)
(-1.00000000e+00,  5.20000000e-01)
(-1.00000000e+00,  5.40000000e-01)
(-1.00000000e+00,  5.60000000e-01)
(-1.00000000e+00,  5.80000000e-01)
(-1.00000000e+00,  6.00000000e-01)
(-1.00000000e+00,  6.20000000e-01)
(-1.00000000e+00,  6.40000000e-01)
(-1.00000000e+00,  6.60000000e-01)
(-1.00000000e+00,  6.80000000e-01)
(-1.00000000e+00,  7.00000000e-01)
(-1.00000000e+00,  7.20000000e-01)
(-1.00000000e+00,  7.40000000e-01)
(-1.00000000e+00,  7.60000000e-01)
(-1.00000000e+00,  7.80000000e-01)
(-1.00000000e+00,  8.00000000e-01)
(-1.00000000e+00,  8.20000000e-01)
(-1.00000000e+00,  8.40000000e-01)
(-1.00000000e+00,  8.60000000e-01)
(-1.00000000e+00,  8.80000000e-01)
(-1.00000000e+00,  9.00000000e-01)
(-1.00000000e+00,  9.20000000e-01)
(-1.00000000e+00,  9.40000000e-01)
(-1.00000000e+00,  9.60000000e-01)
(-1.00000000e+00,  9.80000000e-01)
( 1.00000000e+00,  2.00000000e-02)
( 1.00000000e+00,  4.00000000e-02)
( 1.00000000e+00,  6.00000000e-02)
( 1.00000000e+00,  8.00000000e-02)
( 1.00000000e+00,  1.00000000e-01)
( 1.00000000e+00,  1.20000000e-01)
( 1.00000000e+00,  1.40000000e-01)
( 1.00000000e+00,  1.60000000e-01)
( 1.00000000e+00,  1.80000000e-01)
( 1.00000000e+00,  2.00000000e-01)
( 1.00000000e+00,  2.20000000e-01)
( 1.00000000e+00,  2.40000000e-01)
( 1.00000000e+00,  2.60000000e-01)
( 1.00000000e+00,  2.80000000e-01)
( 1.00000000e+00,  3.00000000e-01)
( 1.00000000e+00,  3.20000000e-01)
( 1.00000000e+00,  3.40000000e-01)
( 1.00000000e+00,  3.60000000e-01)
( 1.00000000e+00,  3.80000000e-01)
( 1.00000000e+00,  4.00000000e-01)
( 1.00000000e+00,  4.20000000e-01)
( 1.00000000e+00,  4.40000000e-01)
( 1.00000000e+00,  4.60000000e-01)
( 1.00000000e+00,  4.80000000e-01)
( 1.00000000e+00,  5.00000000e-01)
( 1.00000000e+00,  5.20000000e-01)
( 1.00000000e+00,  5.40000000e-01)
( 1.00000000e+00,  5.60000000e-01)
( 1.00000000e+00,  5.80000000e-01)
( 1.00000000e+00,  6.00000000e-01)
( 1.00000000e+00,  6.20000000e-01)
( 1.00000000e+00,  6.40000000e-01)
( 1.00000000e+00,  6.60000000e-01)
( 1.00000000e+00,  6.80000000e-01)
( 1.00000000e+00,  7.00000000e-01)
( 1.00000000e+00,  7.20000000e-01)
( 1.00000000e+00,  7.40000000e-01)
( 1.00000000e+00,  7.60000000e-01)
( 1.00000000e+00,  7.80000000e-01)
( 1.00000000e+00,  8.00000000e-01)
( 1.00000000e+00,  8.20000000e-01)
( 1.00000000e+00,  8.40000000e-01)
( 1.00000000e+00,  8.60000000e-01)
( 1.00000000e+00,  8.80000000e-01)
( 1.00000000e+00,  9.00000000e-01)
( 1.00000000e+00,  9.20000000e-01)
( 1.00000000e+00,  9.40000000e-01)
( 1.00000000e+00,  9.60000000e-01)
( 1.00000000e+00,  9.80000000e-01)
(-9.60000000e-01,  0.00000000e+00)
(-9.20000000e-01,  0.00000000e+00)
(-8.80000000e-01,  0.00000000e+00)
(-8.40000000e-01,  0.00000000e+00)
(-8.00000000e-01,  0.00000000e+00)
(-7.60000000e-01,  0.00000000e+00)
(-7.20000000e-01,  0.00000000e+00)
(-6.80000000e-01,  0.00000000e+00)
(-6.40000000e-01,  0.00000000e+00)
(-6.00000000e-01,  0.00000000e+00)
(-5.60000000e-01,  0.00000000e+00)
(-5.20000000e-01,  0.00000000e+00)
(-4.80000000e-01,  0.00000000e+00)
(-4.40000000e-01,  0.00000000e+00)
(-4.00000000e-01,  0.00000000e+00)
(-3.60000000e-01,  0.00000000e+00)
(-3.20000000e-01,  0.00000000e+00)
(-2.80000000e-01,  0.00000000e+00)
(-2.40000000e-01,  0.00000000e+00)
(-2.00000000e-01,  0.00000000e+00)
(-1.60000000e-01,  0.00000000e+00)
(-1.20000000e-01,  0.00000000e+00)
(-8.00000000e-02,  0.00000000e+00)
(-4.00000000e-02,  0.00000000e+00)
( 0.00000000e+00,  0.00000000e+00)
( 4.00000000e-02,  0.00000000e+00)
( 8.00000000e-02,  0.00000000e+00)
( 1.20000000e-01,  0.00000000e+00)
( 1.60000000e-01,  0.00000000e+00)
( 2.00000000e-01,  0.00000000e+00)
( 2.40000000e-01,  0.00000000e+00)
( 2.80000000e-01,  0.00000000e+00)
( 3.20000000e-01,  0.00000000e+00)
( 3.60000000e-01,  0.00000000e+00)
( 4.00000000e-01,  0.00000000e+00)
( 4.40000000e-01,  0.00000000e+00)
( 4.80000000e-01,  0.00000000e+00)
( 5.20000000e-01,  0.00000000e+00)
( 5.60000000e-01,  0.00000000e+00)
( 6.00000000e-01,  0.00000000e+00)
( 6.40000000e-01,  0.00000000e+00)
( 6.80000000e-01,  0.00000000e+00)
( 7.20000000e-01,  0.00000000e+00)
( 7.60000000e-01,  0.00000000e+00)
( 8.00000000e-01,  0.00000000e+00)
( 8.40000000e-01,  0.00000000e+00)
( 8.80000000e-01,  0.00000000e+00)
( 9.20000000e-01,  0.00000000e+00)
( 9.60000000e-01,  0.00000000e+00)
(-4.76204097e-01,  9.51071830e-02)
(-4.26393645e-01,  2.94333072e-01)
(-3.76010317e-01,  4.95698499e-01)
(-3.26509758e-01,  6.94370258e-01)
(-2.76166065e-01,  8.95303253e-01)
(-9.58291625e-01,  1.00000000e+00)
(-9.16583250e-01,  1.00000000e+00)
(-8.74874875e-01,  1.00000000e+00)
(-8.33166500e-01,  1.00000000e+00)
(-7.91458125e-01,  1.00000000e+00)
(-7.49749750e-01,  1.00000000e+00)
(-7.08041375e-01,  1.00000000e+00)
(-6.66333000e-01,  1.00000000e+00)
(-6.24624625e-01,  1.00000000e+00)
(-5.82916250e-01,  1.00000000e+00)
(-5.41207875e-01,  1.00000000e+00)
(-4.99499499e-01,  1.00000000e+00)
(-4.57791124e-01,  1.00000000e+00)
(-4.16082749e-01,  1.00000000e+00)
(-3.74374374e-01,  1.00000000e+00)
(-3.32665999e-01,  1.00000000e+00)
(-2.90957624e-01,  1.00000000e+00)
(-2.49249249e-01,  1.00000000e+00)
(-2.10210210e-01,  1.00000000e+00)
(-1.71171171e-01,  1.00000000e+00)
(-1.32132132e-01,  1.00000000e+00)
(-9.30930931e-02,  1.00000000e+00)
(-5.40540541e-02,  1.00000000e+00)
(-1.50150150e-02,  1.00000000e+00)
( 2.40240240e-02,  1.00000000e+00)
( 6.30630631e-02,  1.00000000e+00)
( 1.02102102e-01,  1.00000000e+00)
( 1.41141141e-01,  1.00000000e+00)
( 1.80180180e-01,  1.00000000e+00)
( 2.19219219e-01,  1.00000000e+00)
( 2.58258258e-01,  1.00000000e+00)
( 2.97297297e-01,  1.00000000e+00)
( 3.36336336e-01,  1.00000000e+00)
( 3.75375375e-01,  1.00000000e+00)
( 4.14414414e-01,  1.00000000e+00)
( 4.53453453e-01,  1.00000000e+00)
( 4.92492492e-01,  1.00000000e+00)
( 5.31531532e-01,  1.00000000e+00)
( 5.70570571e-01,  1.00000000e+00)
( 6.09609610e-01,  1.00000000e+00)
( 6.48648649e-01,  1.00000000e+00)
( 6.87687688e-01,  1.00000000e+00)
( 7.26726727e-01,  1.00000000e+00)
( 7.65765766e-01,  1.00000000e+00)
( 8.04804805e-01,  1.00000000e+00)
( 8.43843844e-01,  1.00000000e+00)
( 8.82882883e-01,  1.00000000e+00)
( 9.21921922e-01,  1.00000000e+00)
( 9.60960961e-01,  1.00000000e+00)};\label{line:burgers_xc}

\end{axis}
\end{tikzpicture}

%% file: figs/cylinder/control_points/cy_xc.tikz
\begin{tikzpicture}
\begin{axis}[
axis equal image,
axis lines=none,
width=1.\textwidth,
ymax=3.0,
xmax=4.202421881177826,
xmin=-4.202421881177824,
ymin=0.45399049973954675]
\addplot []
graphics [xmin=-4.202421881177824,xmax=4.202421881177826,ymin=0.45399049973954675,ymax=3.0] { 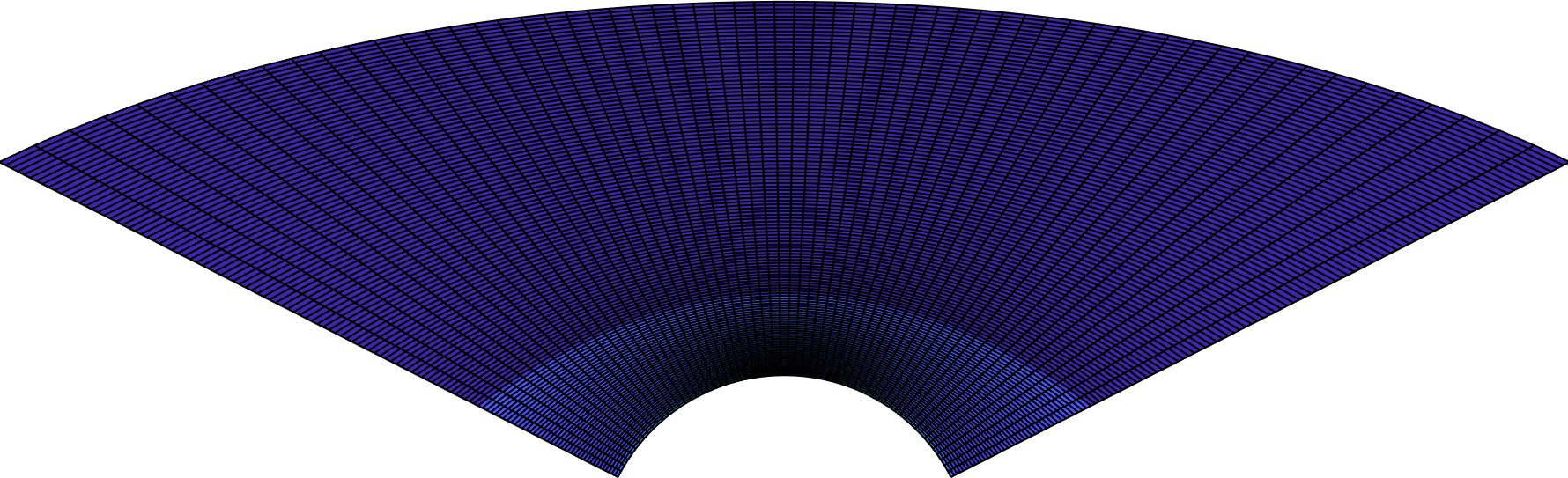};

\addplot [red, only marks, mark options={solid, thin}, mark=square*, mark size=0.5, mark repeat=0]
coordinates {
( 8.91006524e-01,  4.53990500e-01)
(-8.91006524e-01,  4.53990500e-01)
( 4.20242188e+00,  2.14124090e+00)
(-4.20242188e+00,  2.14124090e+00)
( 8.70183755e-01,  4.92727342e-01)
( 8.47677936e-01,  5.30511184e-01)
( 8.23532598e-01,  5.67268949e-01)
( 7.97794440e-01,  6.02929542e-01)
( 7.70513243e-01,  6.37423990e-01)
( 7.41741773e-01,  6.70685577e-01)
( 7.11535677e-01,  7.02649970e-01)
( 6.79953379e-01,  7.33255346e-01)
( 6.47055962e-01,  7.62442511e-01)
( 6.12907054e-01,  7.90155012e-01)
( 5.77572703e-01,  8.16339251e-01)
( 5.41121252e-01,  8.40944582e-01)
( 5.03623202e-01,  8.63923417e-01)
( 4.65151078e-01,  8.85231311e-01)
( 4.25779292e-01,  9.04827052e-01)
( 3.85583992e-01,  9.22672740e-01)
( 3.44642923e-01,  9.38733858e-01)
( 3.03035270e-01,  9.52979342e-01)
( 2.60841506e-01,  9.65381639e-01)
( 2.18143241e-01,  9.75916762e-01)
( 1.75023059e-01,  9.84564335e-01)
( 1.31564359e-01,  9.91307631e-01)
( 8.78511966e-02,  9.96133609e-01)
( 4.39681183e-02,  9.99032935e-01)
( 2.83276945e-16,  1.00000000e+00)
(-4.39681183e-02,  9.99032935e-01)
(-8.78511966e-02,  9.96133609e-01)
(-1.31564359e-01,  9.91307631e-01)
(-1.75023059e-01,  9.84564335e-01)
(-2.18143241e-01,  9.75916762e-01)
(-2.60841506e-01,  9.65381639e-01)
(-3.03035270e-01,  9.52979342e-01)
(-3.44642923e-01,  9.38733858e-01)
(-3.85583992e-01,  9.22672740e-01)
(-4.25779292e-01,  9.04827052e-01)
(-4.65151078e-01,  8.85231311e-01)
(-5.03623202e-01,  8.63923417e-01)
(-5.41121252e-01,  8.40944582e-01)
(-5.77572703e-01,  8.16339251e-01)
(-6.12907054e-01,  7.90155012e-01)
(-6.47055962e-01,  7.62442511e-01)
(-6.79953379e-01,  7.33255346e-01)
(-7.11535677e-01,  7.02649970e-01)
(-7.41741773e-01,  6.70685577e-01)
(-7.70513243e-01,  6.37423990e-01)
(-7.97794440e-01,  6.02929542e-01)
(-8.23532598e-01,  5.67268949e-01)
(-8.47677936e-01,  5.30511184e-01)
(-8.70183755e-01,  4.92727342e-01)
( 3.97143577e+00,  2.24876065e+00)
( 3.74509850e+00,  2.34383432e+00)
( 3.52458945e+00,  2.42782151e+00)
( 3.31062264e+00,  2.50198810e+00)
( 3.10355717e+00,  2.56748578e+00)
( 2.90348909e+00,  2.62534527e+00)
( 2.71032480e+00,  2.67647807e+00)
( 2.52383812e+00,  2.72168336e+00)
( 2.34371361e+00,  2.76165741e+00)
( 2.16957878e+00,  2.79700411e+00)
( 2.00102758e+00,  2.82824542e+00)
( 1.83763715e+00,  2.85583130e+00)
( 1.67897954e+00,  2.88014877e+00)
( 1.52462953e+00,  2.90152998e+00)
( 1.37416972e+00,  2.92025930e+00)
( 1.22719355e+00,  2.93657947e+00)
( 1.08330680e+00,  2.95069680e+00)
( 9.42128023e-01,  2.96278563e+00)
( 8.03288259e-01,  2.97299209e+00)
( 6.66430169e-01,  2.98143719e+00)
( 5.31206821e-01,  2.98821934e+00)
( 3.97280224e-01,  2.99341646e+00)
( 2.64319691e-01,  2.99708755e+00)
( 1.32000083e-01,  2.99927391e+00)
( 8.49830835e-16,  3.00000000e+00)
(-1.32000083e-01,  2.99927391e+00)
(-2.64319691e-01,  2.99708755e+00)
(-3.97280224e-01,  2.99341646e+00)
(-5.31206821e-01,  2.98821934e+00)
(-6.66430169e-01,  2.98143719e+00)
(-8.03288259e-01,  2.97299209e+00)
(-9.42128023e-01,  2.96278563e+00)
(-1.08330680e+00,  2.95069680e+00)
(-1.22719355e+00,  2.93657947e+00)
(-1.37416972e+00,  2.92025930e+00)
(-1.52462953e+00,  2.90152998e+00)
(-1.67897954e+00,  2.88014877e+00)
(-1.83763715e+00,  2.85583130e+00)
(-2.00102758e+00,  2.82824542e+00)
(-2.16957878e+00,  2.79700411e+00)
(-2.34371361e+00,  2.76165741e+00)
(-2.52383812e+00,  2.72168336e+00)
(-2.71032480e+00,  2.67647807e+00)
(-2.90348909e+00,  2.62534527e+00)
(-3.10355717e+00,  2.56748578e+00)
(-3.31062264e+00,  2.50198810e+00)
(-3.52458945e+00,  2.42782151e+00)
(-3.74509850e+00,  2.34383432e+00)
(-3.97143577e+00,  2.24876065e+00)
(-1.49902378e+00,  9.27241390e-01)
(-1.24751572e+00,  1.08293722e+00)
(-9.67064629e-01,  1.22264208e+00)
(-6.76022159e-01,  1.32993358e+00)
(-3.62347088e-01,  1.40139872e+00)
(-3.62150912e-02,  1.43011332e+00)
( 2.85446291e-01,  1.41566701e+00)
( 6.09615681e-01,  1.35091310e+00)
( 9.29542048e-01,  1.23770473e+00)
( 1.22197445e+00,  1.09574603e+00)
( 1.49063100e+00,  9.36745443e-01)
( 9.65513370e-01,  4.91953634e-01)
( 1.04002022e+00,  5.29916768e-01)
( 1.11452706e+00,  5.67879902e-01)
( 1.18903391e+00,  6.05843036e-01)
( 1.26354075e+00,  6.43806170e-01)
( 1.33804760e+00,  6.81769304e-01)
( 1.41255444e+00,  7.19732437e-01)
( 1.48706129e+00,  7.57695571e-01)
( 1.56156813e+00,  7.95658705e-01)
( 1.63607498e+00,  8.33621839e-01)
( 1.70023365e+00,  8.66312316e-01)
( 1.76439232e+00,  8.99002792e-01)
( 1.82855100e+00,  9.31693269e-01)
( 1.89270967e+00,  9.64383745e-01)
( 1.95686834e+00,  9.97074222e-01)
( 2.02102701e+00,  1.02976470e+00)
( 2.08518569e+00,  1.06245517e+00)
( 2.14934436e+00,  1.09514565e+00)
( 2.21350303e+00,  1.12783613e+00)
( 2.27766170e+00,  1.16052660e+00)
( 2.34182038e+00,  1.19321708e+00)
( 2.40597905e+00,  1.22590756e+00)
( 2.47013772e+00,  1.25859803e+00)
( 2.53429640e+00,  1.29128851e+00)
( 2.59845507e+00,  1.32397899e+00)
( 2.66261374e+00,  1.35666946e+00)
( 2.72677241e+00,  1.38935994e+00)
( 2.79093109e+00,  1.42205042e+00)
( 2.85508976e+00,  1.45474089e+00)
( 2.91924843e+00,  1.48743137e+00)
( 2.98340710e+00,  1.52012185e+00)
( 3.04756578e+00,  1.55281232e+00)
( 3.11172445e+00,  1.58550280e+00)
( 3.17588312e+00,  1.61819327e+00)
( 3.24004179e+00,  1.65088375e+00)
( 3.30420047e+00,  1.68357423e+00)
( 3.36835914e+00,  1.71626470e+00)
( 3.43251781e+00,  1.74895518e+00)
( 3.49667648e+00,  1.78164566e+00)
( 3.56083516e+00,  1.81433613e+00)
( 3.62499383e+00,  1.84702661e+00)
( 3.68915250e+00,  1.87971709e+00)
( 3.75331117e+00,  1.91240756e+00)
( 3.81746985e+00,  1.94509804e+00)
( 3.88162852e+00,  1.97778852e+00)
( 3.94578719e+00,  2.01047899e+00)
( 4.00994586e+00,  2.04316947e+00)
( 4.07410454e+00,  2.07585995e+00)
( 4.13826321e+00,  2.10855042e+00)
(-9.65513370e-01,  4.91953634e-01)
(-1.04002022e+00,  5.29916768e-01)
(-1.11452706e+00,  5.67879902e-01)
(-1.18903391e+00,  6.05843036e-01)
(-1.26354075e+00,  6.43806170e-01)
(-1.33804760e+00,  6.81769304e-01)
(-1.41255444e+00,  7.19732437e-01)
(-1.48706129e+00,  7.57695571e-01)
(-1.56156813e+00,  7.95658705e-01)
(-1.63607498e+00,  8.33621839e-01)
(-1.70023365e+00,  8.66312316e-01)
(-1.76439232e+00,  8.99002792e-01)
(-1.82855100e+00,  9.31693269e-01)
(-1.89270967e+00,  9.64383745e-01)
(-1.95686834e+00,  9.97074222e-01)
(-2.02102701e+00,  1.02976470e+00)
(-2.08518569e+00,  1.06245517e+00)
(-2.14934436e+00,  1.09514565e+00)
(-2.21350303e+00,  1.12783613e+00)
(-2.27766170e+00,  1.16052660e+00)
(-2.34182038e+00,  1.19321708e+00)
(-2.40597905e+00,  1.22590756e+00)
(-2.47013772e+00,  1.25859803e+00)
(-2.53429640e+00,  1.29128851e+00)
(-2.59845507e+00,  1.32397899e+00)
(-2.66261374e+00,  1.35666946e+00)
(-2.72677241e+00,  1.38935994e+00)
(-2.79093109e+00,  1.42205042e+00)
(-2.85508976e+00,  1.45474089e+00)
(-2.91924843e+00,  1.48743137e+00)
(-2.98340710e+00,  1.52012185e+00)
(-3.04756578e+00,  1.55281232e+00)
(-3.11172445e+00,  1.58550280e+00)
(-3.17588312e+00,  1.61819327e+00)
(-3.24004179e+00,  1.65088375e+00)
(-3.30420047e+00,  1.68357423e+00)
(-3.36835914e+00,  1.71626470e+00)
(-3.43251781e+00,  1.74895518e+00)
(-3.49667648e+00,  1.78164566e+00)
(-3.56083516e+00,  1.81433613e+00)
(-3.62499383e+00,  1.84702661e+00)
(-3.68915250e+00,  1.87971709e+00)
(-3.75331117e+00,  1.91240756e+00)
(-3.81746985e+00,  1.94509804e+00)
(-3.88162852e+00,  1.97778852e+00)
(-3.94578719e+00,  2.01047899e+00)
(-4.00994586e+00,  2.04316947e+00)
(-4.07410454e+00,  2.07585995e+00)
(-4.13826321e+00,  2.10855042e+00)};\label{line:cy_xc}

\end{axis}
\end{tikzpicture}